\documentclass[final,3p,times]{elsarticle}

\usepackage{amssymb}
\usepackage{amsthm}
\usepackage{amsmath}
\usepackage{color}
\usepackage{subfigure}
\usepackage[percent]{overpic}

\usepackage{algorithm}
\usepackage{algorithmic}
\usepackage{vruler}

\allowdisplaybreaks

\newtheorem{theorem}{Theorem}
\newtheorem{lemma}[theorem]{Lemma}

\newtheorem{remark}[theorem]{Remark}

\journal{Computer Methods in Applied Mechanics and Engineering}
\usepackage{hyperref}
\begin{document}
\begin{frontmatter}



\title{Improved stabilization technique for frictional contact problems solved with $hp$-BEM\tnoteref{t1}}
\tnotetext[t1]{Dedicated to Professor Ernst P.~Stephan on the occasion of his 70th birthday}


\author[lab2]{Lothar~Banz\corref{cor1}}
\ead{lothar.banz@sbg.ac.at}

\author[lab2]{Gregor~Milicic}
\ead{gregor.milicic@sbg.ac.at}

\author[lab1]{Nina~Ovcharova}
\ead{nina.ovcharova@unibw.de}

\address[lab2]{Department of Mathematics, University of Salzburg, Hellbrunner Stra{\ss}e 34, 5020 Salzburg, Austria}

\address[lab1]{Universit\"at der Bundeswehr M\"unchen, D-85577 Neubiberg/Munich, Germany}
\cortext[cor1]{Corresponding author}

\begin{abstract}
We improve the residual based stabilization technique for Signorini contact problems with Tresca friction in linear elasticity solved with $hp$-mixed BEM which has been recently analyzed by Banz et al.~in Numer.~Math.~135 (2017) pp.~217--263. The stabilization allows us to circumvent the discrete inf-sup condition and thus the primal and dual sets can be discretized independently. Compared to the above mentioned paper we are able to remove the dependency of the scaling parameter on the unknown Sobolev regularity of the exact solution and can thus also improve the convergence rate in the a priori error estimate. The second improvement is a rigorous a priori and a posteriori error analysis when the boundary integral operators in the stabilization term are approximated. The latter is of fundamental importance to keep the computational time small.
We present numerical results in two and three dimensions to underline our theoretical findings, show the superiority of the $hp$-adaptive stabilized mixed scheme and the effect induced by approximating the stabilization term. Moreover, we show the applicability of the proposed method to the Coulomb frictional case for which we extend the a posteriori error analysis.
\end{abstract}


\end{frontmatter}


\section{Introduction}
\label{sec:Introduction}

Simulations of frictional contact problems have a vast amount of practical applications and are used daily in the design of new products, e.g.~in tire engineering. In these applications, the contact stresses are often the main quantities of interest. Mixed methods as analyzed in e.g.~\cite{Kikuchi1988Contact,Hlavacek1988,Hild2002,boffi2013Mixed} among many others have the advantage that the product engineer has direct access to the contact stresses and can control explicitly the approximation quality of these. However, the drawback of mixed methods is that also at the discrete level the (discrete) inf-sup condition must be satisfied to guarantee uniqueness of the discrete solution. For a priori error analysis even the asymptotic behavior of the discrete inf-sup constant must be known. The latter is a measure for how well the discrete primal space and the discrete dual space are balanced out, in particular, if the discrete dual space is sufficiently small compared to the discrete primal space. In contact problems, that can be achieved by coarsening sufficiently the mesh and/or decreasing sufficiently the polynomial degree for the dual variable as has been shown in \cite{Schroeder09SIN}. However, what sufficiently means remains open. It is well known that the desirable case of using the same mesh for both variables and reducing the polynomial degree by at most one for the dual variable is insufficient.\\
In finite element methods a stabilization technique based on the square of the discrete residual to circumvent the discrete inf-sup condition goes back to \cite{barbosahughes1991,barbosahughes1992}. That strategy has been used recently in \cite{hildrenard2010} to stabilize a lowest order $h$-FE discretization of contact problems and in \cite{biermann2013stabilization} for a lowest order $h$-FE discretization of an obstacle problem. Only very recently that strategy has been generalized to $hp$-discretizations and to boundary element methods in \cite{Banz2015Stab}. The major difficulty in stabilizing the BEM is that even though the discrete residual lies in $L^2$, going to the limit of the mapping properties of the boundary integral operators \cite{costabel1988boundary}, the natural norm to estimate it is the dual norm $H^{-1/2}$. That becomes even more severe when the Poincar\'e-Steklov operator $S: H^{1/2} \rightarrow H^{-1/2}$ is used to describe the residual, as done in \cite{Banz2015Stab}, since on the discrete level the operator $S$ must be discretized as well, i.e.~$S$ becomes $S_{hp}$. Consequently, the consistency error $E_{hp}=S-S_{hp}$ must be measured in the stronger $L^2$-norm due to the stabilization term. In \cite{Banz2015Stab} to guarantee convergence, the authors increase the $h$ and $p$-dependency of the stabilization scaling parameter compared to what is necessary by the polynomial inverse estimate. On the one hand, that enables them to proof a priori error convergence rates, on the other hand these are even more suboptimal than usual in contact problems, and, moreover, the optimal scaling of the stabilization term becomes dependent on the unknown Sobolev regularity of the exact solution.\\
Independently of whether the mixed formulation is stabilized or not, and if a mixed formulation as in e.g.~\cite{hueber2005mixed,wohlmuth2011variationally} is used at all or rather a variational inequality of first or of second kind as in e.g.~\cite{gwinner1993boundary,gwinner2009p} is used, the solution of the contact problem is typically of reduced regularity at the interfaces from contact to non-contact and from stick to slip. The locations of these interfaces are however not known a priori, and thus one must resort to automatic mesh refinements, $h$-adaptivity, or automatic polynomial degree adaptation, $p$-adaptivity, or best to $hp$-adaptivity which is a combination of both, to obtain improved or even optimal convergence rates. The necessary a posteriori error estimate for $h$-adaptivity has been derived in e.g.~\cite{veeser2001efficient,bartels2004averaging,hild2009residual} for contact and obstacle problems, and extended to the $hp$-setting in e.g.~\cite{maischak2005adaptive,maischak2007adaptive,schroder2012posteriori,banz2013posteriori,banz2013hp,Banz2014BEM,banz2015biorthogonal} and extended to the even more complicated class of hemivariational inequalities in \cite{NinaLothar}.\\
In this paper we use a weak formulation not based on the Poincar\'e-Steklov operator but based on the entire Calder\'on projector, which leads to a slightly different stabilization term compared to \cite{Banz2015Stab}. In particular, the approximation of the inverse of the single layer potential does not appear in the stabilization term. That allows us to proof convergence rates in Theorem~\ref{thm:AprioriConvergencerate} where the scaling of the stabilization terms does not depend on the Sobolev regularity of the exact solution contrary to \cite{Banz2015Stab}. As our discrete solution is very close to the discrete solution of \cite{Banz2015Stab}, c.f.~Lemma~\ref{lem:dist_2_NuMa_Paper}, we can improve the a priori error estimate in \cite[Thm.~16]{Banz2015Stab} and remove the dependency of the scaling of the stabilization terms on the Sobolev regularity of the exact solution, c.f.~Theorem~\ref{thm:improvedConvergenceRates}.\\
In \cite[Sec.~6]{Banz2015Stab} an implementation of the highly problematic term $\langle W u_{hp}, W v_{hp}\rangle$ which appears in the stabilization term is discussed. Their numerical results in \cite[Sec.~8.3]{Banz2015Stab} indicate that an approximation of the stabilization term works in practice very well, but a rigorous error analysis for that is missing. We fill that gap by allowing the boundary integral operators in the bilinear form to differ from those used in the stabilization term. Naturally, the arguments for the proofs are quite close to those in \cite{Banz2015Stab} and it would have been sufficient to only point out the required modifications in the proofs, but we state the entire proofs here to be self-contained.\\
The rest of the paper is structured as follows: In Section~\ref{sec:Amixedboundaryintegralformulation} we introduce the boundary integral formulation for a Tresca frictional contact problem based on the entire Calder\'on projector. Section~\ref{sec:stabalizedDiscretization} is devoted to the $hp$-discretization and stabilization of the mixed problem. In particular, the important result of existence, uniqueness and stability of a discrete solution is shown in Theorem~\ref{thm:discreteExistenceUniquenessStability} without using the discrete inf-sup condition. An a priori error analysis is carried out in Section~\ref{sec:Apriorierrorestimate} with the two major results Theorem~\ref{thm:AprioriConvergencerate}, a priori convergence rates for our problem, and Theorem~\ref{thm:improvedConvergenceRates} improved convergence rate for the formulation in \cite{Banz2015Stab}. In Section~\ref{sec:Aposteriorierrorestimate} the reliable a posteriori error estimate of Theorem~\ref{thm:aposterioriErrorEst} is proven. Sections~\ref{sec:CoulombFriction} and \ref{sec:approximation} deal with the extension to Coulomb friction and a practical approximation of the boundary integral operators in the stabilization term. Our theoretical results are underlined by the numerical results, presented in Section~\ref{sec:Numerics}, for discretizations which do not satisfy the discrete inf-sup condition.\\

\textbf{Notation:} We denote by $X^\prime$ the dual space or adjoint operator to $X$ and by $\langle \cdot,\cdot \rangle_{\Gamma}$ the duality pairing with integration domain $\Gamma$. Moreover, we use the notation $\| v \|_{\widetilde{H}^{1/2}(\Gamma_\Sigma)} = \| v_0\|_{H^{1/2}(\Gamma)}$ with $v_0$ the extension of $v$ by zero from $\Gamma_\Sigma \subset \Gamma$ to the whole $\Gamma$. If the meaning is clear we omit the index zero which indicates the extension. We denote the Euclidean norm by $| \cdot |$.
If the number of indices become too many some of them move to the exponent, e.g.~$v_{hp,n}= v_{n}^{hp} = v_{hp} \cdot n$. We use generic constants $C$ which take different values at different places.

\section{A mixed boundary integral formulation for Tresca frictional contact problems}
\label{sec:Amixedboundaryintegralformulation}

Let $\Omega \subset \mathbb{R}^d$ ($d=2,3$) be a bounded polygonal domain with boundary $\Gamma$ and outward unit normal $n$. We assume that $\Gamma$ is already sufficiently scaled such that $\operatorname{cap}(\Gamma)<1$ if $d=2$. Furthermore, let $\overline{\Gamma}=\overline{\Gamma}_D \cup \overline{\Gamma}_N \cup \overline{\Gamma}_C$ be decomposed into non-overlapping, Dirichlet, Neumann and contact boundary parts, and denote by {$\overline{\Gamma}_\Sigma:=\overline{\Gamma}_N \cup \overline{\Gamma}_C$} the union of the latter two.
For the ease of presentation we assume $\overline{\Gamma}_D\cap \overline{\Gamma}_C = \emptyset$. For given gap function $0 \leq g \in H^{1/2}(\Gamma_C)$, friction threshold $0< \mathcal{F} \in L^\infty(\Gamma_C)$, Neumann data $f\in \widetilde{H}^{-1/2}(\Gamma_N)$ and elasticity tensor $\mathcal{C}$ the considered Tresca frictional contact problem is to find a weak solution $u$ to
\begin{subequations}  \label{prob:strong_formulation}
\begin{alignat}{2}
-\operatorname{div} \sigma(u)&=0 & \quad & \text{in } \label{eq1} \Omega\\
  \sigma(u)&=\mathcal{C}:\epsilon(u) & \quad & \text{in } \Omega \label{eq:prog:materialLaw} \\
	u&=0 & \quad & \text{on } \Gamma_D \\
  \sigma(u)n&={f} & \quad & \text{on } \Gamma_N \\
  \sigma_n\leq 0,\ u_n\leq g,\ \sigma_n(u_n-g)&=0 & \quad & \text{on } \Gamma_C\\
  | \sigma_t |  \leq \mathcal{F},\ \sigma_tu_t+\mathcal{F}| u_t | &=0 & \quad & \text{on } \Gamma_C. \label{eq:prog:FrictionLaw}
\end{alignat}
\end{subequations}
As usual, $\sigma_n$, $u_n\in \mathbb{R}$, $\sigma_t$, $u_t \in \mathbb{R}^{d-1}$ are the normal, tangential components of $\sigma(u)n$, $u$, respectively and \eqref{eq:prog:materialLaw} describes Hooke's law with the linearized strain tensor $\epsilon(u)=\frac{1}{2} ( \nabla u + \nabla u^\top )$. Equation \eqref{eq:prog:FrictionLaw} may be written equivalently in the form
\begin{align} \label{eq:frictionCondition}
 | \sigma_t |  \leq \mathcal{F}, \qquad
 | \sigma_t |  < \mathcal{F} \Rightarrow u_t=0, \qquad
 | \sigma_t |  = \mathcal{F} \Rightarrow \exists \,  \alpha \geq 0 : u_t=-\alpha \sigma_t.
\end{align}
For the solution $u(x)$ of \eqref{eq1}-\eqref{eq:prog:materialLaw} with $x \in \Omega\backslash \Gamma$ 
we have the following representation formula, also known as Somigliana's identity, see e.g.~\cite{Kleiber}:
\begin{align} \label{eq2}
u(x)=\int_{\Gamma} G(x,y) T_y u(y)\, ds_y -\int_\Gamma T_y G(x,y)u(y) \,ds_y , 
\end{align}
where $G(x,y)$ is the fundamental solution of the Navier-Lam\'{e} equation defined by
\begin{align*}
G(x,y)=\begin{cases}\frac{\lambda +3\mu}{4\pi \mu(\lambda+2\mu)}\left(\log |x-y|{\rm I} + \frac{\lambda + \mu}{\lambda +3\mu}\frac{(x-y)(x-y)^\top}{|x-y|^2}\right), \quad &\text{ if d=2}\\[0.2cm]
\frac{\lambda +3\mu}{8\pi \mu(\lambda+2\mu)}\left( |x-y|^{-1}{\rm I} + \frac{\lambda + \mu}{\lambda +3\mu}\frac{(x-y)(x-y)^\top}{|x-y|^3}\right), \quad &\text{ if d=3}
\end{cases}
\end{align*}
with the Lam\'{e} constants $\lambda, \mu >0$ depending on the material parameters, i.e.~the modulus of elasticity~$E$ and the Poisson's ratio~$\nu$:
$$
\lambda = \frac{E\nu}{1-\nu^2} \,, \qquad \mu= \frac{E}{1+\nu}\, . 
$$
 Here, $T_y$  stands for the traction operator with respect to $y$ defined by $T_y (u) := \sigma(u(y))\cdot n_y$, where $n_y$ is  the unit outer normal vector  at $y\in \Gamma$. 
Letting $\Omega\backslash \Gamma \ni x \to \Gamma$ in (\ref{eq2}), we obtain the well-known  system of boundary integral equations 
\[
 \left(\begin{array}{c}
          u\\[0.2cm]
         T_x u
         \end{array}
 \right)
=  \left( \begin{array} {cc}
           \frac{1}{2}I-K & V\\[0.2cm]
           W &      \frac{1}{2}I+K'  
           \end{array}
   \right) 
\left( \begin{array}{c}
u\\[0.2cm]
T_x u  
\end{array}
 \right)
\]       
with the single layer potential $V$, the double layer potential
$K$, its formal adjoint $K'$, and the hypersingular integral operator $W$ defined for $x \in \Gamma$ as follows:
  \begin{alignat*}{2}
   ( V \phi ) (x) & :=  \int\limits_\Gamma G(x,y) \phi (y) \, ds_y, & \qquad
   (K v  ) (x) &:=
    \int_{\Gamma} \left(T_y G(x,y)\right)^\top  v (y) \, ds_y \\ 
   ( K' \phi ) (x) & := T_x \int\limits_\Gamma G 
    (x,y)  \phi (y) \, ds_y, & \qquad
    (W v  ) (x) &:= - T_x   (K v ) (x).
  \end{alignat*}
By introducing the additional unknown $\lambda = -\sigma(u) \cdot n$ on $\Gamma_C$ and using the entire Calder\'on projector we straightforwardly obtain the following weak formulation defined on the boundary only. That is: Find a triple $(u,\phi,\lambda) \in \widetilde{H}^{1/2}(\Gamma_\Sigma) \times H^{-1/2}(\Gamma) \times M^+(\mathcal{F}) $ such that \pagebreak[1]
\begin{subequations} \label{eq:MixedProblem}
\begin{alignat}{6}
 &\langle Wu,v \rangle_{\Gamma} &\,+\,& \langle (K+1/2)^\prime \phi,v \rangle_{\Gamma} + \langle \lambda,v \rangle_{\Gamma_C} & = & \; \langle {f},v \rangle_{\Gamma_N} &\quad &\forall v\in \widetilde{H}^{1/2}(\Gamma_\Sigma), \label{eq:WeakMixedVarEq1} \\
-&\langle (K+ 1/2) u,\psi \rangle_{\Gamma} &\, +\,&\langle V \phi,\psi \rangle_{\Gamma}  & = & \; 0 &\quad &\forall \psi \in H^{-1/2}(\Gamma), \label{eq:WeakMixedVarEq2} \\ 
& \langle u, \mu -\lambda \rangle_{\Gamma_C} & & & \leq & \;  \langle g,\mu_n-\lambda_n \rangle_{\Gamma_C} &\quad &\forall \mu \in M^+(\mathcal{F})  \label{eq:ContContactConstraints}
\end{alignat}
\end{subequations}
with the set of admissible Lagrange multipliers 
\begin{align}
M^+(\mathcal{F}):=\left\lbrace \mu \in \widetilde{H}^{-1/2}(\Gamma_C): \langle \mu,v\rangle_{\Gamma_C} \leq \langle \mathcal{F},|v_t| \rangle_{\Gamma_C} \  \forall v \in \widetilde{H}^{1/2}(\Gamma_\Sigma), v_n \leq 0 \right\rbrace .
\end{align}
By condensing out $\phi$, we obtain the weak formulation considered in e.g.~\cite{Banz2014BEM,Banz2015Stab} with the symmetric Poincar\'e-Steklov operator $S=W+(K+\frac{1}{2})^\prime V^{-1}(K+\frac{1}{2})$. That is: Find $(u,\lambda) \in \widetilde{H}^{1/2}(\Gamma_\Sigma) \times M^+(\mathcal{F})$ such that
\begin{subequations} \label{eq:SteklovFormulation}
\begin{alignat}{2}
 \langle Su,v \rangle_{\Gamma} + \langle \lambda,v \rangle_{\Gamma_C} &= \langle {f},v \rangle_{\Gamma_N} &\quad &\forall v\in \widetilde{H}^{1/2}(\Gamma_\Sigma),   \\
\langle u, \mu -\lambda \rangle_{\Gamma_C} & \leq \langle g,\mu_n-\lambda_n \rangle_{\Gamma_C} &\quad &\forall \mu \in M^+(\mathcal{F}) . 
\end{alignat}
\end{subequations}
Consequently, all the existence, uniqueness, stability and Lipschitz-dependency on the data results carry over one by one to our situation here. 
\begin{theorem}
\begin{enumerate}
 \item The continuous inf-sup condition is satisfied, i.e.~there exists a constant $\beta>0$  such that
 	\begin{align}\label{eq:lame:infsup}
	  \beta \| \mu \|_{\widetilde{H}^{-1/2}(\Gamma_C)} \leq \sup_{v \in \widetilde{H}^{1/2}(\Gamma_\Sigma) \setminus \{0\}} \frac{\langle \mu, v \rangle_{\Gamma_C} }{ \|v\|_{\widetilde{H}^{1/2}(\Gamma_\Sigma)} } \qquad \forall \mu \in \widetilde{H}^{-1/2}(\Gamma_C) \ .
	\end{align}
 \item There exists a unique solution $(u,\phi,\lambda) \in \widetilde{H}^{1/2}(\Gamma_\Sigma) \times H^{-1/2}(\Gamma) \times M^+(\mathcal{F}) $ to \eqref{eq:MixedProblem}.
 \item There exists a constant $C>0$ such that for the solution of \eqref{eq:MixedProblem} it holds
 \begin{align*}
   \| u \|_{\widetilde{H}^{1/2}(\Gamma_\Sigma)} + \| \phi \|_{H^{-1/2}(\Gamma)} + \| \lambda \|_{\widetilde{H}^{-1/2}(\Gamma_C)} \leq C \left( \| f\|_{\widetilde{H}^{-1/2}(\Gamma_N)} + \|g\|_{H^{1/2}(\Gamma_C)} \right).
 \end{align*}
\item The solution of \eqref{eq:MixedProblem} depends Lipschitz continuous in the data $f$ and $g$, i.e.~there exists a constant $C>0$ such that
 \begin{align*}
   \| u_1-u_2 \|_{\widetilde{H}^{1/2}(\Gamma_\Sigma)} + \| \phi_1 - \phi_2 \|_{H^{-1/2}(\Gamma)} + \| \lambda_1-\lambda_2 \|_{\widetilde{H}^{-1/2}(\Gamma_C)} \leq C \left( \| f_1-f_2\|_{\widetilde{H}^{-1/2}(\Gamma_N)} + \|g_1-g_2\|_{H^{1/2}(\Gamma_C)} \right).
 \end{align*}
 with $(u_i,\phi_i,\lambda_i)$ the solution to the data $(f_i,g_i)$.
\end{enumerate}
\end{theorem}
The latter weak formulation \eqref{eq:SteklovFormulation} has for our kind of residual based stabilization some theoretical drawbacks, as the stabilization factor for the formulation \eqref{eq:SteklovFormulation}, see \cite{Banz2015Stab}, depends on the unknown regularity of the solution $u$ and the proven theoretical convergence rate is far from being optimal, see~\cite[Thm.~16]{Banz2015Stab}. The reason is that in \cite{Banz2015Stab} the approximation error of $V^{-1}$ needs to be measured in the $L^2(\Gamma)$-norm instead of the weaker $H^{-1/2}(\Gamma)$-norm due to the stabilization term. That problem can be avoided when stabilizing \eqref{eq:MixedProblem} instead of \eqref{eq:SteklovFormulation}.
By introducing the coercive but non-symmetric bilinear form 
\begin{align*}
 B(u,\phi;v,\psi)= \langle Wu,v \rangle_{\Gamma} + \langle (K+1/2)^\prime \phi,v \rangle_{\Gamma} + \langle V \phi,\psi \rangle_{\Gamma} - \langle (K+ 1/2) u,\psi \rangle_{\Gamma}
\end{align*}
we may write \eqref{eq:MixedProblem} more compactly as finding $(u,\phi,\lambda) \in \widetilde{H}^{1/2}(\Gamma_\Sigma) \times H^{-1/2}(\Gamma) \times M^+(\mathcal{F}) $ such that
\begin{subequations} \label{eq:compactWeakForm}
\begin{alignat}{2}
 B(u,\phi;v,\psi) + \langle \lambda,v \rangle_{\Gamma_C} &= \langle {f},v \rangle_{\Gamma_N} &\quad &\forall (v,\psi)\in \widetilde{H}^{1/2}(\Gamma_\Sigma) \times H^{-1/2}(\Gamma), \label{eq:compactWeakForm1}  \\
\langle u, \mu -\lambda \rangle_{\Gamma_C} & \leq \langle g,\mu_n-\lambda_n \rangle_{\Gamma_C} &\quad &\forall \mu \in M^+(\mathcal{F}) .  \label{eq:compactWeakForm2}
\end{alignat}
\end{subequations}
For the stabilization we make the trivial observation that $\lambda + Wu + (K + 1/2)^\prime \phi =0$ on $\Gamma_C$ which is the residual we use to stabilize the discrete problem allowing us to circumvent the discrete inf-sup condition.

\section{\texorpdfstring{Stabilized mixed $hp$-boundary element discretization}{Stabilized mixed hp-boundary element discretization}}
\label{sec:stabalizedDiscretization}
As usual let $\mathcal{T}_h$, $\mathcal{T}_k^C$ be locally quasi-uniform meshes for $\Gamma$, $\Gamma_C$ with the mesh size $h$, $k$ respectively. Let $p$, $q$ be locally quasi-uniform polynomial degree distributions over the meshes $\mathcal{T}_h$, $\mathcal{T}_k^C$, respectively, and let $ \Psi_T : \hat{T} \rightarrow T $ be the (bi)-linear, bijective mapping from the reference element $\hat{T}$ to the physical element $T$. For the discretization we use the following finite dimensional sets:
\begin{align}
X_{hp}&:=\left\{v_{hp} \in   C^0(\Gamma) : v_{hp}|_T \circ \Psi_T \in \left[\mathbb{P}_{p_T}(\hat{T})\right]^{d} \; \forall T \in \mathcal{T}_h,\ v_{hp}|_{\Gamma_D}=0 \right\}  ,\\
Y_{hp}&:=\left\{\phi_{hp} \in  L^{2}(\Gamma) : \phi_{hp}|_T \circ \Psi_T \in \left[\mathbb{P}_{p_T-1}(\hat{T})\right]^d \; \forall T \in \mathcal{T}_{h} \right\}, \\
M_{kq}^+(\mathcal{F})&:= \left\{ \mu^{kq} \in  L^2(\Gamma_C) : \mu^{kq}|_E \circ \Psi_T \in \left[\mathbb{P}_{q_T}(\hat{T})\right]^d \; \forall T \in \mathcal{T}_k^C,\ \mu^{kq}_n(x) \geq 0, \ |\mu^{kq}_t(x)| \leq \mathcal{F}(x)\; \text{for }x\in G_{kq} \right\} ,
\end{align}
where $G_{kq}$ is a set of discrete points on $\Gamma_C$, e.g.~the quadrature points used. Enforcing constraints of the primal variable~$u$ in such a finite set of points has been applied successfully in e.g.~\cite{gwinner2013hp, krebs2007p,maischak2005adaptive}. Let $\gamma$ be a piecewise constant function on $\Gamma_C$ with respect to the mesh $\mathcal{T}_h$. More precisely, we set $\gamma|_T=\bar{\gamma} h_Tp_T^{-2}$ for $T \in \mathcal{T}_h\cap \Gamma_C$ with global constant $\bar{\gamma}>0$. Let $\widetilde{W}$ and $\widetilde{K}^\prime$ be easily computable approximations of $W$, $K^\prime$ respectively, such that
\begin{alignat}{6}
  \widetilde{W}&:H^1(\Gamma_C) \rightarrow L^2(\Gamma_C) &\quad &\text{and} \quad  \| \gamma^{1/2} \widetilde{W} v_{hp} \|_{L^2(\Gamma)} &\leq& \; \bar{\gamma}^{1/2} C \| v_{hp}\|_{\widetilde{H}^{1/2}(\Gamma_\Sigma)} &\quad &\forall v_{hp} \in X_{hp} \label{eq:theo_W_tilde}\\
 (\widetilde{K}+1/2)^\prime&:L^2(\Gamma_C) \rightarrow L^2(\Gamma_C) &\quad &\text{and} \quad  \| \gamma^{1/2} (\widetilde{K}+1/2)^\prime \psi_{hp} \|_{L^2(\Gamma)} &\leq& \; \bar{\gamma}^{1/2} C \| \psi_{hp}\|_{H^{-1/2}(\Gamma)} &\quad &\forall v_{hp} \in Y_{hp}.  \label{eq:theo_K_tilde}
\end{alignat}
We discuss possible choices of $\widetilde{W}$ and $(\widetilde{K}+1/2)^\prime$ in Section~\ref{sec:approximation}, but $\widetilde{W}=W$ and $(\widetilde{K}+1/2)^\prime=(K+1/2)^\prime$ as chosen in \cite{Banz2015Stab} also works.
The discrete and stabilized mixed formulation is to find $(u_{hp},\phi_{hp},\lambda_{kq}) \in X_{hp} \times Y_{hp} \times M_{kq}^+(\mathcal{F})$ such that
\begin{subequations} \label{eq:DiscreteProblem}
\begin{gather}
 B(u_{hp},\phi_{hp};v_{hp},\psi_{hp})  + \langle \lambda_{kq},v_{hp} \rangle_{\Gamma_C} - \langle \gamma ( \lambda_{kq} + \widetilde{W}u_{hp} + (\widetilde{K} + 1/2)^\prime \phi_{hp}),  \widetilde{W}v_{hp} + (\widetilde{K} + 1/2)^\prime \psi_{hp} \rangle_{\Gamma_C} = \langle {f},v_{hp} \rangle_{\Gamma_N}  \label{eq:DiscreteMixedVarEq1} \\
\langle u_{hp}, \mu_{kq} -\lambda_{kq} \rangle_{\Gamma_C} - \langle \lambda_{kq}+ \widetilde{W}u_{hp} + (\widetilde{K} + 1/2)^\prime \phi_{hp}, \gamma ( \mu_{kq} - \lambda_{kq}) \rangle_{\Gamma_C}  \leq \langle g,\mu^{kq}_n-\lambda^{kq}_n \rangle_{\Gamma_C}  \label{eq:DiscreteContactConstraints}
\end{gather}
\end{subequations}
for all $(v_{hp},\psi_{hp},\mu_{kq}) \in X_{hp} \times Y_{hp} \times M_{kq}^+(\mathcal{F})$.

\begin{lemma} \label{lem:coercive}
There exists a $\bar{\gamma}_{max}$ such that for all $\bar{\gamma} \in (0,\bar{\gamma}_{max})$ there exists a constant $\alpha>0$ independent of $h$, $p$, $k$, $q$ and $\bar{\gamma}$ such that
\begin{align}
 B(v_{hp},\psi_{hp};v_{hp},\psi_{hp}) - \langle \gamma ( \widetilde{W}v_{hp} + (\widetilde{K} + 1/2)^\prime \psi_{hp}),  \widetilde{W}v_{hp} + (\widetilde{K} + 1/2)^\prime \psi_{hp} \rangle_{\Gamma_C} \geq \alpha \left( \| v_{hp} \|^2_{\widetilde{H}^{1/2}(\Gamma_\Sigma)} + \| \psi_{hp} \|^2_{H^{-1/2}(\Gamma)} \right)
\end{align}
for all $(v_{hp},\psi_{hp}) \in X_{hp} \times Y_{hp}$.
\begin{proof}
 By triangle inequality, Young's inequality and \eqref{eq:theo_W_tilde}-\eqref{eq:theo_K_tilde} there holds
 \begin{align*}
  \langle \gamma ( \widetilde{W}v_{hp} + (\widetilde{K} + 1/2)^\prime \psi_{hp}),  \widetilde{W}v_{hp} + (\widetilde{K} + 1/2)^\prime \psi_{hp} \rangle_{\Gamma_C} &\leq 2 \| \gamma^{1/2} \widetilde{W} v_{hp} \|_{L^2(\Gamma)}^2 + 2\| \gamma^{1/2} (\widetilde{K}+1/2)^\prime \psi_{hp} \|_{L^2(\Gamma)}^2 \\
  & \leq C \bar{\gamma} \left(\| v_{hp}\|_{\widetilde{H}^{1/2}(\Gamma_\Sigma)}^2 + \| \psi_{hp}\|_{H^{-1/2}(\Gamma)}^2 \right).
 \end{align*}
Thus we obtain with the coercivity of $W$ and $V$ that 
\begin{align*}
 B(v_{hp},\psi_{hp};v_{hp},\psi_{hp}) - \langle \gamma ( \widetilde{W}v_{hp} + (\widetilde{K} + 1/2)^\prime \psi_{hp}),  \widetilde{W}v_{hp} + (\widetilde{K} + 1/2)^\prime \psi_{hp} \rangle_{\Gamma_C} &\geq (\alpha_W-C \bar{\gamma})\| v_{hp}\|_{\widetilde{H}^{1/2}(\Gamma_\Sigma)}^2 \\
 & \quad + (\alpha_V -C \bar{\gamma})\| \psi_{hp}\|_{H^{-1/2}(\Gamma)}^2.
\end{align*}
Choosing $\bar{\gamma}$ sufficiently small yields the assertion when setting $\bar{\gamma}_{max}=\min\{\alpha_V,\alpha_W\}/2C$ for example.
\end{proof}
\end{lemma}
\begin{flushleft}
Consequently, we assume $\bar{\gamma}$ to be sufficiently small from now on for Lemma~\ref{lem:coercive} to be applicable.\end{flushleft}

\begin{theorem} \label{thm:discreteExistenceUniquenessStability}
There exists a unique solution $(u_{hp},\phi_{hp},\lambda_{kq}) \in X_{hp} \times Y_{hp} \times M_{kq}^+(\mathcal{F})$ to \eqref{eq:DiscreteProblem}. Moreover, $(u_{hp},\phi_{hp},\lambda_{kq})$ satisfy 
\begin{align} \label{eq:discrete_Stab_est1}
\ \| u_{hp} \|_{\widetilde{H}^{1/2}(\Gamma_\Sigma)}^2 + \| \phi_{hp} \|_{H^{-1/2}(\Gamma)}^2 + \| \gamma^{1/2} \lambda_{kq} \|_{L^2(\Gamma_C)}^2 \leq C \left( \| f\|_{\widetilde{H}^{-1/2}(\Gamma_N)}^2 -\langle g,\lambda^{kq}_n \rangle_{\Gamma_C} \right) 
\end{align}
and if $g \in C^0(\Gamma_C)$, $G_{kq}$ are the Gauss-Legendre quadrature points and thus $\hat{T}=[-1,1]$ in 2d and $\hat{T}=[-1,1]^2$ in 3d there even holds
\begin{align} \label{eq:discrete_Stab_est2}
\ \| u_{hp} \|_{\widetilde{H}^{1/2}(\Gamma_\Sigma)}^2 + \| \phi_{hp} \|_{H^{-1/2}(\Gamma)}^2 + \| \gamma^{1/2} \lambda_{kq} \|_{L^2(\Gamma_C)}^2 \leq C \left( \| f\|_{\widetilde{H}^{-1/2}(\Gamma_N)}^2 + \| \gamma^{-1/2} (g-\mathcal{I}_{kq}g) \|_{L^2(\Gamma_C)}^2 \right) 
\end{align}
each with a constant $C>0$ independent of $h$, $p$, $k$, $q$ and $\bar{\gamma}$ and with $\mathcal{I}_{kq}$ the classical interpolation operator in the points $G_{kq}$.
\begin{proof}
\textbf{Unique existence:} Since $B(\cdot,\cdot;\cdot,\cdot)$ is not symmetric we cannot simply argue with the saddle point theory, but work with Banach's fixed point theorem. \eqref{eq:DiscreteContactConstraints} can be written as
\begin{align*}
\langle \lambda_{kq} - \lambda_{kq} + r \left( \gamma \lambda_{kq}- u_{hp} + g\cdot n + \gamma \widetilde{W}u_{hp} + \gamma (\widetilde{K} + 1/2)^\prime \phi_{hp}  \right), \mu_{kq} -\lambda_{kq} \rangle_{\Gamma_C} \geq 0  \quad \forall \mu_{kq} \in M_{kq}^+(\mathcal{F})
\end{align*}
with an arbitrary constant $r>0$. And thus, as $M_{kq}^+(\mathcal{F})$ is convex, \eqref{eq:DiscreteContactConstraints} can be written as the projection problem
\begin{align} \label{eq:projection_problem}
 \lambda_{kq} = \mathcal{P}_{M_{kq}^+(\mathcal{F})} \left( \lambda_{kq} + r \left( u_{hp} - g\cdot n - \gamma \widetilde{W}u_{hp} - \gamma (\widetilde{K} + 1/2)^\prime \phi_{hp}  -\gamma \lambda_{kq}\right) \right).
\end{align}
For a given $\lambda_i \in M_{kq}^+(\mathcal{F})$, let $(u_i,\phi_i)\in X_{hp} \times Y_{hp}$ be the unique solution to \eqref{eq:DiscreteMixedVarEq1}. Existence and uniqueness of $(u_i,\phi_i)$ follow from Lemma~\ref{lem:coercive} as \eqref{eq:DiscreteMixedVarEq1} is a finite dimensional positive definite system of linear equations. 
To shorten the notation let $\delta_\lambda=\lambda_1-\lambda_2$, $\delta_u=u_1-u_2$, $\delta_\phi=\phi_1 - \phi_2$, $\gamma_{\max}=\max_{\Gamma_C} \gamma$, $\gamma_{\min}=\min_{\Gamma_C} \gamma$ and $\Delta^2 = \left( \| \delta_u \|^2_{\widetilde{H}^{1/2}(\Gamma_\Sigma)} + \| \delta_\phi \|^2_{H^{-1/2}(\Gamma)} \right)\| \delta_\lambda \|_{L^2(\Gamma_C)}^{-2}$.
By \eqref{eq:DiscreteMixedVarEq1}, Lemma~\ref{lem:coercive}, Cauchy-Schwarz inequality, \eqref{eq:theo_W_tilde}-\eqref{eq:theo_K_tilde} and Young's inequality with $\epsilon>0$ arbitrary there holds
\begin{align*}
 \alpha \left(\| \delta_u \|^2_{\widetilde{H}^{1/2}(\Gamma_\Sigma)} + \| \delta_\phi \|^2_{H^{-1/2}(\Gamma)}  \right)  &\leq B(\delta_u,\delta_\phi; \delta_u, \delta_\phi) - \langle \gamma (\widetilde{W}\delta_u + (\widetilde{K} + 1/2)^\prime \delta_\phi) , \widetilde{W}\delta_u + (\widetilde{K} + 1/2)^\prime \delta_\phi \rangle_{\Gamma_C} \\
 &\leq - \langle \delta_\lambda, \delta_u \rangle_{\Gamma_C} + \langle \gamma \delta_\lambda ,  \widetilde{W}\delta_u + (\widetilde{K} + 1/2)^\prime \delta_\phi \rangle_{\Gamma_C} \\
 &\leq \frac{1/4 +C \bar{\gamma} \gamma_{\max} }{\epsilon}\| \delta_\lambda \|_{L^2(\Gamma_C)}^2 + 2\epsilon \| \delta_u \|^2_{\widetilde{H}^{1/2}(\Gamma_\Sigma)}  + \epsilon \| \delta_\phi \|^2_{H^{-1/2}(\Gamma)}.
\end{align*}
Since $\gamma_{\max}$ can be bounded by $\bar{\gamma} | \Gamma_C|$, there exists a constant $C>0$ depending on $\bar{\gamma}$ such that
\begin{align} \label{eq:elliptic_stability_discrete}
 \| \delta_u \|^2_{\widetilde{H}^{1/2}(\Gamma_\Sigma)} + \| \delta_\phi \|^2_{H^{-1/2}(\Gamma)} \leq C \| \delta_\lambda \|_{L^2(\Gamma_C)}^2.
\end{align}
For $r>0$ sufficiently small, the projection \eqref{eq:projection_problem} is a contraction, as
\begin{align*}
 \| \delta_\lambda \|_{L^2(\Gamma_C)}^2 &\leq \| (1-r\gamma) \delta_\lambda + r((I- \gamma \widetilde{W})\delta_u - \gamma (\widetilde{K}+1/2)^\prime \delta_\phi  ) \|_{L^2(\Gamma_C)}^2\\
 &= \| (1-r \gamma) \delta_\lambda \|_{L^2(\Gamma_C)}^2 +2r \langle (1-r \gamma) \delta_\lambda, (I- \gamma \widetilde{W})\delta_u - \gamma (\widetilde{K}+1/2)^\prime \delta_\phi \rangle_{\Gamma_C} + r^2 \|(I- \gamma \widetilde{W})\delta_u - \gamma (\widetilde{K}+1/2)^\prime \delta_\phi  \|_{L^2(\Gamma_C)}^2 \\
 &=\| (1-r \gamma) \delta_\lambda \|_{L^2(\Gamma_C)}^2 -2r \left( B(\delta_u,\delta_\phi; \delta_u,\delta_\phi) - \langle \gamma ( \widetilde{W}\delta_u +(\widetilde{K}+1/2)^\prime \delta_\phi ) ,\widetilde{W}\delta_u +(\widetilde{K}+1/2)^\prime \delta_\phi \rangle_{\Gamma_C} \right) \\
 & \quad +2r^2 \langle \gamma \delta_\lambda, - (I- \gamma \widetilde{W})\delta_u - \gamma (\widetilde{K}+1/2)^\prime \delta_\phi \rangle_{\Gamma_C} + r^2 \|(I- \gamma \widetilde{W})\delta_u - \gamma (\widetilde{K}+1/2)^\prime \delta_\phi  \|_{L^2(\Gamma_C)}^2\\
 & \leq \| (1-r \gamma) \delta_\lambda \|_{L^2(\Gamma_C)}^2 - 2r\alpha \left( \| \delta_u \|^2_{\widetilde{H}^{1/2}(\Gamma_\Sigma)} + \| \delta_\phi \|^2_{H^{-1/2}(\Gamma)} \right) + r^2\| \gamma \delta_\lambda \|_{L^2(\Gamma_C)}^2 + 2r^2 \|(I- \gamma \widetilde{W})\delta_u - \gamma (\widetilde{K}+1/2)^\prime \delta_\phi  \|_{L^2(\Gamma_C)}^2\\
 & \leq \left( (1-r \gamma_{\min})^2 + r^2 \gamma_{\max}^2 \right)\|  \delta_\lambda \|_{L^2(\Gamma_C)}^2 - 2r\alpha \left( \| \delta_u \|^2_{\widetilde{H}^{1/2}(\Gamma_\Sigma)} + \| \delta_\phi \|^2_{H^{-1/2}(\Gamma)} \right) \\
 & \quad + 6r^2  \left((1+ \gamma_{\max} \bar{\gamma}C) \|\delta_u\|_{\widetilde{H}^{1/2}(\Gamma_\Sigma)}^2 + \gamma_{\max} \bar{\gamma}C \|\delta_\phi\|_{H^{-1/2}(\Gamma)}^2\right)\\
 & \leq  \left[ \left( (1-r \gamma_{\min})^2 + r^2 \gamma_{\max}^2 \right) - 2r\alpha \left( \| \delta_u \|^2_{\widetilde{H}^{1/2}(\Gamma_\Sigma)} + \| \delta_\phi \|^2_{H^{-1/2}(\Gamma)} \right)\|  \delta_\lambda \|_{L^2(\Gamma_C)}^{-2} \right. \\
 & \left. \quad + 6r^2(1+ \gamma_{\max} \bar{\gamma}C)  \left( \|\delta_u\|_{\widetilde{H}^{1/2}(\Gamma_\Sigma)}^2 + \|\delta_\phi\|_{H^{-1/2}(\Gamma)}^2\right)\|  \delta_\lambda \|_{L^2(\Gamma_C)}^{-2} \right] \|  \delta_\lambda \|_{L^2(\Gamma_C)}^2\\
 &=\left[1 -2r (\gamma_{\min}+\alpha\Delta^2) + r^2( \gamma_{\min}^2+ \gamma_{\max}^2 +  6(1+ \gamma_{\max} \bar{\gamma}C)\Delta^2) \right] \|  \delta_\lambda \|_{L^2(\Gamma_C)}^2\\
 & \leq \left[1 -2r\gamma_{\min} + r^2( \gamma_{\min}^2+C) \right] \|  \delta_\lambda \|_{L^2(\Gamma_C)}^2
\end{align*}
where the third line follows from  \eqref{eq:DiscreteMixedVarEq1}, the fifth line from Lemma~\ref{lem:coercive}, Cauchy-Schwarz inequality and Young's inequality, the sixth line from \eqref{eq:theo_W_tilde}-\eqref{eq:theo_K_tilde} and the final line from \eqref{eq:elliptic_stability_discrete}.
The square bracket expression is minimal for $r=\gamma_{\min}( \gamma_{\min}^2+ C)^{-1}$ and takes the value $ 0 \leq 1-\gamma_{\min}^2( \gamma_{\min}^2+ C)^{-1} \leq \bar{C} <1 $. 
Thus, existence and uniqueness of $\lambda_{kq}$ is a consequence of Banach's fixed point theorem which in turn yields existence and uniqueness of $(u_{hp},\phi_{hp})$ by Lemma~\ref{lem:coercive}.
We point out that this argument is also valid for the non-stabilized case, i.e.~$\bar{\gamma}=0$ if $\Delta^2 \geq C >0$ which generally requires the discrete inf-sup condition to hold. The second noteworthy point is that the above argument only works for the finite dimensional case, i.e.~$\gamma_{\min}>0$, and not for the limit case $h=0$ or $p=\infty$.\\ 
\textbf{Stability:} Choosing $v_{hp}=u_{hp}$, $\psi_{hp}=\phi_{hp}$, $\mu_{kq}=0$ in \eqref{eq:DiscreteProblem} we obtain from \eqref{eq:DiscreteContactConstraints} 
\begin{align*}
 \langle u_{hp}, \lambda_{kq} \rangle_{\Gamma_C} \geq   \langle \gamma  \lambda_{kq}, \lambda_{kq}+ \widetilde{W}u_{hp} + (\widetilde{K} + 1/2)^\prime \phi_{hp} \rangle_{\Gamma_C}  + \langle g,\lambda^{kq}_n \rangle_{\Gamma_C} 
\end{align*}
which inserted into \eqref{eq:DiscreteMixedVarEq1} yields with Lemma~\ref{lem:coercive}
\begin{align*}
 \langle {f},u_{hp} \rangle_{\Gamma_N} &= B(u_{hp},\phi_{hp};u_{hp},\phi_{hp})  + \langle \lambda_{kq},u_{hp} \rangle_{\Gamma_C} - \langle \gamma ( \lambda_{kq} + \widetilde{W}u_{hp} + (\widetilde{K} + 1/2)^\prime \phi_{hp}),  \widetilde{W}u_{hp} + (\widetilde{K} + 1/2)^\prime \phi_{hp} \rangle_{\Gamma_C}\\
 &\geq  B(u_{hp},\phi_{hp};u_{hp},\phi_{hp}) - \langle \gamma ( \widetilde{W}u_{hp} + (\widetilde{K} + 1/2)^\prime \phi_{hp}),  \widetilde{W}u_{hp} + (\widetilde{K} + 1/2)^\prime \phi_{hp} \rangle_{\Gamma_C} +\langle \gamma  \lambda_{kq}, \lambda_{kq}\rangle_{\Gamma_C}+ \langle g,\lambda^{kq}_n \rangle_{\Gamma_C} \\
 & \geq \alpha \left( \| u_{hp} \|^2_{\widetilde{H}^{1/2}(\Gamma_\Sigma)} + \| \phi_{hp} \|^2_{H^{-1/2}(\Gamma)} \right) + \| \gamma^{1/2} \lambda_{kq} \|_{L^2(\Gamma_C)}^2+ \langle g,\lambda^{kq}_n \rangle_{\Gamma_C}.
\end{align*}
Applying Cauchy-Schwarz and Young's inequality yields the first stability estimate \eqref{eq:discrete_Stab_est1}. 
For $\hat{T}=[-1,1]$ in 2d and $\hat{T}=[-1,1]^2$ in 3d, the right hand side in \eqref{eq:discrete_Stab_est1} can be further estimated, if $g \in C^0(\Gamma_C)$ and $G_{kq}$ (the points in which constraints on $\lambda_{kq}$ are enforced) are the Gauss-Legendre quadrature points. Let $\mathcal{I}_{kq}$ be the classical interpolation operator in the points $G_{kq}$, than for $\epsilon >0$ arbitrary we have
\begin{align*}
 \langle -g,\lambda^{kq}_n \rangle_{\Gamma_C} = \langle \mathcal{I}_{kq}g-g -\mathcal{I}_{kq}g,\lambda^{kq}_n \rangle_{\Gamma_C} \leq \langle \mathcal{I}_{kq}g-g,\lambda^{kq}_n \rangle_{\Gamma_C} \leq \frac{1}{4 \epsilon} \| \gamma^{-1/2} (g-\mathcal{I}_{kq}g) \|_{L^2(\Gamma_C)}^2 + \epsilon \| \gamma^{1/2} \lambda_{kq} \|_{L^2(\Gamma_C)}^2
\end{align*}
as the Gauss-Legendre quadrature on $[-1,1]$ with $q+1$ points is exact for polynomials of order $2q+1$ and has positive weights, and $\lambda_n^{kq}$ and $g$ are non negative in these quadrature points. 
\end{proof}
\end{theorem}

\section{A priori error estimate}
\label{sec:Apriorierrorestimate}

For the proof of the a priori error estimate we need the following two lemmas.
Subtracting \eqref{eq:DiscreteMixedVarEq1} from \eqref{eq:compactWeakForm1} yields by the conformity of the discretization $X_{hp} \times Y_{hp} \subset \widetilde{H}^{1/2}(\Gamma_\Sigma) \times H^{-1/2}(\Gamma)$ that:
\begin{lemma} \label{lem:GalerkinOrtho}
 Let $(u,\phi,\lambda)$, $(u_{hp},\phi_{hp},\lambda_{kq})$ be the solution to \eqref{eq:compactWeakForm}, \eqref{eq:DiscreteProblem}, respectively. Then there holds
 \begin{align*}
  B(u-u_{hp},\phi - \phi_{hp};v_{hp},\psi_{hp})  + \langle \lambda - \lambda_{kq},v_{hp} \rangle_{\Gamma_C} + \langle \gamma ( \lambda_{kq} + \widetilde{W}u_{hp} + (\widetilde{K} + 1/2)^\prime \phi_{hp}),  \widetilde{W}v_{hp} + (\widetilde{K} + 1/2)^\prime \psi_{hp} \rangle_{\Gamma_C} = 0
 \end{align*}
for all $(v_{hp},\psi_{hp}) \in X_{hp} \times Y_{hp}$.
\end{lemma}

\begin{lemma} \label{lem:aprioriErrorLambda}
 Let $(u,\phi,\lambda)$, $(u_{hp},\phi_{hp},\lambda_{kq})$ be the solution to \eqref{eq:compactWeakForm}, \eqref{eq:DiscreteProblem}, respectively. If $\lambda$, $(Wu)|_{\Gamma_C}$ and $((K+1/2)^\prime \phi)|_{\Gamma_C} \in  L^2(\Gamma_C)$, then there holds
\begin{align*}
 &2^{-1} \| \gamma^{1/2} ( \lambda - \lambda_{kq} ) \|_{L^2(\Gamma_C)}^2 \\
 & \leq   \langle u_n-g, \lambda_n^{kq} - \mu_n \rangle_{\Gamma_C} + \langle u_t, \lambda_t^{kq} - \mu_t \rangle_{\Gamma_C} + \langle g, \mu_n^{kq} - \lambda_n  \rangle_{\Gamma_C} +\| \lambda - \mu_{kq} \|_{\widetilde{H}^{-1/2}(\Gamma_C)} \| u_{hp} \|_{\widetilde{H}^{1/2}(\Gamma_\Sigma)}+ 4^{-1}  \| \gamma^{1/2} (\mu_{kq}-\lambda) \|_{L^2(\Gamma_C)}^2\\
 & \quad  + C\| \gamma^{1/2}(\lambda- \mu_{kq}) \|_{L^2(\Gamma_C)} \left( \bar{\gamma}^{1/2}\|u_{hp}\|_{\widetilde{H}^{1/2}(\Gamma_\Sigma)} + \bar{\gamma}^{1/2}\| \phi_{hp}\|_{H^{-1/2}(\Gamma)} + \|\gamma^{1/2}\lambda_{kq}\|_{L^2(\Gamma_C)}\right) +\langle \lambda -\lambda_{kq},u-u_{hp} \rangle_{\Gamma_C}    \\
 & \quad+ 2\left(  C \bar{\gamma} \|u- v_{hp} \|^2_{\widetilde{H}^{1/2}(\Gamma_\Sigma)} + C \bar{\gamma} \|u - u_{hp} \|^2_{\widetilde{H}^{1/2}(\Gamma_\Sigma)}  + C \bar{\gamma} \|\phi - \psi_{hp} \|^2_{H^{-1/2}(\Gamma)} + C \bar{\gamma} \|\phi - \phi_{hp} \|^2_{H^{-1/2}(\Gamma)}\right)\\
 &\quad  + 4 \left( \| \gamma^{1/2} W(u-v_{hp}) \|_{L^2(\Gamma_C)}^2+\| \gamma^{1/2} (K + 1/2)^\prime (\phi-\psi_{hp}) \|_{L^2(\Gamma_C)}^2+\| \gamma^{1/2}  (\widetilde{W}-W)u_{hp} \|_{L^2(\Gamma_C)}^2
  +\| \gamma^{1/2}  (\widetilde{K}^\prime - K^\prime) \phi_{hp} \|_{L^2(\Gamma_C)}^2 \right)
\end{align*}
for all $(v_{hp},\psi_{hp},\mu_{kq}) \in X_{hp} \times Y_{hp} \times M_{kq}^+(\mathcal{F})$ and all $\mu \in M^+(\mathcal{F})$ with a constant $C>0$ independent of $h$, $p$, $k$, $q$ and $\bar{\gamma}$.
\begin{proof}
There holds
\begin{gather*}
  \langle u_{hp}, \lambda_{kq} - \mu_{kq}  \rangle_{\Gamma_C}  + \langle u , \lambda  - \mu  \rangle_{\Gamma_C} =  \langle \lambda -\lambda_{kq},u-u_{hp} \rangle_{\Gamma_C} + \langle u, \lambda_{kq} - \mu \rangle_{\Gamma_C} + \langle u_{hp}, \lambda - \mu_{kq} \rangle_{\Gamma_C}, \\
\langle \gamma \lambda , \lambda \rangle_{\Gamma_C} -2 \langle \gamma \lambda , \lambda_{kq} \rangle_{\Gamma_C}  + \langle \gamma \mu_{kq}, \lambda_{kq}  \rangle_{\Gamma_C} 
=  \langle \gamma (\lambda-\lambda_{kq}) , \lambda \rangle_{\Gamma_C} + \langle \gamma \lambda_{kq} , \mu_{kq}- \lambda \rangle_{\Gamma_C}
\end{gather*}
and 
\begin{align*}
  \langle \mu_{kq} - \lambda_{kq}, \gamma ( \widetilde{W}u_{hp} + (\widetilde{K} + 1/2)^\prime \phi_{hp} )\rangle_{\Gamma_C}
 &= \langle \mu_{kq}-\lambda  , \gamma ( Wu_{hp} + (K + 1/2)^\prime \phi_{hp} )\rangle_{\Gamma_C} \\
 & \quad + \langle  \lambda - \lambda_{kq}, \gamma ( Wu_{hp} + (K + 1/2)^\prime \phi_{hp} )\rangle_{\Gamma_C} \\
 & \quad 
 +\langle \mu_{kq} - \lambda_{kq}, \gamma ( (\widetilde{W}-W)u_{hp} + (\widetilde{K}^\prime - K^\prime) \phi_{hp} )\rangle_{\Gamma_C}.
\end{align*}
Rearranging \eqref{eq:DiscreteContactConstraints} and adding \eqref{eq:compactWeakForm2} to the right hand side yields
\begin{align*}
 \langle \gamma \lambda_{kq},\lambda_{kq} \rangle_{\Gamma_C} &\leq \langle g, \mu_n^{kq} - \lambda_n + \mu_n - \lambda_n^{kq}  \rangle_{\Gamma_C} + \langle \gamma (\mu_{kq} - \lambda_{kq}),  \widetilde{W}u_{hp} + (\widetilde{K} + 1/2)^\prime \phi_{hp} \rangle_{\Gamma_C} + \langle \gamma \mu_{kq}, \lambda_{kq}  \rangle_{\Gamma_C}\\
 &\quad + \langle u_{hp}, \lambda_{kq} - \mu_{kq}  \rangle_{\Gamma_C}  + \langle u , \lambda  - \mu  \rangle_{\Gamma_C}\\
 &= \langle g, \mu_n^{kq} - \lambda_n + \mu_n - \lambda_n^{kq}  \rangle_{\Gamma_C} + \langle \gamma (\mu_{kq} - \lambda_{kq}),  \widetilde{W}u_{hp} + (\widetilde{K} + 1/2)^\prime \phi_{hp} \rangle_{\Gamma_C} + \langle \gamma \mu_{kq}, \lambda_{kq}  \rangle_{\Gamma_C}\\
 &\quad + \langle \lambda -\lambda_{kq},u-u_{hp} \rangle_{\Gamma_C} + \langle u, \lambda_{kq} - \mu \rangle_{\Gamma_C} + \langle u_{hp}, \lambda - \mu_{kq} \rangle_{\Gamma_C}.
\end{align*}
Thus
\begin{align*}
 \| \gamma^{1/2} ( \lambda - \lambda_{kq} ) \|_{L^2(\Gamma_C)}^2 &= \langle \gamma \lambda , \lambda \rangle_{\Gamma_C} -2 \langle \gamma \lambda , \lambda_{kq} \rangle_{\Gamma_C} + \langle \gamma \lambda_{kq} , \lambda_{kq} \rangle_{\Gamma_C} \\
 &\leq \langle \gamma \lambda , \lambda \rangle_{\Gamma_C} -2 \langle \gamma \lambda , \lambda_{kq} \rangle_{\Gamma_C} +\langle g, \mu_n^{kq} - \lambda_n + \mu_n - \lambda_n^{kq}  \rangle_{\Gamma_C} + \langle \mu_{kq} - \lambda_{kq}, \gamma ( \widetilde{W}u_{hp} + (\widetilde{K} + 1/2)^\prime \phi_{hp} )\rangle_{\Gamma_C} \\
 &\quad + \langle \gamma \mu_{kq}, \lambda_{kq}  \rangle_{\Gamma_C} +\langle \lambda -\lambda_{kq},u-u_{hp} \rangle_{\Gamma_C} + \langle u, \lambda_{kq} - \mu \rangle_{\Gamma_C} + \langle u_{hp}, \lambda - \mu_{kq} \rangle_{\Gamma_C}\\
 &=\langle \gamma (\lambda-\lambda_{kq}) , \lambda \rangle_{\Gamma_C} + \langle \gamma \lambda_{kq} , \mu_{kq}- \lambda \rangle_{\Gamma_C} +\langle g, \mu_n^{kq} - \lambda_n + \mu_n - \lambda_n^{kq}  \rangle_{\Gamma_C} \\
 & \quad + \langle \mu_{kq}-\lambda  , \gamma ( Wu_{hp} + (K + 1/2)^\prime \phi_{hp} )\rangle_{\Gamma_C}  + \langle  \lambda - \lambda_{kq}, \gamma ( Wu_{hp} + (K + 1/2)^\prime \phi_{hp} )\rangle_{\Gamma_C} \\
 & \quad +\langle \mu_{kq} - \lambda_{kq}, \gamma ( (\widetilde{W}-W)u_{hp} + (\widetilde{K}^\prime - K^\prime) \phi_{hp} )\rangle_{\Gamma_C} \\
 &\quad  +\langle \lambda -\lambda_{kq},u-u_{hp} \rangle_{\Gamma_C} + \langle u, \lambda_{kq} - \mu \rangle_{\Gamma_C} + \langle u_{hp}, \lambda - \mu_{kq} \rangle_{\Gamma_C}\\
 &= \langle g, \mu_n^{kq} - \lambda_n + \mu_n - \lambda_n^{kq}  \rangle_{\Gamma_C} + \underbrace{\langle \mu_{kq}-\lambda  , -u_{hp}+\gamma (\lambda_{kq}+ Wu_{hp} + (K + 1/2)^\prime \phi_{hp} )\rangle_{\Gamma_C}}_{\rm I} \\
 & \quad + \underbrace{\langle \gamma (\lambda - \lambda_{kq}), \lambda + Wu_{hp} + (K + 1/2)^\prime \phi_{hp} \rangle_{\Gamma_C}}_{\rm II}
 +\underbrace{\langle \mu_{kq} - \lambda_{kq}, \gamma ( (\widetilde{W}-W)u_{hp} + (\widetilde{K}^\prime - K^\prime) \phi_{hp} )\rangle_{\Gamma_C}}_{\rm III} \\
 &\quad  +\langle \lambda -\lambda_{kq},u-u_{hp} \rangle_{\Gamma_C} + \langle u, \lambda_{kq} - \mu \rangle_{\Gamma_C}.
 \end{align*}
 With Cauchy-Schwarz inequality, triangle inequality and Young's inequality we obtain
 \begin{align*}
  {\rm III} &\leq \left( \| \gamma^{1/2} (\mu_{kq}-\lambda) \|_{L^2(\Gamma_C)} + \| \gamma^{1/2} (\lambda-\lambda_{kq}) \|_{L^2(\Gamma_C)} \right) \left( \| \gamma^{1/2}  (\widetilde{W}-W)u_{hp} \|_{L^2(\Gamma_C)}
  +\| \gamma^{1/2}  (\widetilde{K}^\prime - K^\prime) \phi_{hp} \|_{L^2(\Gamma_C)} \right) \\
  & \leq 2\epsilon \left( \| \gamma^{1/2} (\mu_{kq}-\lambda) \|_{L^2(\Gamma_C)}^2 + \| \gamma^{1/2} (\lambda-\lambda_{kq}) \|_{L^2(\Gamma_C)}^2 \right) + \frac{1}{2\epsilon} \left( \| \gamma^{1/2}  (\widetilde{W}-W)u_{hp} \|_{L^2(\Gamma_C)}^2
  +\| \gamma^{1/2}  (\widetilde{K}^\prime - K^\prime) \phi_{hp} \|_{L^2(\Gamma_C)}^2 \right).
 \end{align*} 
 Analogously but with $\lambda = -Wu - (K+1/2)^\prime \phi$ a.e.~on $\Gamma_C$ by \eqref{eq:compactWeakForm1} there holds
 \begin{align*}
 {\rm II}
  &= \langle \gamma (\lambda - \lambda_{kq}),  W(u_{hp}-u) + (K + 1/2)^\prime (\phi_{hp}-\phi) \rangle_{\Gamma_C} \\
 & \leq  \frac{1}{4\epsilon}\left( \| \gamma^{1/2} W(u_{hp}-u) \|_{L^2(\Gamma_C)}^2 + \| \gamma^{1/2} (K + 1/2)^\prime (\phi_{hp}-\phi) \|_{L^2(\Gamma_C)}^2\right) +2\epsilon \| \gamma^{1/2} (\lambda-\lambda_{kq}) \|_{L^2(\Gamma_C)}^2.
 \end{align*}
Inserting these two estimates into the estimate for $\| \gamma^{1/2} ( \lambda - \lambda_{kq} ) \|_{L^2(\Gamma_C)}^2$ yields
\begin{align*}
 (1-4\epsilon)\| \gamma^{1/2} ( \lambda - \lambda_{kq} ) \|_{L^2(\Gamma_C)}^2
 &\leq   \langle g, \mu_n^{kq} - \lambda_n + \mu_n - \lambda_n^{kq}  \rangle_{\Gamma_C} + \langle \mu_{kq}-\lambda  , -u_{hp}+\gamma (\lambda_{kq}+ Wu_{hp} + (K + 1/2)^\prime \phi_{hp} )\rangle_{\Gamma_C}  \\
 &\quad  +\langle \lambda -\lambda_{kq},u-u_{hp} \rangle_{\Gamma_C} + \langle u, \lambda_{kq} - \mu \rangle_{\Gamma_C}+ 2\epsilon  \| \gamma^{1/2} (\mu_{kq}-\lambda) \|_{L^2(\Gamma_C)}^2 \\
 & \quad+ \frac{1}{4\epsilon}\left( \| \gamma^{1/2} W(u_{hp}-u) \|_{L^2(\Gamma_C)}^2 + \| \gamma^{1/2} (K + 1/2)^\prime (\phi_{hp}-\phi) \|_{L^2(\Gamma_C)}^2\right)\\
 &\quad   + \frac{1}{2\epsilon} \left( \| \gamma^{1/2}  (\widetilde{W}-W)u_{hp} \|_{L^2(\Gamma_C)}^2
  +\| \gamma^{1/2}  (\widetilde{K}^\prime - K^\prime) \phi_{hp} \|_{L^2(\Gamma_C)}^2 \right).
\end{align*}
Due to triangle inequality and \eqref{eq:theo_W_tilde} there holds for all $v_{hp} \in X_{hp}$
\begin{align}
   \| \gamma^{1/2} W(u_{hp}-u) \|_{L^2(\Gamma_C)}^2 &\leq 2 \| \gamma^{1/2} W(u-v_{hp}) \|_{L^2(\Gamma_C)}^2+2\| \gamma^{1/2} W(v_{hp}-u_{hp}) \|_{L^2(\Gamma_C)}^2 \nonumber \\
 &\leq    2 \| \gamma^{1/2} W(u-v_{hp}) \|_{L^2(\Gamma_C)}^2+ C \bar{\gamma} \|u- v_{hp} \|^2_{\widetilde{H}^{1/2}(\Gamma_\Sigma)} + C \bar{\gamma} \|u - u_{hp} \|^2_{\widetilde{H}^{1/2}(\Gamma_\Sigma)}, \label{eq:scaled_W_error_est}
\end{align}
and, analogously with \eqref{eq:theo_K_tilde} and $\psi_{hp} \in Y_{hp}$ arbitrary
\begin{align}
 \| \gamma^{1/2} (K + 1/2)^\prime (\phi_{hp}-\phi) \|_{L^2(\Gamma_C)}^2 \leq 2 \| \gamma^{1/2} (K + 1/2)^\prime (\phi-\psi_{hp}) \|_{L^2(\Gamma_C)}^2 + C \bar{\gamma} \|\phi - \psi_{hp} \|^2_{H^{-1/2}(\Gamma)} + C \bar{\gamma} \|\phi - \phi_{hp} \|^2_{H^{-1/2}(\Gamma)}.  \label{eq:scaled_K_error_est}
\end{align}
 By the inverse inequality \cite[Proof of Thm.~5]{Banz2015Stab} there holds
\begin{align*}
  {\rm I}  \leq \| \lambda - \mu_{kq} \|_{\widetilde{H}^{-1/2}(\Gamma_C)} \| u_{hp} \|_{\widetilde{H}^{1/2}(\Gamma_\Sigma)} + \| \gamma^{1/2}(\lambda- \mu_{kq}) \|_{L^2(\Gamma_C)} \left( C \bar{\gamma}^{1/2}\|u_{hp}\|_{\widetilde{H}^{1/2}(\Gamma_\Sigma)} + C \bar{\gamma}^{1/2}\| \phi_{hp}\|_{H^{-1/2}(\Gamma)} + \|\gamma^{1/2}\lambda_{kq}\|_{L^2(\Gamma_C)}\right).
\end{align*}
Finally,  
\begin{align*}
 \langle g, \mu_n^{kq} - \lambda_n + \mu_n - \lambda_n^{kq}  \rangle_{\Gamma_C}+\langle u, \lambda_{kq} - \mu \rangle_{\Gamma_C} = \langle u_n-g, \lambda_n^{kq} - \mu_n \rangle_{\Gamma_C} + \langle u_t, \lambda_t^{kq} - \mu_t \rangle_{\Gamma_C} + \langle g, \mu_n^{kq} - \lambda_n  \rangle_{\Gamma_C}.
\end{align*}
Thus
\begin{align*}
 (1- &4\epsilon)  \| \gamma^{1/2} ( \lambda - \lambda_{kq} ) \|_{L^2(\Gamma_C)}^2 \\
 &  \leq   \langle u_n-g, \lambda_n^{kq} - \mu_n \rangle_{\Gamma_C} + \langle u_t, \lambda_t^{kq} - \mu_t \rangle_{\Gamma_C} + \langle g, \mu_n^{kq} - \lambda_n  \rangle_{\Gamma_C} +\| \lambda - \mu_{kq} \|_{\widetilde{H}^{-1/2}(\Gamma_C)} \| u_{hp} \|_{\widetilde{H}^{1/2}(\Gamma_\Sigma)}+\langle \lambda -\lambda_{kq},u-u_{hp} \rangle_{\Gamma_C}\\
 & \quad  + \| \gamma^{1/2}(\lambda- \mu_{kq}) \|_{L^2(\Gamma_C)} \left( C \bar{\gamma}^{1/2}\|u_{hp}\|_{\widetilde{H}^{1/2}(\Gamma_\Sigma)} + C \bar{\gamma}^{1/2}\| \phi_{hp}\|_{H^{-1/2}(\Gamma)} + \|\gamma^{1/2}\lambda_{kq}\|_{L^2(\Gamma_C)}\right) + 2\epsilon  \| \gamma^{1/2} (\mu_{kq}-\lambda) \|_{L^2(\Gamma_C)}^2    \\
 & \quad  + \frac{1}{4\epsilon}\left(  C \bar{\gamma} \|u- v_{hp} \|^2_{\widetilde{H}^{1/2}(\Gamma_\Sigma)} + C \bar{\gamma} \|u - u_{hp} \|^2_{\widetilde{H}^{1/2}(\Gamma_\Sigma)}  + C \bar{\gamma} \|\phi - \psi_{hp} \|^2_{H^{-1/2}(\Gamma)} + C \bar{\gamma} \|\phi - \phi_{hp} \|^2_{H^{-1/2}(\Gamma)}\right)\\
 &\quad   + \frac{1}{2\epsilon} \left( \| \gamma^{1/2} W(u-v_{hp}) \|_{L^2(\Gamma_C)}^2+\| \gamma^{1/2} (K + 1/2)^\prime (\phi-\psi_{hp}) \|_{L^2(\Gamma_C)}^2+\| \gamma^{1/2}  (\widetilde{W}-W)u_{hp} \|_{L^2(\Gamma_C)}^2 \right. \\
 & \quad \left.  +\| \gamma^{1/2}  (\widetilde{K}^\prime - K^\prime) \phi_{hp} \|_{L^2(\Gamma_C)}^2 \right)
\end{align*}
and choosing $\epsilon=8^{-1}$ yields the assertion.
\end{proof}
\end{lemma}

Now, we can proof a C\'ea-Lemma like a priori error estimate.
\begin{theorem} \label{thm:AprioriError}
 Let $(u,\phi,\lambda)$, $(u_{hp},\phi_{hp},\lambda_{kq})$ be the solution to \eqref{eq:compactWeakForm}, \eqref{eq:DiscreteProblem}, respectively. If $\bar{\gamma}$ is sufficiently small, potentially smaller than for Lemma~\ref{lem:coercive} needed, and if $\lambda$, $(Wu)|_{\Gamma_C}$ and $((K+1/2)^\prime \phi)|_{\Gamma_C} \in  L^2(\Gamma_C)$, then there holds
\begin{align*}
 C & \left(\| u-u_{hp}\|_{\widetilde{H}^{1/2}(\Gamma_\Sigma)}^2 +  \| \phi-\phi_{hp} \|_{H^{-1/2}(\Gamma)}^2 +\| h^{1/2}p^{-1} ( \lambda - \lambda_{kq} ) \|_{L^2(\Gamma_C)}^2 \right)\\
&  \leq \| u-v_{hp}\|_{\widetilde{H}^{1/2}(\Gamma_\Sigma)}^2 +  \frac{p^2}{h}\|(u-v_{hp})\|_{L^2(\Gamma_C)}^2+  \frac{h}{p^2} \| u-v_{hp} \|^2_{H^1(\Gamma)}   +\| \phi - \psi_{hp} \|^2_{H^{-1/2}(\Gamma)} + \frac{h}{p^2} \| \phi-\psi_{hp} \|^2_{L^2(\Gamma)}\\
& \quad    +  \frac{h}{p^2} \| (\mu_{kq}-\lambda) \|_{L^2(\Gamma_C)}^2 + \frac{h}{p^2} \| (\widetilde{W}-W)u\|_{L^2(\Gamma_C)}^2  +   \frac{h}{p^2}\| (\widetilde{K}^\prime-K^\prime) \psi\|_{L^2(\Gamma_C)}^2 +\langle u_n-g, \lambda_n^{kq} - \mu_n \rangle_{\Gamma_C} + \langle u_t, \lambda_t^{kq} - \mu_t \rangle_{\Gamma_C} \\
& \quad + \| \lambda - \mu_{kq} \|_{\widetilde{H}^{-1/2}(\Gamma_C)} \left( \| u_{hp} \|_{\widetilde{H}^{1/2}(\Gamma_\Sigma)} + \|g\|_{H^{1/2}(\Gamma_C)} \right)+  \frac{h}{p^2} \| (\lambda- \mu_{kq}) \|_{L^2(\Gamma_C)} \left( \|u_{hp}\|_{\widetilde{H}^{1/2}(\Gamma_\Sigma)} + \| \phi_{hp}\|_{H^{-1/2}(\Gamma)} +  \|h^{1/2}p^{-1}\lambda_{kq}\|_{L^2(\Gamma_C)}  \right)
 \end{align*}
for all $(v_{hp},\psi_{hp},\mu_{kq}) \in X_{hp} \times Y_{hp} \times M_{kq}^+(\mathcal{F})$ and all $\mu \in M^+(\mathcal{F})$ with a constant $C>0$ independent of $h$, $p$, $k$ and $q$.
\begin{proof}
 Choosing the test functions $-u_{hp}+v_{hp}$ and $-\phi_{hp}+\psi_{hp}$ in Lemma~\ref{lem:GalerkinOrtho} yields together with Lemma~\ref{lem:aprioriErrorLambda} that
\begin{align*}
  & \qquad \alpha_W \| u-u_{hp}\|_{\widetilde{H}^{1/2}(\Gamma_\Sigma)}^2 + \alpha_V \| \phi-\phi_{hp} \|_{H^{-1/2}(\Gamma)}^2 +2^{-1}\| \gamma^{1/2} ( \lambda - \lambda_{kq} ) \|_{L^2(\Gamma_C)}^2 \\
  &\leq B(u-u_{hp},\phi-\phi_{hp};u-u_{hp},\phi-\phi_{hp}) +2^{-1}\| \gamma^{1/2} ( \lambda - \lambda_{kq} ) \|_{L^2(\Gamma_C)}^2 \\
  &= B(u-u_{hp},\phi-\phi_{hp};u-v_{hp},\phi-\psi_{hp}) + \langle \lambda - \lambda_{kq}, u_{hp} - v_{hp} \rangle_{\Gamma_C} +2^{-1}\| \gamma^{1/2} ( \lambda - \lambda_{kq} ) \|_{L^2(\Gamma_C)}^2\\
  & \quad + \langle \gamma (\lambda_{kq} + \widetilde{W}u_{hp} + (\widetilde{K}+1/2)^\prime \phi_{hp}),\widetilde{W}(u_{hp}-v_{hp}) + (\widetilde{K}+1/2)^\prime (\phi_{hp}-\psi_{hp}) \rangle_{\Gamma_C} \\
  & \leq \underbrace{B(u-u_{hp},\phi-\phi_{hp};u-v_{hp},\phi-\psi_{hp})}_{\rm I} + \underbrace{\langle \lambda - \lambda_{kq}, u - v_{hp} \rangle_{\Gamma_C}}_{\rm II} \\
  & \quad + \underbrace{\langle \gamma (\lambda_{kq} + \widetilde{W}u_{hp} + (\widetilde{K}+1/2)^\prime \phi_{hp}),\widetilde{W}(u_{hp}-v_{hp}) + (\widetilde{K}+1/2)^\prime (\phi_{hp}-\psi_{hp}) \rangle_{\Gamma_C}}_{\rm III}\\
& \quad +\langle u_n-g, \lambda_n^{kq} - \mu_n \rangle_{\Gamma_C} + \langle u_t, \lambda_t^{kq} - \mu_t \rangle_{\Gamma_C} + \langle g, \mu_n^{kq} - \lambda_n  \rangle_{\Gamma_C} + \| \lambda - \mu_{kq} \|_{\widetilde{H}^{-1/2}(\Gamma_C)} \| u_{hp} \|_{\widetilde{H}^{1/2}(\Gamma_\Sigma)}\\
 & \quad + C\| \gamma^{1/2}(\lambda- \mu_{kq}) \|_{L^2(\Gamma_C)} \left( \bar{\gamma}^{1/2}\|u_{hp}\|_{\widetilde{H}^{1/2}(\Gamma_\Sigma)} + \bar{\gamma}^{1/2}\| \phi_{hp}\|_{H^{-1/2}(\Gamma)} + \|\gamma^{1/2}\lambda_{kq}\|_{L^2(\Gamma_C)}\right)    + 4^{-1}  \| \gamma^{1/2} (\mu_{kq}-\lambda) \|_{L^2(\Gamma_C)}^2 \\
 & \quad+ 2\left(  C \bar{\gamma} \|u- v_{hp} \|^2_{\widetilde{H}^{1/2}(\Gamma_\Sigma)} + C \bar{\gamma} \|u - u_{hp} \|^2_{\widetilde{H}^{1/2}(\Gamma_\Sigma)}  + C \bar{\gamma} \|\phi - \psi_{hp} \|^2_{H^{-1/2}(\Gamma)} + C \bar{\gamma} \|\phi - \phi_{hp} \|^2_{H^{-1/2}(\Gamma)}\right)\\
 &\quad  + 4 \left( \| \gamma^{1/2} W(u-v_{hp}) \|_{L^2(\Gamma_C)}^2+\| \gamma^{1/2} (K + 1/2)^\prime (\phi-\psi_{hp}) \|_{L^2(\Gamma_C)}^2+\| \gamma^{1/2}  (\widetilde{W}-W)u_{hp} \|_{L^2(\Gamma_C)}^2
  +\| \gamma^{1/2}  (\widetilde{K}^\prime - K^\prime) \phi_{hp} \|_{L^2(\Gamma_C)}^2 \right).
\end{align*}
 By the continuous mapping of the boundary integral operators and Young's inequality we obtain
 \begin{align*}
 {\rm I} \leq \epsilon \| u-u_{hp}\|_{\widetilde{H}^{1/2}(\Gamma_\Sigma)}^2 + \epsilon\| \phi - \phi_{hp} \|^2_{H^{-1/2}(\Gamma)} + \frac{C}{\epsilon}\| u-v_{hp}\|_{\widetilde{H}^{1/2}(\Gamma_\Sigma)}^2+\frac{C}{\epsilon}\| \phi - \psi_{hp} \|^2_{H^{-1/2}(\Gamma)}.
 \end{align*}
There also holds trivially that
\begin{align*}
  {\rm II} \leq \epsilon \| \gamma^{1/2} (\lambda- \lambda_{kq}) \|_{L^2(\Gamma_C)}^2 + \frac{1}{4\epsilon} \| \gamma^{-1/2} (u-v_{hp})\|_{L^2(\Gamma_C)}^2.
\end{align*}
With triangle inequality, $\lambda = -Wu - (K+1/2)^\prime \phi$ a.e.~on $\Gamma_C$ by  \eqref{eq:compactWeakForm1}, Young's inequality, \eqref{eq:scaled_W_error_est} and \eqref{eq:scaled_K_error_est} we obtain
 \begin{align*}
  & \| \gamma^{1/2}(\lambda_{kq} + \widetilde{W}u_{hp} + (\widetilde{K}+1/2)^\prime \phi_{hp}) \|_{L^2(\Gamma_C)}^2 \leq  \left(\| \gamma^{1/2}(\lambda_{kq} -\lambda) \|_{L^2(\Gamma_C)} +  \| \gamma^{1/2}(\lambda + Wu_{hp} + (K+1/2)^\prime \phi_{hp}) \|_{L^2(\Gamma_C)} \right.\\
   & \qquad \left. +\| \gamma^{1/2}(\widetilde{W}-W)u_{hp}\|_{L^2(\Gamma_C)} + \| \gamma^{1/2}(\widetilde{K}^\prime-K^\prime) \phi_{hp} \|_{L^2(\Gamma_C)} \right)^2 \\
  & \quad \leq  5\| \gamma^{1/2}(\lambda_{kq} -\lambda) \|_{L^2(\Gamma_C)}^2 +  10 \| \gamma^{1/2} W(u-v_{hp}) \|_{L^2(\Gamma_C)}^2+ C \bar{\gamma} \|u- v_{hp} \|^2_{\widetilde{H}^{1/2}(\Gamma_\Sigma)} + C \bar{\gamma} \|u - u_{hp} \|^2_{\widetilde{H}^{1/2}(\Gamma_\Sigma)} \\
  & \qquad + 10 \| \gamma^{1/2} (K + 1/2)^\prime (\phi-\psi_{hp}) \|_{L^2(\Gamma_C)}^2 + C \bar{\gamma} \|\phi - \psi_{hp} \|^2_{H^{-1/2}(\Gamma)} + C \bar{\gamma} \|\phi - \phi_{hp} \|^2_{H^{-1/2}(\Gamma)} \\
   & \qquad +5\| \gamma^{1/2}(\widetilde{W}-W)u_{hp}\|_{L^2(\Gamma_C)}^2 + 5\| \gamma^{1/2}(\widetilde{K}^\prime-K^\prime) \phi_{hp} \|_{L^2(\Gamma_C)}^2.
 \end{align*} 
Due to \eqref{eq:theo_W_tilde}-\eqref{eq:theo_K_tilde} there holds
 \begin{align*}
  &\| \gamma^{1/2}\widetilde{W}(u_{hp}-v_{hp})  + \gamma^{1/2} (\widetilde{K}+1/2)^\prime (\phi_{hp}-\psi_{hp}) \|_{L^2(\Gamma_C)}^2 \\
  &\quad \leq \bar{\gamma} C\left( \| u_{hp}-v_{hp}  \|_{\widetilde{H}^{1/2}(\Gamma_\Sigma)} + 
   \|\phi_{hp}-\psi_{hp} \|_{H^{-1/2}(\Gamma)} \right)^2\\
   &\quad \leq \bar{\gamma} C\left( \| u-v_{hp}  \|_{\widetilde{H}^{1/2}(\Gamma_\Sigma)}^2 + \| u-u_{hp}  \|_{\widetilde{H}^{1/2}(\Gamma_\Sigma)}^2 + 
   \|\phi-\psi_{hp} \|_{H^{-1/2}(\Gamma)}^2 + \|\phi-\phi_{hp} \|_{H^{-1/2}(\Gamma)} ^2\right).
 \end{align*}
With these two estimates there holds
 \begin{align*} 
    {\rm III}  & \leq \epsilon \| \gamma^{1/2}(\lambda_{kq} + \widetilde{W}u_{hp} + (\widetilde{K}+1/2)^\prime \phi_{hp}) \|_{L^2(\Gamma_C)}^2 + \frac{1}{4\epsilon}   \| \gamma^{1/2}\widetilde{W}(u_{hp}-v_{hp}) + \gamma^{1/2} (\widetilde{K}+1/2)^\prime (\phi_{hp}-\psi_{hp}) \|_{L^2(\Gamma_C)}^2 \\
   &  \leq 5\epsilon\| \gamma^{1/2}(\lambda_{kq} -\lambda) \|_{L^2(\Gamma_C)}^2+ \bar{\gamma} C (\epsilon+\epsilon^{-1})\| u-u_{hp}  \|_{\widetilde{H}^{1/2}(\Gamma_\Sigma)}^2
  +\bar{\gamma} C (\epsilon+\epsilon^{-1}) \|\phi-\phi_{hp} \|_{H^{-1/2}(\Gamma)} ^2\\
   & \quad  +  10\epsilon \| \gamma^{1/2} W(u-v_{hp}) \|_{L^2(\Gamma_C)}^2+\bar{\gamma} C (\epsilon+\epsilon^{-1})\| u-v_{hp}  \|_{\widetilde{H}^{1/2}(\Gamma_\Sigma)}^2    
  + 10\epsilon \| \gamma^{1/2} (K + 1/2)^\prime (\phi-\psi_{hp}) \|_{L^2(\Gamma_C)}^2 \\
   & \quad +\bar{\gamma} C (\epsilon+\epsilon^{-1})\|\phi-\psi_{hp} \|_{H^{-1/2}(\Gamma)}^2  +5\epsilon\| \gamma^{1/2}(\widetilde{W}-W)u_{hp}\|_{L^2(\Gamma_C)}^2 + 5\epsilon\| \gamma^{1/2}(\widetilde{K}^\prime-K^\prime) \phi_{hp} \|_{L^2(\Gamma_C)}^2.
 \end{align*}
We estimate the approximation of $W$ further by
\begin{align*}
 \| \gamma^{1/2}(\widetilde{W}-W)u_{hp}\|_{L^2(\Gamma_C)} &\leq \| \gamma^{1/2}(\widetilde{W}-W)(u-u_{hp})\|_{L^2(\Gamma_C)} + \| \gamma^{1/2}(\widetilde{W}-W)u\|_{L^2(\Gamma_C)}\\
 &\leq \bar{\gamma}^{1/2} C \frac{h^{1/2}}{p} \| u-u_{hp} \|_{H^1(\Gamma)}  + \| \gamma^{1/2}(\widetilde{W}-W)u\|_{L^2(\Gamma_C)} \\
 &\leq \bar{\gamma}^{1/2} C \frac{h^{1/2}}{p} \| u-v_{hp} \|_{H^1(\Gamma)} +\bar{\gamma}^{1/2} C  \| v_{hp}-u_{hp} \|_{H^{1/2}(\Gamma)} + \| \gamma^{1/2}(\widetilde{W}-W)u\|_{L^2(\Gamma_C)}\\
 &\leq \bar{\gamma}^{1/2} C \left(\frac{h^{1/2}}{p} \| u-v_{hp} \|_{H^1(\Gamma)} +  \| u-v_{hp} \|_{H^{1/2}(\Gamma)} + \| u-u_{hp} \|_{H^{1/2}(\Gamma)} \right)+ \| \gamma^{1/2}(\widetilde{W}-W)u\|_{L^2(\Gamma_C)} 
\end{align*}
to remove the $u_{hp}$ dependency. Analogously, we obtain
\begin{align*}
\| \gamma^{1/2}(\widetilde{K}^\prime-K^\prime) \phi_{hp}) \|_{L^2(\Gamma_C)} \leq \bar{\gamma}^{1/2} C \left(\frac{h^{1/2}}{p} \| \phi-\psi_{hp} \|_{L^2(\Gamma)} +  \| \phi-\psi_{hp} \|_{H^{-1/2}(\Gamma)} +\| \phi-\phi_{hp} \|_{H^{1/2}(\Gamma)} \right) + \| \gamma^{1/2}(\widetilde{K}^\prime-K^\prime) \phi\|_{L^2(\Gamma_C)} .
\end{align*} 
Combining all the estimates, choosing $\epsilon$ sufficiently small, and then $\bar{\gamma}$ sufficiently small such that $\bar{\gamma} \epsilon^{-1}$ is sufficiently small, we obtain after lengthy and non-enlightening calculations that
 \begin{align*}
 C & \left(\| u-u_{hp}\|_{\widetilde{H}^{1/2}(\Gamma_\Sigma)}^2 +  \| \phi-\phi_{hp} \|_{H^{-1/2}(\Gamma)}^2 +\| h^{1/2}p^{-1} ( \lambda - \lambda_{kq} ) \|_{L^2(\Gamma_C)}^2 \right)\\
& \quad \leq \| u-v_{hp}\|_{\widetilde{H}^{1/2}(\Gamma_\Sigma)}^2 +  \frac{p^2}{h}\|u-v_{hp}\|_{L^2(\Gamma_C)}^2+  \frac{h}{p^2} \| u-v_{hp} \|^2_{H^1(\Gamma)}  +   \frac{h}{p^2} \| W(u-v_{hp}) \|_{L^2(\Gamma_C)}^2 +\| \phi - \psi_{hp} \|^2_{H^{-1/2}(\Gamma)} \\
& \qquad  + \frac{h}{p^2} \| \phi-\psi_{hp} \|^2_{L^2(\Gamma)}  +   \frac{h}{p^2} \|  (K + 1/2)^\prime (\phi-\psi_{hp}) \|_{L^2(\Gamma_C)}^2  +  \frac{h}{p^2} \| \mu_{kq}-\lambda \|_{L^2(\Gamma_C)}^2 + \frac{h}{p^2} \| (\widetilde{W}-W)u\|_{L^2(\Gamma_C)}^2\\
   & \qquad     +   \frac{h}{p^2}\| (\widetilde{K}^\prime-K^\prime) \phi\|_{L^2(\Gamma_C)}^2 +\langle u_n-g, \lambda_n^{kq} - \mu_n \rangle_{\Gamma_C} + \langle u_t, \lambda_t^{kq} - \mu_t \rangle_{\Gamma_C} + \langle g, \mu_n^{kq} - \lambda_n  \rangle_{\Gamma_C} \\
& \qquad + \| \lambda - \mu_{kq} \|_{\widetilde{H}^{-1/2}(\Gamma_C)} \| u_{hp} \|_{\widetilde{H}^{1/2}(\Gamma_\Sigma)}+  \frac{h}{p^2} \| \lambda- \mu_{kq} \|_{L^2(\Gamma_C)} \left( \|u_{hp}\|_{\widetilde{H}^{1/2}(\Gamma_\Sigma)} + \| \phi_{hp}\|_{H^{-1/2}(\Gamma)}  +  \|h^{1/2}p^{-1}\lambda_{kq}\|_{L^2(\Gamma_C)} \right).
 \end{align*}
Now, the assertion follows with the mapping properties of the boundary integral operators, i.e.~$W: H^1(\Gamma) \to L^2(\Gamma)$ and $K': L^2(\Gamma) \to L^2(\Gamma)$ continuously, and Cauchy-Schwarz inequality for $\langle g, \mu_n^{kq} - \lambda_n  \rangle_{\Gamma_C}$.
\end{proof}
\end{theorem}

With that theorem at hand we can proof convergence rates under assumed regularity.
\begin{theorem} \label{thm:AprioriConvergencerate}
 Let $(u,\phi,\lambda)$, $(u_{hp},\phi_{hp},\lambda_{kq})$ be the solution to \eqref{eq:compactWeakForm}, \eqref{eq:DiscreteProblem}, respectively. If $u \in H^\beta(\Gamma)$, $u_n|_{\Gamma_C}$, $u_t|_{\Gamma_C}\in C^0(\Gamma_C)$, $\phi \in H^{\beta-1}(\Gamma)$, $\lambda \in H^{\beta-1}(\Gamma) \cap C^0(\Gamma_C)$, $g \in H^\beta(\Gamma_C)\cap C^0(\Gamma_C)$ for $\beta \geq 1$, and if further $\bar{\gamma}$ is sufficiently small, potentially smaller than for Lemma~\ref{lem:coercive} needed, $G_{kq}$ are the Gauss-Legendre quadrature points and thus $\hat{T}=[-1,1]$ in 2d and $\hat{T}=[-1,1]^2$ in 3d and if the mesh size $k$ is sufficiently small, then there holds 
\begin{align*}
C  \left(\| u-u_{hp}\|_{\widetilde{H}^{1/2}(\Gamma_\Sigma)}^2 +  \| \phi-\phi_{hp} \|_{H^{-1/2}(\Gamma)}^2 +\| h^{1/2}p^{-1} ( \lambda - \lambda_{kq} ) \|_{L^2(\Gamma_C)}^2 \right) & \leq \frac{h^{2\beta-1}}{p^{2\beta-2}} + \frac{p}{h^{1/2}} \left(\frac{k}{q}\right)^{\beta} + \frac{p^2}{h } \left(\frac{k}{q}\right)^{2\beta} +  \left(\frac{k}{q}\right)^{\beta-1} \\  & \quad + \frac{h}{p^2} \| (\widetilde{W}-W)u\|_{L^2(\Gamma_C)}^2  +   \frac{h}{p^2}\| (\widetilde{K}^\prime-K^\prime) \phi\|_{L^2(\Gamma_C)}^2
\end{align*}
 with a constant $C>0$ independent of $h$, $p$, $k$ and $q$.
\begin{proof} 
 We estimate the right hand side terms of Theorem~\ref{thm:AprioriError} individually. For the Scott-Zhang or Cl\'ement interpolation operator $\mathfrak{I}_{hp}$ onto $X_{hp}$ we have the a priori error estimates ($v_{hp}:=\mathfrak{I}_{hp}u$)
 \begin{align} \label{int1}
  \| u-v_{hp}\|_{\widetilde{H}^{1/2}(\Gamma_\Sigma)}^2 \leq C \left(\frac{h}{p} \right)^{2\beta - 1} \| u \|_{H^\beta(\Gamma)}^2, \quad
  \frac{p^2}{h}\| u-v_{hp}\|^2_{L^2(\Gamma_C)} \leq C  \frac{h^{2\beta-1}}{p^{2\beta-2}} \| u \|_{H^\beta(\Gamma)}^2, \quad
  \frac{h}{p^2}\|  u-v_{hp}\|^2_{H^1(\Gamma)} \leq   C  \frac{h^{2\beta-1}}{p^{2\beta}} \| u \|_{H^\beta(\Gamma)}^2.
 \end{align}
For $\psi_{hp}$ the $L^2$-projection of $\phi$ we have
\begin{align} \label{int2}
 \| \phi - \psi_{hp} \|_{H^{-1/2}(\Gamma)}^2 \leq C \left(\frac{h}{p}\right)^{2 \beta -1} \| \phi \|_{H^{\beta-1}(\Gamma)}^2, \quad
  \frac{h}{p^2}\| \phi-\psi_{hp}\|_{L^2(\Gamma_C)}^2 \leq   C \frac{h^{2\beta -1}}{p^{2\beta}}\| \phi \|_{H^{\beta-1}(\Gamma)}^2.
\end{align}
If $\lambda \in C^0(\Gamma_C)$, we can use the classical interpolation in the Gauss-Lobatto points denoted by $\mathcal{I}_{kq} \lambda$ to obtain \cite[Thm.~3.4]{bernardi1992polynomial}
\begin{align}
  \| \lambda - \mu_{kq} \|_{\widetilde{H}^{-1/2}(\Gamma_C)} \leq \| \lambda - \mu_{kq} \|_{L^{2}(\Gamma_C)} \leq C \left(\frac{k}{q}\right)^{\beta-1} \| \lambda \|_{H^{\beta-1}(\Gamma_C)}, \quad 
 \frac{h}{p^2} \| \lambda- \mu_{kq} \|_{L^2(\Gamma_C)} \leq C \frac{h}{p^2} \left(\frac{k}{q}\right)^{\beta-1} \| \lambda \|_{H^{\beta-1}(\Gamma_C)}
\end{align}
and
\begin{align*}
\|(u_n-g)-\mathcal{I}_{kq}(u_n-g)\|_{L^2(\Gamma_C)}\leq C \frac{k^\beta}{q^\beta}\|u_n-g\|_{H^{\beta}(\Gamma_C)} \leq C \frac{k^\beta}{q^\beta}\left(\|u_n\|_{H^{\beta}(\Gamma_C)} + \|g\|_{H^{\beta}(\Gamma_C)} \right), \  \|u_t-\mathcal{I}_{kq} u_t)\|_{L^2(\Gamma_C)}\leq C \frac{k^\beta}{q^\beta}\|u_t\|_{H^{\beta}(\Gamma_C)}.
\end{align*}
From the stability estimate \eqref{eq:discrete_Stab_est2} for a fixed $\bar{\gamma}$ we obtain with the additional regularity of the gap function $g$  that
\begin{align*} 
\ \| u_{hp} \|_{\widetilde{H}^{1/2}(\Gamma_\Sigma)}^2 + \| \phi_{hp} \|_{H^{-1/2}(\Gamma)}^2 + \| h^{1/2}p^{-1} \lambda_{kq} \|_{L^2(\Gamma_C)}^2 \leq C \left( \| f\|_{\widetilde{H}^{-1/2}(\Gamma_N)}^2 + \frac{p^2}{h} \frac{k^{2\beta}}{q^{2\beta}}\| g\|_{H^\beta(\Gamma_C)}^2 \right) \leq C\left(1+\frac{p^2}{h} \frac{k^{2\beta}}{q^{2\beta}}\right) .
\end{align*} 
Since $\int_T \lambda_{kq}\mathcal{I}_{kq}u$ is an integration of a polynomial of degree $2q+1$ on the reference element \cite[Lem.~14]{Banz2015Stab} holds verbatim for the three dimensional case. The assumptions of \cite[Lem.~14]{Banz2015Stab} are satisfied as $ \lambda \in C^0(\Gamma_C)$ and the mesh size $k$ is sufficiently small, c.f.~\cite[Rem.~15.2]{Banz2015Stab}. Hence, we obtain with $\mu_n=0$ that
\begin{align*}
  \langle u_n-g, \lambda_n^{kq} - \mu_n \rangle_{\Gamma_C} &= \int _{\Gamma_C} \lambda_n^{kq}(u_n-g-\mathcal{I}_{kq}(u_n-g)) \, ds + \int_{\Gamma_C} \lambda_n^{kq} \mathcal{I}_{kq}(u_n-g) \, ds \\
	& \leq \int _{\Gamma_C} \lambda_n^{kq}(u_n-g-\mathcal{I}_{kq}(u_n-g)) \, ds \\
	&\leq \|h^{1/2}p^{-1}\lambda_n^{kq}\|_{L^2(\Gamma_C)} \left\|h^{-1/2}p \left [(u_n-g) - \mathcal{I}_{kq}(u_n-g)\right ]\right\|_{L^2(\Gamma_C)} \\
		& \leq C\|h^{1/2}p^{-1}\lambda_n^{kq}\|_{L^2(\Gamma_C)} \frac{p}{h^{1/2}} \left( \frac{k}{q} \right)^\beta \left( \|u_n\|_{H^\beta(\Gamma_C)} + \|g\|_{H^\beta(\Gamma_C)} \right) \\
		&\leq	 C\left(1+\frac{p^2}{h} \frac{k^{2\beta}}{q^{2\beta}}\right)^{1/2} \left( \frac{k}{q} \right)^\beta \frac{p}{h^{1/2}}.
	\end{align*}
Following  the proof of \cite[Lem.~14]{Banz2015Stab} we have
\begin{align*}
\langle \lambda_t^{kq} - \mu_t, u_t \rangle_{L^2(\Gamma_C)} &\leq  \|h^{1/2}p^{-1}\lambda^{kq}\|_{L^2(\Gamma_C)} \| h^{-1/2}p(u_t-\mathcal{I}_{kq} u_t)\|_{L^2(\Gamma_C)} + \|\mu_t\|\, \|u_t-\mathcal{I}_{kq} u_t\|_{L^2(\Gamma_C)} \\
&\leq C\|h^{1/2}p^{-1}\lambda^{kq}\|_{L^2(\Gamma_C)} \frac{p}{h^{1/2}} \left( \frac{k}{q} \right)^\beta\|u_t\|_{H^\beta(\Gamma_C)} + C\|\mu_t\|_{L^2(\Gamma_C)}\left( \frac{k}{q} \right)^\beta\|u_t\|_{H^\beta(\Gamma_C)} \\
& \leq C\|h^{1/2}p^{-1}\lambda^{kq}\|_{L^2(\Gamma_C)} \frac{p}{h^{1/2}} \left( \frac{k}{q} \right) ^\beta\|u_t\|_{H^\beta(\Gamma_C)} + C\|\mathcal{F}\|_{L^\infty(\Gamma_C)} \left |\Gamma_C \right |^{1/2}\left( \frac{k}{q} \right)^\beta\|u_t\|_{H^\beta(\Gamma_C)} \\
&\leq C\left(1+\frac{p^2}{h} \frac{k^{2\beta}}{q^{2\beta}}\right)^{1/2} \left( \frac{k}{q} \right)^\beta \frac{p}{h^{1/2}} + C \left( \frac{k}{q} \right)^\beta.
\end{align*}
Hence,
\[ 
 \inf_{\mu \in M^+(\mathcal{F})} \langle u_n-g, \lambda_n^{kq} - \mu_n \rangle_{\Gamma_C} + \langle u_t, \lambda_t^{kq} - \mu_t \rangle_{\Gamma_C} \leq
 C\left(1+\frac{p^2}{h} \frac{k^{2\beta}}{q^{2\beta}}\right)^{1/2} \left( \frac{k}{q} \right)^\beta \frac{p}{h^{1/2}} + C \left( \frac{k}{q} \right)^\beta.
\]
Thus,
\begin{align*}
 C & \left(\| u-u_{hp}\|_{\widetilde{H}^{1/2}(\Gamma_\Sigma)}^2 +  \| \phi-\phi_{hp} \|_{H^{-1/2}(\Gamma)}^2 +\| h^{1/2}p^{-1} ( \lambda - \lambda_{kq} ) \|_{L^2(\Gamma_C)}^2 \right) \\
 & \leq  \frac{h^{2\beta-1}}{p^{2\beta-2}}  
 +  \frac{h^{2 \beta -1}}{p^{2 \beta -1}}   + \frac{h}{p^2} \frac{k^{\beta-1}}{q^{\beta-1}} + \left ( 1 + \frac{p^2}{h} \frac{k^{2\beta}}{q ^{2\beta}}\right)^{1/2} \frac{k^\beta}{q^\beta} \frac{p}{h^{1/2}} + \frac{k^\beta}{q^\beta} + \frac{k^{\beta-1}}{q^{\beta-1}} \left ( 1 + \frac{p^2}{h} \frac{k^{2\beta}}{q ^{2\beta}}\right)^{1/2} + \frac{h}{p^2} \frac{k^{\beta-1}}{q^{\beta-1}} \left ( 1 + \frac{p^2}{h} \frac{k^{2\beta}}{q ^{2\beta}}\right)^{1/2} \\
&\quad +\frac{h}{p^2} \| (\widetilde{W}-W)u\|_{L^2(\Gamma_C)}^2  +   \frac{h}{p^2}\| (\widetilde{K}^\prime-K^\prime) \phi\|_{L^2(\Gamma_C)}^2.
 \end{align*}
Removing the dominated convergence rates we obtain
\begin{align*}
C & \left(\| u-u_{hp}\|_{\widetilde{H}^{1/2}(\Gamma_\Sigma)}^2 +  \| \phi-\phi_{hp} \|_{H^{-1/2}(\Gamma)}^2 +\| h^{1/2}p^{-1} ( \lambda - \lambda_{kq} ) \|_{L^2(\Gamma_C)}^2 \right) \\
& \leq \frac{h^{2\beta-1}}{p^{2\beta-2}} + \frac{p}{h^{1/2}} \left(\frac{k}{q}\right)^{\beta} + \frac{p^2}{h } \left(\frac{k}{q}\right)^{2\beta} +  \left(\frac{k}{q}\right)^{\beta-1}  + \frac{h}{p^2} \| (\widetilde{W}-W)u\|_{L^2(\Gamma_C)}^2  +   \frac{h}{p^2}\| (\widetilde{K}^\prime-K^\prime) \phi\|_{L^2(\Gamma_C)}^2,
\end{align*}
which completes the proof.
\end{proof}
\end{theorem}

\begin{remark}
 \begin{enumerate}
  \item For the desirable case $h=k$, $p=q$ or $p=q+1$ we get the convergence rate
 \begin{align*}
 C  \left(\| u-u_{hp}\|_{\widetilde{H}^{1/2}(\Gamma_\Sigma)}^2 +  \| \phi-\phi_{hp} \|_{H^{-1/2}(\Gamma)}^2 +\| h^{1/2}p^{-1} ( \lambda - \lambda_{kq} ) \|_{L^2(\Gamma_C)}^2 \right) & 
 \leq  \frac{h^{\beta-1}}{p^{\beta-1}} + \frac{h}{p^2} \| (\widetilde{W}-W)u\|_{L^2(\Gamma_C)}^2 \\
 & \quad +   \frac{h}{p^2}\| (\widetilde{K}^\prime-K^\prime) \phi\|_{L^2(\Gamma_C)}^2.
 \end{align*}  
  \item If $\widetilde{W}$ and $\widetilde{K}^\prime$ are as defined in Lemma~\ref{lem:ApproxStabTerme}, the estimate above becomes
 \begin{align*}
 C  \left(\| u-u_{hp}\|_{\widetilde{H}^{1/2}(\Gamma_\Sigma)}^2 +  \| \phi-\phi_{hp} \|_{H^{-1/2}(\Gamma)}^2 +\| h^{1/2}p^{-1} ( \lambda - \lambda_{kq} ) \|_{L^2(\Gamma_C)}^2 \right) & 
 \leq \frac{h^{\beta-1}}{p^{\beta-1}} + \frac{h}{p^2} \left(\frac{H}{Q}\right)^{2\beta-2} \| Wu \|_{H^{\beta-1}(\Gamma)}^2 \\
 & \quad +   \frac{h}{p^2} \left(\frac{H}{Q}\right)^{2\beta-2} \| K^\prime \phi \|_{H^{\beta-1}(\Gamma)}^2.
 \end{align*} 
 \item We have the for contact problems typical loss of a factor $2$ in the convergence rate, and an additional loss of $1/2$ from the estimates $\| \lambda - \mu_{kq} \|_{\widetilde{H}^{-1/2}(\Gamma_C)} \leq \| \lambda - \mu_{kq} \|_{L^{2}(\Gamma_C)}$ and $\inf_{\mu \in M^+(\mathcal{F})} \langle u_n-g, \lambda_n^{kq} - \mu_n \rangle_{\Gamma_C} + \langle u_t, \lambda_t^{kq} - \mu_t \rangle_{\Gamma_C} \leq
 C\left(1+\frac{p^2}{h} \frac{k^{2\beta}}{q^{2\beta}}\right)^{1/2} \left( \frac{k}{q} \right)^\beta \frac{p}{h^{1/2}} + C \left( \frac{k}{q} \right)^\beta$. 
 \end{enumerate}
\end{remark}

As mentioned at the end of Section~\ref{sec:Amixedboundaryintegralformulation}, the weak formulation with the discrete Poincar\'e-Steklov operator $S_{hp}=W+(K+ 1/2)^\prime V_{hp}^{-1}(K+ 1/2)$ has been stabilized in \cite{Banz2015Stab}. In our notation the discrete weak formulation \cite[Eq.~13]{Banz2015Stab} with the parameters $\beta=\eta=0$ used in \cite[p.~223]{Banz2015Stab}, i.e.~the scaling function $\gamma$ is the same to the one here, is to find the pair $(u_{hp}^{(S)},\lambda_{kq}^{(S)}) \in X_{hp} \times M_{kq}^+(\mathcal{F})$ such that
\begin{subequations} \label{eq:DiscreteProblemSteklov}
\begin{alignat}{2}
 \langle S_{hp} u_{hp}^{(S)},v_{hp}\rangle_{\Gamma}  + \langle \lambda_{kq}^{(S)},v_{hp} \rangle_{\Gamma_C} - \langle \gamma ( \lambda_{kq}^{(S)} + S_{hp} u_{hp}^{(S)}), S_{hp} v_{hp} \rangle_{\Gamma_C} &= \langle {f},v_{hp} \rangle_{\Gamma_N} &\quad& \forall v_{hp} \in X_{hp} \label{eq:DiscreteSteklovMixedVarEq1} \\
\langle u_{hp}^{(S)}, \mu_{kq} -\lambda_{kq}^{(S)} \rangle_{\Gamma_C} - \langle \lambda_{kq}^{(S)}+ S_{hp} u_{hp}^{(S)}, \gamma ( \mu_{kq} - \lambda_{kq}^{(S)}) \rangle_{\Gamma_C}  &\leq \langle g,\mu^{kq}_n-\lambda_{kq,n}^{(S)} \rangle_{\Gamma_C} &\quad& \forall \mu_{kq} \in  M_{kq}^+(\mathcal{F}). \label{eq:DiscreteSteklovContactConstraints}
\end{alignat}
\end{subequations}
According to \cite[Thm.~16]{Banz2015Stab} the proven convergence rate is zero as $\beta=\eta=0$. We improve that estimate in Theorem~\ref{thm:improvedConvergenceRates} by combining our a priori error estimate with the following lemma, which gives a bound on the distance of the two discrete solutions. 

\begin{lemma} \label{lem:dist_2_NuMa_Paper}
 For $\widetilde{W}=W$ and $\widetilde{K}^\prime=K^\prime$ and $\bar{\gamma}$ sufficiently small, potentially smaller than needed for Lemma~\ref{lem:coercive}, let $(u,\phi,\lambda)$, $(u_{hp},\phi_{hp},\lambda_{kq})$, $(u_{hp}^{(S)},\lambda_{kq}^{(S)})$  be the solution to \eqref{eq:compactWeakForm}, \eqref{eq:DiscreteProblem}, \eqref{eq:DiscreteProblemSteklov} respectively. Then there exists a constant $C>0$ independent
 of $h$, $p$, $k$, $q$ and $\bar{\gamma}$ such that 
\begin{align}
  C &\left(  \| u_{hp}^{(S)}-u_{hp} \|^2_{\widetilde{H}^{1/2}(\Gamma_\Sigma)}  +  \| \phi_{hp}^{(S)}-\phi_{hp} \|^2_{H^{-1/2}(\Gamma)} +  \| \gamma^{1/2} (\lambda_{kq}^{(S)}-\lambda_{kq}) \|_{L^2(\Gamma_C)}^2  \right) \nonumber \\
 & \quad \leq \bar{\gamma}^{1/2} \| \gamma^{1/2}(\lambda-\lambda_{kq}) \|_{L^2(\Gamma_C)}^2 + (\bar{\gamma}^{1/2} + \bar{\gamma} )  \left( \|u-u_{hp}\|_{\widetilde{H}^{1/2}(\Gamma_\Sigma)}^2 +  \| \phi-\phi_{hp} \|^2_{H^{-1/2}(\Gamma)}  +  \| \phi - \psi_{hp}  \|_{H^{-1/2}(\Gamma)} ^2  \right) \nonumber \\
  & \qquad  +\bar{\gamma} \left( \frac{h}{p^2} \| u-v_{hp}\|_{H^1(\Gamma)}^2 + \| u-v_{hp} \|^2_{\widetilde{H}^{1/2}(\Gamma_\Sigma)} + \frac{h}{p^2} \| \phi-\psi_{hp}  \|_{L^2(\Gamma)}^2  \right)
\end{align}
for all $(v_{hp},\psi_{hp}) \in X_{hp} \times Y_{hp}$ and with $\phi_{hp}^{(S)} =  V_{hp}^{-1} (K+1/2) u_{hp}^{(S)}$. 
 \begin{proof}
 Recall that $\langle S_{hp} u_{hp}^{(S)},v_{hp}\rangle_{\Gamma}$ is in fact
\begin{align}
 \langle W u_{hp}^{(S)} + (K+1/2)^\prime \phi_{hp}^{(S)},v_{hp}\rangle_{\Gamma} \quad \text{with} \quad  \langle V \phi_{hp}^{(S)},\psi_{hp}\rangle_{\Gamma} - \langle (K+1/2) u_{hp}^{(S)},\psi_{hp}\rangle_{\Gamma} =0 \quad \forall \psi_{hp} \in Y_{hp}. \label{representation-Steklov}
\end{align}
 Subtracting \eqref{eq:DiscreteMixedVarEq1} from \eqref{eq:DiscreteSteklovMixedVarEq1} each with the test functions $v_{hp}=u_{hp}^{(S)}-u_{hp}$ and $\psi_{hp}=\phi_{hp}^{(S)}-\phi_{hp}$, and using the representation \eqref{representation-Steklov} yields
 \begin{align*}
 0 &= \langle W(u_{hp}^{(S)}-u_{hp}),u_{hp}^{(S)}-u_{hp} \rangle_\Gamma + \langle V(\phi_{hp}^{(S)}-\phi_{hp}),\phi_{hp}^{(S)}-\phi_{hp} \rangle_\Gamma + \langle \lambda_{kq}^{(S)} - \lambda_{kq}, u_{hp}^{(S)}-u_{hp} \rangle_{\Gamma_C} \\
 & \quad  - \langle \gamma ( \lambda_{kq}^{(S)} + S_{hp} u_{hp}^{(S)}), S_{hp} (u_{hp}^{(S)}-u_{hp})\rangle_{\Gamma_C} 
  + \langle \gamma ( \lambda_{kq} + Wu_{hp} + (K + 1/2)^\prime \phi_{hp}),  W(u_{hp}^{(S)}-u_{hp}) \\
  & \quad + (K + 1/2)^\prime (\phi_{hp}^{(S)}-\phi_{hp}) \rangle_{\Gamma_C} .
 \end{align*}
Analogously, adding \eqref{eq:DiscreteContactConstraints} and \eqref{eq:DiscreteSteklovContactConstraints} with the test function $\lambda_{kq}^{(S)}$, $\lambda_{kq}$, respectively, yields
\begin{align*}
 \langle u_{hp}^{(S)}-u_{hp},\lambda_{kq} - \lambda_{kq}^{(S)} \rangle_{\Gamma_C} - \langle \lambda_{kq}^{(S)}-\lambda_{kq}+ S_{hp} u_{hp}^{(S)} - Wu_{hp} - (K + 1/2)^\prime \phi_{hp}, \gamma ( \lambda_{kq} - \lambda_{kq}^{(S)}) \rangle_{\Gamma_C} \leq 0
\end{align*}
which we write as
\begin{align*}
 \| \gamma^{1/2} (\lambda_{kq}^{(S)}-\lambda_{kq}) \|_{L^2(\Gamma_C)}^2 +  \langle  S_{hp} u_{hp}^{(S)} - Wu_{hp} - (K + 1/2)^\prime \phi_{hp}, \gamma ( \lambda_{kq}^{(S)} - \lambda_{kq}) \rangle_{\Gamma_C} \leq \langle u_{hp}^{(S)}-u_{hp}, \lambda_{kq}^{(S)} -\lambda_{kq} \rangle_{\Gamma_C}.
\end{align*}
Hence, we obtain with the coercivity of $W$ and $V$ that
\begin{align}
 &\alpha_W \| u_{hp}^{(S)}-u_{hp} \|^2_{\widetilde{H}^{1/2}(\Gamma_\Sigma)} + \alpha_V \| \phi_{hp}^{(S)}-\phi_{hp} \|^2_{H^{-1/2}(\Gamma)} + \| \gamma^{1/2} (\lambda_{kq}^{(S)}-\lambda_{kq}) \|_{L^2(\Gamma_C)}^2 \nonumber \\
 &\quad \leq \langle \gamma ( \lambda_{kq}^{(S)} + S_{hp} u_{hp}^{(S)}), S_{hp} (u_{hp}^{(S)}-u_{hp})\rangle_{\Gamma_C} 
 -\langle \gamma ( \lambda_{kq} + Wu_{hp} + (K + 1/2)^\prime \phi_{hp}),  W(u_{hp}^{(S)}-u_{hp}) + (K + 1/2)^\prime (\phi_{hp}^{(S)}-\phi_{hp}) \rangle_{\Gamma_C} \nonumber \\
 & \qquad - \langle  S_{hp} u_{hp}^{(S)} - Wu_{hp} - (K + 1/2)^\prime \phi_{hp}, \gamma ( \lambda_{kq}^{(S)} - \lambda_{kq}) \rangle_{\Gamma_C}. \label{eq:dist_2_NuMa_Paper1}
\end{align}
Using the representation formula \eqref{representation-Steklov}, the definition of $S_{hp}$ and with $\phi_{hp}^{(S)}=V_{hp}^{-1}(K+1/2) u_{hp}^{(S)}$ we have
\begin{align*}
 &  \langle \gamma ( \lambda_{kq}^{(S)} + S_{hp} u_{hp}^{(S)}), S_{hp} (u_{hp}^{(S)}-u_{hp})\rangle_{\Gamma_C}  \\
 &\quad =\langle \gamma ( \lambda_{kq}^{(S)} + W u_{hp}^{(S)} + (K+1/2)^\prime  \phi_{hp}^{(S)} ), W(u_{hp}^{(S)}-u_{hp}) + (K+1/2)^\prime V_{hp}^{-1}(K+1/2) (u_{hp}^{(S)}-u_{hp})\rangle_{\Gamma_C}  \\
&\quad = \langle \gamma ( \lambda_{kq}^{(S)} + W u_{hp}^{(S)} + (K+1/2)^\prime  \phi_{hp}^{(S)} ), W(u_{hp}^{(S)}-u_{hp}) + (K+1/2)^\prime (\phi_{hp}^{(S)}-V_{hp}^{-1}(K+1/2) u_{hp})\rangle_{\Gamma_C}
\end{align*}
and, hence, 
using analogon of \eqref{eq:theo_W_tilde}-\eqref{eq:theo_K_tilde} for the operators $W$ and $(K+1/2)'$
\begin{align*}
 &\langle \gamma ( \lambda_{kq}^{(S)} + S_{hp} u_{hp}^{(S)}), S_{hp} (u_{hp}^{(S)}-u_{hp})\rangle_{\Gamma_C} 
 -\langle \gamma ( \lambda_{kq} + Wu_{hp} + (K + 1/2)^\prime \phi_{hp}),  W(u_{hp}^{(S)}-u_{hp}) + (K + 1/2)^\prime (\phi_{hp}^{(S)}-\phi_{hp}) \rangle_{\Gamma_C}\\
 &\quad = \langle \gamma ( \lambda_{kq}^{(S)}-\lambda_{kq} + W (u_{hp}^{(S)} - u_{hp}) + (K+1/2)^\prime ( \phi_{hp}^{(S)}-\phi_{hp}) ), W(u_{hp}^{(S)}-u_{hp}) + (K + 1/2)^\prime (\phi_{hp}^{(S)}-\phi_{hp})\rangle_{\Gamma_C} \\
 & \qquad +  \langle \gamma ( \lambda_{kq}^{(S)} + W u_{hp}^{(S)} + (K+1/2)^\prime  \phi_{hp}^{(S)} ), (K+1/2)^\prime (\phi_{hp} - V_{hp}^{-1}(K+1/2) u_{hp})\rangle_{\Gamma_C} \\
 & \quad \leq \bar{\gamma}C  \| u_{hp}^{(S)}-u_{hp} \|^2_{\widetilde{H}^{1/2}(\Gamma_\Sigma)} + \bar{\gamma}C  \| \phi_{hp}^{(S)}-\phi_{hp} \|_{H^{-1/2}(\Gamma)} +\langle \gamma ( \lambda_{kq}^{(S)}-\lambda_{kq} ), W(u_{hp}^{(S)}-u_{hp}) + (K + 1/2)^\prime (\phi_{hp}^{(S)}-\phi_{hp})\rangle_{\Gamma_C} \\
 & \qquad +  \langle \gamma ( \lambda_{kq}^{(S)} + W u_{hp}^{(S)} + (K+1/2)^\prime  \phi_{hp}^{(S)} ), (K+1/2)^\prime (\phi_{hp} - V_{hp}^{-1}(K+1/2) u_{hp})\rangle_{\Gamma_C}.
\end{align*}
Consequently, using \eqref{representation-Steklov}, namely, $\langle S_{hp}u_{hp}^{(S)}, \gamma(\lambda_k^{(S)}-\lambda_{kq})\rangle _\Gamma = \langle Wu_{hp}^{(S)}+(K+1/2)'\phi_{hp}^{(S)}, \gamma(\lambda_k^{(S)}-\lambda_{kq})\rangle _\Gamma$, the estimation \eqref{eq:dist_2_NuMa_Paper1} simplifies to
\begin{align}
 (\alpha_W-\bar{\gamma}C) & \| u_{hp}^{(S)}-u_{hp} \|^2_{\widetilde{H}^{1/2}(\Gamma_\Sigma)}  + (\alpha_V-\bar{\gamma}C) \| \phi_{hp}^{(S)}-\phi_{hp} \|^2_{H^{-1/2}(\Gamma)} + \| \gamma^{1/2} (\lambda_{kq}^{(S)}-\lambda_{kq}) \|_{L^2(\Gamma_C)}^2 \nonumber \\
 & \quad  \leq  \langle \gamma ( \lambda_{kq}^{(S)} + W u_{hp}^{(S)} + (K+1/2)^\prime  \phi_{hp}^{(S)} ), (K+1/2)^\prime (\phi_{hp} - V_{hp}^{-1}(K+1/2) u_{hp})\rangle_{\Gamma_C}. \label{eq:dist_2_NuMa_Paper2}
\end{align}
Setting for now $b:=(K+1/2)^\prime (\phi_{hp} - V_{hp}^{-1}(K+1/2) ) u_{hp}$ and adding $\lambda +Wu + (K+1/2)'\phi=0$ a.e.~on $\Gamma_C$ in the third line below, we obtain with  analogon of \eqref{eq:theo_W_tilde}-\eqref{eq:theo_K_tilde} for the operators $W$ and $(K+ 1/2)^\prime$ that
 \begin{align*}
  &\langle \gamma ( \lambda_{kq}^{(S)} + W u_{hp}^{(S)} + (K+1/2)^\prime  \phi_{hp}^{(S)} ), (K+1/2)^\prime (\phi_{hp} - V_{hp}^{-1}(K+1/2) u_{hp})\rangle_{\Gamma_C}\\
  &\quad =\langle \gamma^{1/2} ( \lambda_{kq}^{(S)}-\lambda_{kq} + W (u_{hp}^{(S)}-u_{hp}) + (K+1/2)^\prime ( \phi_{hp}^{(S)}- \phi_{hp})), \gamma^{1/2} b \rangle_{\Gamma_C} \\
  & \qquad + \langle \gamma^{1/2} ( \lambda_{kq}-\lambda + W (u_{hp}-u) + (K+1/2)^\prime  (\phi_{hp}-\phi) ), \gamma^{1/2} b\rangle_{\Gamma_C} \\
  &\quad \leq \left( \| \gamma^{1/2}(\lambda_{kq}^{(S)}-\lambda_{kq}) \|_{L^2(\Gamma_C)}  + \bar{\gamma}^{1/2}C \| u_{hp}^{(S)}-u_{hp} \|_{\widetilde{H}^{1/2}(\Gamma_\Sigma)} + \bar{\gamma}^{1/2}C \| \phi_{hp}^{(S)}-\phi_{hp} \|_{H^{-1/2}(\Gamma)} + \| \gamma^{1/2}(\lambda-\lambda_{kq}) \|_{L^2(\Gamma_C)}  \right. \\
  & \qquad \left. + \| \gamma^{1/2}W(u-u_{hp}) \|_{L^2(\Gamma_C)} + \| \gamma^{1/2} (K+1/2)^\prime(\phi-\phi_{hp}) \|_{L^2(\Gamma_C)}\right)  \cdot \| \gamma^{1/2} b \|_{L^2(\Gamma_C)}.
 \end{align*}
We further obtain with $v_{hp} \in X_{hp}$ arbitrary and using  analogon of \eqref{eq:theo_W_tilde} for the operator $W$ that
 \begin{align*}
  \| \gamma^{1/2}W(u-u_{hp}) \|_{L^2(\Gamma_C)} &\leq \| \gamma^{1/2}W(u-v_{hp}) \|_{L^2(\Gamma_C)} + \| \gamma^{1/2}W(v_{hp}-u_{hp}) \|_{L^2(\Gamma_C)} \\
  & \leq \bar{\gamma}^{1/2} \frac{h^{1/2}}{p} \| u-v_{hp}\|_{H^1(\Gamma)} + \bar{\gamma}^{1/2}C \| v_{hp}-u_{hp} \|_{\widetilde{H}^{1/2}(\Gamma_\Sigma)} \\
  &\leq \bar{\gamma}^{1/2} \frac{h^{1/2}}{p} \| u-v_{hp}\|_{H^1(\Gamma)} + \bar{\gamma}^{1/2}C \| u-v_{hp} \|_{\widetilde{H}^{1/2}(\Gamma_\Sigma)} + \bar{\gamma}^{1/2}C \| u-u_{hp} \|_{\widetilde{H}^{1/2}(\Gamma_\Sigma)}
 \end{align*}
and, analogously with $\psi_{hp} \in Y_{hp}$ arbitrary, that
\begin{align*}
 \| \gamma^{1/2} (K+1/2)^\prime(\phi-\phi_{hp} ) \|_{L^2(\Gamma_C)} \leq \bar{\gamma}^{1/2} \frac{h^{1/2}}{p} \| \phi-\psi_{hp}  \|_{L^2(\Gamma)} + \bar{\gamma}^{1/2}C \| \phi-\psi_{hp} \|_{H^{-1/2}(\Gamma)} + \bar{\gamma}^{1/2}C \| \phi-\phi_{hp} \|_{H^{-1/2}(\Gamma)}.
\end{align*}
Since $V_{hp}^{-1}(K+1/2) u_{hp} \in Y_{hp}$ we have from the $(K+1/2)^\prime$ analogon of \eqref{eq:theo_K_tilde}, $\phi = V^{-1}(K+1/2)u$ and bijectivity of $V$ that
\begin{align*}
 \| \gamma^{1/2} b \|_{L^2(\Gamma_C)} &\leq \bar{\gamma}^{1/2} C \| \phi_{hp} - V_{hp}^{-1}(K+1/2) u_{hp} \|_{H^{-1/2}(\Gamma)} \\
 &= \bar{\gamma}^{1/2} C \| \phi_{hp} - \phi - (V_{hp}^{-1}-V^{-1})(K+1/2) u_{hp} + V^{-1}(K+1/2) (u-u_{hp}) \|_{H^{-1/2}(\Gamma)} \\
 &\leq \bar{\gamma}^{1/2} C \left( \| \phi_{hp} - \phi\|_{H^{-1/2}(\Gamma)} + \|u-u_{hp}\|_{\widetilde{H}^{1/2}(\Gamma_\Sigma)} + \| (V_{hp}^{-1}-V^{-1})(K+1/2) u_{hp}  \|_{H^{-1/2}(\Gamma)} \right).
\end{align*}
Hence, \eqref{eq:dist_2_NuMa_Paper2} becomes 
\begin{align*}
  (\alpha_W &-\bar{\gamma}C)  \| u_{hp}^{(S)}-u_{hp} \|^2_{\widetilde{H}^{1/2}(\Gamma_\Sigma)}  + (\alpha_V-\bar{\gamma}C) \| \phi_{hp}^{(S)}-\phi_{hp} \|_{H^{-1/2}(\Gamma)}^2 + 2^{-1} \| \gamma^{1/2} (\lambda_{kq}^{(S)}-\lambda_{kq}) \|_{L^2(\Gamma_C)}^2 \\
 & \leq \bar{\gamma}^{1/2} \| \gamma^{1/2}(\lambda-\lambda_{kq}) \|_{L^2(\Gamma_C)}^2 + (\bar{\gamma}^{1/2} + \bar{\gamma} ) C \left( \|u-u_{hp}\|_{\widetilde{H}^{1/2}(\Gamma_\Sigma)}^2 +  \| \phi-\phi_{hp} \|^2_{H^{-1/2}(\Gamma)}  +  \| (V_{hp}^{-1}-V^{-1})(K+1/2) u_{hp}  \|_{H^{-1/2}(\Gamma)}^2  \right) \\
  & \quad  +\bar{\gamma} C \left( \frac{h}{p^2} \| u-v_{hp}\|_{H^1(\Gamma)}^2 + \| u-v_{hp} \|^2_{\widetilde{H}^{1/2}(\Gamma_\Sigma)} + \frac{h}{p^2} \| \phi-\psi_{hp}  \|_{L^2(\Gamma)}^2 +  \| \phi-\psi_{hp} \|_{H^{-1/2}(\Gamma)}^2 \right).
\end{align*}
It is well known that
\begin{align*}
 C_V^{-1} \alpha_V \| (V_{hp}^{-1}-V^{-1})(K+1/2) u_{hp}  \|_{H^{-1/2}(\Gamma)} &\leq  \| V^{-1}(K+1/2) u_{hp} - \psi_{hp}  \|_{H^{-1/2}(\Gamma)} \\
 &\leq \| V^{-1}(K+1/2) u - \psi_{hp}  \|_{H^{-1/2}(\Gamma)} + \| V^{-1}(K+1/2) (u- u_{hp}) \|_{H^{-1/2}(\Gamma)} \\
 &\leq \| V^{-1}(K+1/2) u - \psi_{hp}  \|_{H^{-1/2}(\Gamma)} + \| u- u_{hp} \|_{\widetilde{H}^{1/2}(\Gamma_\Sigma)} 
\end{align*}
for all $\psi_{hp} \in Y_{hp}$. Moreover, $V^{-1}(K+1/2)$ is the non-symmetric Steklov operator and thus $V^{-1}(K+1/2) u = \phi$ which completes the proof.
 \end{proof}
\end{lemma}

\begin{theorem} \label{thm:improvedConvergenceRates}
Under the assumptions of Theorem~\ref{thm:AprioriConvergencerate} and Lemma~\ref{lem:dist_2_NuMa_Paper}, there holds
\begin{align*}
  C\left( \| u- u_{hp}^{(S)} \|^2_{\widetilde{H}^{1/2}(\Gamma_\Sigma)}  +  \| \phi- \phi_{hp}^{(S)} \|_{H^{-1/2}(\Gamma)}^2 +  \| h^{1/2}p^{-1} (\lambda - \lambda_{kq}^{(S)}) \|_{L^2(\Gamma_C)}^2 \right) \leq  \frac{h^{2\beta-1}}{p^{2\beta-2}} + \frac{p}{h^{1/2}} \left(\frac{k}{q}\right)^{\beta} + \frac{p^2}{h } \left(\frac{k}{q}\right)^{2\beta} +  \left(\frac{k}{q}\right)^{\beta-1} 
\end{align*}
with a constant $C>0$ independent of $h$, $p$, $k$ and $q$.
\begin{proof}
 By triangle inequality and Lemma~\ref{lem:dist_2_NuMa_Paper} we have
 \begin{align*}
   \| u- u_{hp}^{(S)} \|^2_{\widetilde{H}^{1/2}(\Gamma_\Sigma)}  &+  \| \phi- \phi_{hp}^{(S)} \|_{H^{-1/2}(\Gamma)}^2 +  \| h^{1/2}p^{-1} (\lambda - \lambda_{kq}^{(S)}) \|_{L^2(\Gamma_C)}^2 \\
   &\leq   2 \| u- u_{hp} \|^2_{\widetilde{H}^{1/2}(\Gamma_\Sigma)}  +  2\| \phi- \phi_{hp} \|_{H^{-1/2}(\Gamma)}^2 +  2\| h^{1/2}p^{-1} (\lambda - \lambda_{kq}) \|_{L^2(\Gamma_C)}^2+2 \| u_{hp}^{(S)}- u_{hp} \|^2_{\widetilde{H}^{1/2}(\Gamma_\Sigma)} \\
   & \quad  +  2\| \phi_{hp}^{(S)}- \phi_{hp} \|_{H^{-1/2}(\Gamma)}^2 +  2\| h^{1/2}p^{-1} (\lambda_{kq}^{(S)} - \lambda_{kq}) \|_{L^2(\Gamma_C)}^2 \\
   &\leq  C \left( \| u- u_{hp} \|^2_{\widetilde{H}^{1/2}(\Gamma_\Sigma)}  +  \| \phi- \phi_{hp} \|_{H^{-1/2}(\Gamma)}^2 +  \| h^{1/2}p^{-1} (\lambda - \lambda_{kq}) \|_{L^2(\Gamma_C)}^2+ \frac{h}{p^2} \| u-v_{hp}\|_{H^1(\Gamma)}^2 \right.\\
   & \quad   \left.  + \| u-v_{hp} \|^2_{\widetilde{H}^{1/2}(\Gamma_\Sigma)} + \frac{h}{p^2} \| \phi-\psi_{hp}  \|_{L^2(\Gamma)}^2 +  \| \phi-\psi_{hp} \|_{H^{-1/2}(\Gamma)}^2 \right) \\
   &\leq  C \left( \| u- u_{hp} \|^2_{\widetilde{H}^{1/2}(\Gamma_\Sigma)}  +  \| \phi- \phi_{hp} \|_{H^{-1/2}(\Gamma)}^2 +  \| h^{1/2}p^{-1} (\lambda - \lambda_{kq}) \|_{L^2(\Gamma_C)}^2+  \left(\frac{h}{p} \right)^{2\beta - 1} \left( \| u \|_{H^\beta(\Gamma)}^2 + \| \phi \|_{H^{\beta-1}(\Gamma)}^2 \right) \right)   
 \end{align*}
where the last inequality follows from the approximation properties \eqref{int1} and \eqref{int2} of the spaces $X_{hp}$ and $Y_{hp}$. 
The assertion follows with Theorem~\ref{thm:AprioriConvergencerate}. 
\end{proof}
\end{theorem}

\section{A posteriori error estimate}
\label{sec:Aposteriorierrorestimate}
As the stabilization term is based on the residual, using a residual based a posteriori error estimate is suggesting itself. 
To prove the actual a posteriori error estimate Theorem~\ref{thm:aposterioriErrorEst} we need several intermediate results. In almost all estimates the following local error contributions appear:
\begin{subequations} \label{eq:local_error_contributions}
\begin{gather}
 \eta_{N,T}^2 := \frac{h_T}{p_T} \| f -Wu_{hp} - (K+1/2)^\prime \phi_{hp} \|^2_{L^2(T)}, \quad
 \eta_{C,T}^2 := \frac{h_T}{p_T} \| \lambda_{kq} +Wu_{hp} + (K+1/2)^\prime \phi_{hp} \|_{L^2(T)}^2 \\
 \eta_{V,T}^2 := h_T \| \nabla (V\phi_{hp} - (K+1/2)u_{hp}) \|^2_{L^2(T)}, \quad
 \eta_{W,T}^2 := \frac{h_T}{p_T^2} \| (\widetilde{W}-W)u_{hp}\|_{L^2(T)}^2 , \quad
 \eta_{K,T}^2 := \frac{h_T}{p_T^2} \|(\widetilde{K} -K)^\prime \phi_{hp} \|_{L^2(T)}^2
\end{gather}
\end{subequations}

As the weak inequality constraint \eqref{eq:compactWeakForm2} is independent of the stabilization term, and the discrete inequality constraint \eqref{eq:DiscreteContactConstraints} is not needed for the proof of the a posteriori error estimate there holds analogously to \cite[Lem.~19]{Banz2015Stab}  that:

\begin{lemma} \label{lem:aposterioriContactPart}
Let $(u,\phi,\lambda)$, $(u_{hp},\phi_{hp},\lambda_{kq})$ be the solution to \eqref{eq:compactWeakForm}, \eqref{eq:DiscreteProblem}, respectively. Then there holds
\begin{align*}
\langle \lambda-\lambda_{kq}, u_{hp}-u\rangle_{\Gamma_C} \leq &
\langle (\lambda^{kq}_n)^+,(g-u^{hp}_n)^+\rangle_{\Gamma_C}+ \|\lambda_{kq} -\lambda \|_{\widetilde{H}^{-1/2}(\Gamma_C)} \| (g-u^{hp}_n)^-\|_{H^{1/2}(\Gamma_C)} + \|(\lambda^{kq}_n)^-\|_{\widetilde{H}^{-1/2}(\Gamma_C)} \|u_{hp}-u \|_{\widetilde{H}^{1/2}(\Gamma_\Sigma)} \\
&+ \| (|\lambda^{kq}_t|-\mathcal{F} )^+ \|_{\widetilde{H}^{-1/2}(\Gamma_C)} \| u-u_{hp} \|_{\widetilde{H}^{1/2}(\Gamma_\Sigma)}  -\langle (|\lambda^{kq}_t|-\mathcal{F} )^-, |u^{hp}_t| \rangle_{\Gamma_C} -\langle  \lambda^{kq}_t,u^{hp}_t \rangle_{\Gamma_C}
 +\langle |\lambda^{kq}_t| , |u^{hp}_t| \rangle_{\Gamma_C} 
\end{align*}
where $v^+=\max\{0,v\}$ and $v^-=\min\{0,v\}$, i.e.~$v=v^+ +v^-$. 
\end{lemma}

\begin{lemma} \label{lem:aposteriori_lambda}
Let $(u,\phi,\lambda)$, $(u_{hp},\phi_{hp},\lambda_{kq})$ be the solution to \eqref{eq:compactWeakForm}, \eqref{eq:DiscreteProblem}, respectively. Then there holds
\begin{align*}
C\| \lambda-\lambda_{kq} \|_{\widetilde{H}^{-1/2}(\Gamma_C)}^2 & \leq   \sum_{T \in \mathcal{T}_h \cap \Gamma_N} \eta_{N,T}^2 + \sum_{T \in \mathcal{T}_h \cap \Gamma_C} \left( (1+\bar{\gamma}^2) \eta_{C,T}^2 + \bar{\gamma}^2\eta_{W,T}^2 + \bar{\gamma}^2\eta_{K,T}^2 \right)  +  \|u-u_{hp}\|_{\widetilde{H}^{1/2}(\Gamma_\Sigma)}^2  + \|\phi-\phi_{hp}\|_{{H}^{-1/2}(\Gamma)}^2 
 \end{align*}
with a constant $C>0$ independent of $h$, $p$, $k$, $q$ and $\bar{\gamma}$ and with the local error contributions defined in \eqref{eq:local_error_contributions}.
 \begin{proof}
  Let $v \in \widetilde{H}^{1/2}(\Gamma_\Sigma)$ be arbitrary, from Lemma~\ref{lem:GalerkinOrtho} and \eqref{eq:compactWeakForm1} with $\psi=\psi_{hp}=0$ we obtain
  \begin{align*}
   \langle \lambda - \lambda_{kq} , v \rangle_{\Gamma_C} &= \langle \lambda - \lambda_{kq} , v  - v_{hp}\rangle_{\Gamma_C} + \langle \lambda - \lambda_{kq} ,  v_{hp}\rangle_{\Gamma_C} -   B(u-u_{hp},\phi - \phi_{hp};v_{hp},0) \\
    & \quad - \langle \gamma ( \lambda_{kq} + \widetilde{W}u_{hp} + (\widetilde{K} + 1/2)^\prime \phi_{hp}),  \widetilde{W}v_{hp} \rangle_{\Gamma_C}  + \langle f , v-v_{hp} \rangle_{\Gamma_N} - B(u,\phi;v-v_{hp},0) - \langle \lambda , v - v_{hp} \rangle_{\Gamma_C} \\
&=  \langle f -Wu_{hp} - (K+1/2)^\prime \phi_{hp}, v-v_{hp} \rangle_{\Gamma_N} + \langle -\lambda_{kq} -Wu_{hp} - (K+1/2)^\prime \phi_{hp}, v-v_{hp} \rangle_{\Gamma_C}   \\
& \quad -   B(u-u_{hp},\phi - \phi_{hp};v,0) - \langle \gamma ( \lambda_{kq} + \widetilde{W}u_{hp} + (\widetilde{K} + 1/2)^\prime \phi_{hp}),  \widetilde{W}v_{hp} \rangle_{\Gamma_C} \\
&\leq \sum_{T \in \mathcal{T}_h \cap \Gamma_N} \| f -Wu_{hp} - (K+1/2)^\prime \phi_{hp} \|_{L^2(T)} \| v-v_{hp}\|_{L^2(T)} + 
\sum_{T \in \mathcal{T}_h \cap \Gamma_C} \| \lambda_{kq} +Wu_{hp} + (K+1/2)^\prime \phi_{hp} \|_{L^2(T)} \\
& \quad \cdot \| v-v_{hp}\|_{L^2(T)}  + \| W(u-u_{hp}) \|_{H^{-1/2}(\Gamma)} \|v\|_{H^{1/2}(\Gamma)} + \| (K+1/2)^\prime (\phi-\phi_{hp})\|_{H^{-1/2}(\Gamma)} \|v\|_{H^{1/2}(\Gamma)} \\
& \quad - \langle \gamma ( \lambda_{kq} + \widetilde{W}u_{hp} + (\widetilde{K} + 1/2)^\prime \phi_{hp}),  \widetilde{W}v_{hp} \rangle_{\Gamma_C} \\
& \leq C \left( \sum_{T \in \mathcal{T}_h \cap \Gamma_N} \frac{h_T}{p_T} \| f -Wu_{hp} - (K+1/2)^\prime \phi_{hp} \|^2_{L^2(T)} + \sum_{T \in \mathcal{T}_h \cap \Gamma_C} \frac{h_T}{p_T} \| \lambda_{kq} +Wu_{hp} + (K+1/2)^\prime \phi_{hp} \|_{L^2(T)}^2 \right. \\
& \quad  \left. +  \|u-u_{hp}\|_{\widetilde{H}^{1/2}(\Gamma_\Sigma)}^2  + \|\phi-\phi_{hp}\|^2_{H^{-1/2}(\Gamma)}  \right)^{1/2} \|v\|_{\tilde{H}^{1/2}(\Gamma_\Sigma)} - \langle \gamma ( \lambda_{kq} + \widetilde{W}u_{hp} + (\widetilde{K} + 1/2)^\prime \phi_{hp}),  \widetilde{W}v_{hp} \rangle_{\Gamma_C}
  \end{align*}
with $v_{hp}\in X_{hp}$ the Cl\'{e}ment interpolation of $v \in \widetilde{H}^{1/2}(\Gamma_\Sigma)$ with the property $\|v-v_{hp}\|_{L^2(T)}\leq C \left ( \frac{h_T}{p_T}\right)^{1/2}\|v\|_{{H}^{1/2}(T)}$, and using the inequality $a+b \leq \sqrt{2a^2+2b^2}$.\\
The last term which comes from the stabilization can be bounded with \eqref{eq:theo_W_tilde} by the $H^{1/2}$-stability of the interpolation
\begin{align}
  & \qquad \langle \gamma ( \lambda_{kq} + \widetilde{W}u_{hp} + (\widetilde{K} + 1/2)^\prime \phi_{hp}),  \widetilde{W}v_{hp} \rangle_{\Gamma_C} \nonumber \\
  &\leq    \| \gamma^{1/2} ( \lambda_{kq} + Wu_{hp} + (K + 1/2)^\prime \phi_{hp} + (\widetilde{W}-W)u_{hp} + (\widetilde{K} -K)^\prime \phi_{hp} \|_{L^2(\Gamma_C)} \|\gamma^{1/2} \widetilde{W}v_{hp} \|_{L^2(\Gamma_C)} \nonumber \\
  &\leq \bar{\gamma}^{1/2} C \left( \| \gamma^{1/2} ( \lambda_{kq} + Wu_{hp} + (K + 1/2)^\prime \phi_{hp} )\|_{L^2(\Gamma_C)} + \| \gamma^{1/2}(\widetilde{W}-W)u_{hp}\|_{L^2(\Gamma_C)} + \| \gamma^{1/2}(\widetilde{K} -K)^\prime \phi_{hp} \|_{L^2(\Gamma_C)} \right) \| v_{hp} \|_{\widetilde{H}^{1/2}(\Gamma_\Sigma)} \nonumber \\
  &\leq \bar{\gamma}^{1/2} C \left( \| \gamma^{1/2} ( \lambda_{kq} + Wu_{hp} + (K + 1/2)^\prime \phi_{hp} )\|_{L^2(\Gamma_C)} + \| \gamma^{1/2}(\widetilde{W}-W)u_{hp}\|_{L^2(\Gamma_C)} + \| \gamma^{1/2}(\widetilde{K} -K)^\prime \phi_{hp} \|_{L^2(\Gamma_C)} \right) \| v \|_{\widetilde{H}^{1/2}(\Gamma_\Sigma)}. \label{eq:apostEstStabTerm}
\end{align}
The assertion follows with the continuous inf-sup condition \eqref{eq:lame:infsup} proven in \cite[Thm.~3.2.1]{Chernov2006}. 
 \end{proof}
\end{lemma}

Localizing the constraint $V\phi - (K+1/2)u=0$ error by localizing its $H^{1/2}(\Gamma)$-norm, we need to generalize \cite[Thm.~3.2]{carstensen2001posteriori} as $V\phi_{hp} - (K+1/2)u_{hp}$ has no zero mean value on each element due to the stabilization term.
\begin{lemma} \label{lem:localization_inverse_V_apost}
Let $(u,\phi,\lambda)$, $(u_{hp},\phi_{hp},\lambda_{kq})$ be the solution to \eqref{eq:compactWeakForm}, \eqref{eq:DiscreteProblem}, respectively. Then there exists a constant $C>0$ independent of $h$, $p$, $k$, $q$ and $\bar{\gamma}$ such that
\begin{align}
 B(u_{hp},\phi_{hp};0,\phi-\phi_{hp}) \leq C \|\phi-\phi_{hp}  \|_{H^{-1/2}(\Gamma)} \left( \sum_{T\in \mathcal{T}_h} \eta_{V,T}^2 +\sum_{T \in \mathcal{T}_h \cap \Gamma_C} \bar{\gamma}^2 \eta_{C,T}^2   \right)^{1/2}
\end{align} 
with the local error contributions defined in \eqref{eq:local_error_contributions}.
 \begin{proof}
  Let $\phi_e:=\phi-\phi_{hp}$, there trivially holds
\begin{align*}
 B(u_{hp},\phi_{hp};0,\phi_e)= \langle V\phi_{hp} - (K+1/2)u_{hp} , \phi_e \rangle_\Gamma  \leq \|  V\phi_{hp} - (K+1/2)u_{hp} \|_{H^{1/2}(\Gamma)} \|\phi_e  \|_{H^{-1/2}(\Gamma)}.
\end{align*}
The next step is to localize $\|  V\phi_{hp} - (K+1/2)u_{hp} \|_{H^{1/2}(\Gamma)}$ by \cite[Thm.~3.2]{carstensen2001posteriori} as in \cite[Cor.~4.2]{carstensen2001posteriori} but with the Poincar\'e inequality which explicitly includes the mean value, i.e.~
\begin{align}
 C \| f \|_{L^2(T)}^2 \leq \left( \int_T f \; ds \right)^2 + \| \nabla f \|_{L^2(T)}^2 \quad \forall f \in H^1(T).
\end{align}
Let $f:=V\phi_{hp} - (K+1/2)u_{hp}$, by the mapping properties of the boundary integral operators \cite{costabel1988boundary} $f  \in H^1(\Gamma)$. Let $\Phi:=(\Phi_1,\ldots,\Phi_M)$ be a finite partition of unity of $\Gamma = \partial \Omega$ consisting of hat functions $\Phi_j$ with support $\omega_j=\operatorname{supp}(\Phi_j)$ matching exactly a finite number of elements $T \in \mathcal{T}_h$ and set $d_j=\operatorname{diam}(\omega_j)$.
By \cite[Thm.~3.2]{carstensen2001posteriori} and in light of \cite[Cor.~4.2]{carstensen2001posteriori} there holds
\begin{align*}
 \| f\|^2_{H^{1/2}(\Gamma)} &\leq C \sum_{j=1}^M d_j (1+d_j^2)^{1/2} \| \nabla (\Phi_j f) \|^2_{L^2(\omega_j)} \leq C \sum_{j=1}^M d_j (1+d_j^2)^{1/2} \left( \| \nabla f \|^2_{L^2(\omega_j)} + \left( \int_{\omega_j} f \; ds \right)^2 \right) \\
 &\leq C \sum_{T\in \mathcal{T}_h} h_T \| \nabla f \|^2_{L^2(T)} + C \sum_{j=1}^M d_j (1+d_j^2)^{1/2}  \left( \int_{\omega_j} f \; ds \right)^2 .
\end{align*}
From  \eqref{eq:DiscreteMixedVarEq1} with $v_{hp}=0$, $\psi_{hp}=1_{\omega_j} \in Y_{hp}$, \eqref{eq:theo_K_tilde} and $L^2(\Gamma) \subset H^{-1/2}(\Gamma)$ we obtain
\begin{align*}
 \left(\int_{\omega_j} V \phi_{hp} - (K+1/2)u_{hp} \, ds \right)^2  &= \left(\langle \gamma ( \lambda_{kq} + \widetilde{W}u_{hp} + (\widetilde{K} + 1/2)^\prime \phi_{hp}),  (\widetilde{K} + 1/2)^\prime 1_{\omega_j} \rangle_{\Gamma_C} \right)^2 \\
 & \leq \|  \gamma^{1/2} ( \lambda_{kq} + \widetilde{W}u_{hp} + (\widetilde{K} + 1/2)^\prime \phi_{hp})\|_{L^2(\Gamma_C)}^2 \| \gamma^{1/2} (\widetilde{K} + 1/2)^\prime 1_{\omega_j}\|_{L^2(\Gamma_C)}^2 \\
 & \leq C \bar{\gamma} \|  \gamma^{1/2} ( \lambda_{kq} + \widetilde{W}u_{hp} + (\widetilde{K} + 1/2)^\prime \phi_{hp})\|_{L^2(\Gamma_C)}^2 \|  1_{\omega_j}\|_{L^2(\Gamma_C)}^2.
\end{align*}
Hence, 
\begin{align*}
 \sum_{j=1}^M d_j (1+d_j^2)^{1/2}  \left( \int_{\omega_j} f \; dx \right)^2 &\leq C\bar{\gamma} \|  \gamma^{1/2} ( \lambda_{kq} + \widetilde{W}u_{hp} + (\widetilde{K} + 1/2)^\prime \phi_{hp})\|_{L^2(\Gamma_C)}^2  \sum_{j=1}^M d_j \|  1_{\omega_j}\|_{L^2(\Gamma_C)}^2 \\
 & \leq C\bar{\gamma} \|  \gamma^{1/2} ( \lambda_{kq} + \widetilde{W}u_{hp} + (\widetilde{K} + 1/2)^\prime \phi_{hp})\|_{L^2(\Gamma_C)}^2
\end{align*}
as $d_j \leq C$. Thus
\begin{align*}
 C \|  V\phi_{hp} - (K+1/2)u_{hp} \|_{H^{1/2}(\Gamma)}^2 \leq \sum_{T\in \mathcal{T}_h} h_T \| \nabla (V\phi_{hp} - (K+1/2)u_{hp}) \|^2_{L^2(T)} + \bar{\gamma}\|  \gamma^{1/2} ( \lambda_{kq} + \widetilde{W}u_{hp} + (\widetilde{K} + 1/2)^\prime \phi_{hp})\|_{L^2(\Gamma_C)}^2
\end{align*}
which completes the proof.
 \end{proof}
\end{lemma}

With the previous three lemmas we can proof the reliability of the a posteriori error estimate. The complementarity term $ \langle (\lambda^{kq}_n)^+,(g-u^{hp}_n)^+\rangle_{\Gamma_C}$ does not seem to allow efficiency estimates for $p$ or $q \geq 2$.
\begin{theorem}\label{thm:aposterioriErrorEst}
Let $(u,\phi,\lambda)$, $(u_{hp},\phi_{hp},\lambda_{kq})$ be the solution to \eqref{eq:compactWeakForm}, \eqref{eq:DiscreteProblem}, respectively. Then there holds 
  \begin{align*}
  C&\left(\| u-u_{hp}\|_{\widetilde{H}^{1/2}(\Gamma_\Sigma)}^2 +  \| \phi - \phi_{hp}\|_{H^{-1/2}(\Gamma)}^2 + \| \lambda-\lambda_{kq} \|_{\widetilde{H}^{-1/2}(\Gamma_C)}^2  \right)\\
  & \leq    \sum_{T \in \mathcal{T}_h \cap \Gamma_N} \eta_{N,T}^2  + \sum_{T\in \mathcal{T}_h} \eta_{V,T}^2
   +  \sum_{T \in \mathcal{T}_h \cap \Gamma_C}   (1 + \bar{\gamma}^2) \eta_{C,T}^2 + \bar{\gamma}^2 \eta_{W,T}^2 + \bar{\gamma}^2 \eta_{K,T}^2+    \langle (\lambda^{kq}_n)^+,(g-u^{hp}_n)^+\rangle_{\Gamma_C} + \| (g-u^{hp}_n)^-\|_{H^{1/2}(\Gamma_C)}^2\\
   & \qquad    + \sum_{T\in \mathcal{T}_h \cap \Gamma_C} \left(\|(\lambda^{kq}_n)^-\|_{\widetilde{H}^{-1/2}(T)}^2  
+ \| (|\lambda^{kq}_t|-\mathcal{F} )^+ \|_{\widetilde{H}^{-1/2}(T)}^2 \right)  -\langle (|\lambda^{kq}_t|-\mathcal{F} )^-, |u^{hp}_t| \rangle_{\Gamma_C}  -\langle  \lambda^{kq}_t,u^{hp}_t \rangle_{\Gamma_C}
 +\langle |\lambda^{kq}_t| , |u^{hp}_t| \rangle_{\Gamma_C}
  \end{align*} 
with a constant $C>0$ independent of $h$, $p$, $k$, $q$ and $\bar{\gamma}$ and with the local error contributions defined in \eqref{eq:local_error_contributions}. 
 \begin{proof}
  By Lemma~\ref{lem:GalerkinOrtho} with the test functions $u_{hp} - v_{hp}$ and zero, and by \eqref{eq:compactWeakForm1} there holds
  \begin{align*}
   \| u-u_{hp}\|_W^2 + \| \phi - \phi_{hp}\|_V^2 &= B(u-u_{hp},\phi - \phi_{hp};u-u_{hp},\phi - \phi_{hp})\\
   &= \underbrace{\langle f, u-v_{hp} \rangle_{\Gamma_N} - \langle \lambda_{kq}, u-v_{hp} \rangle_{\Gamma_C} - B(u_{hp},\phi_{hp};u-v_{hp},0)}_{\rm I} + \underbrace{\langle \lambda - \lambda_{kq}, u_{hp} - u \rangle_{\Gamma_C}}_{\rm II} \\
   &\quad + \underbrace{\langle \gamma ( \lambda_{kq} + \widetilde{W}u_{hp} + (\widetilde{K} + 1/2)^\prime \phi_{hp}),  \widetilde{W}(u_{hp}-v_{hp})\rangle_{\Gamma_C}}_{\rm III}+B(u_{hp},\phi_{hp};0,\phi - \phi_{hp}).
  \end{align*}
We bound the individual terms as follows. Let $v_{hp} = u_{hp} + I_{hp}(u-u_{hp})$ with $I_{hp}$ the Scott-Zhang or Cl\'{e}ment interpolant.
By linearity, Cauchy-Schwarz inequality, interpolation error estimate and Young's inequality ($\epsilon>0$ arbitrary) we obtain analogously to the proof of Lemma~\ref{lem:aposteriori_lambda} that
\begin{align*}
 {\rm I} &  = \langle f - Wu_{hp}- (K+1/2)'\phi_{hp}, u-v_{hp}\rangle _{\Gamma_N}  + \langle -\lambda_{kq} - Wu_{hp}- (K+1/2)'\phi_{hp}, u-v_{hp}\rangle _{\Gamma_C} \\ 
& \leq \frac{C}{\epsilon} \! \!\sum_{T \in \mathcal{T}_h \cap \Gamma_N} \frac{h_T}{p_T} \| f -Wu_{hp} - (K+1/2)^\prime \phi_{hp} \|^2_{L^2(T)}  + \frac{C}{\epsilon}\!\! \sum_{T \in \mathcal{T}_h \cap \Gamma_C} \frac{h_T}{p_T} \| \lambda_{kq} +Wu_{hp} + (K+1/2)^\prime \phi_{hp} \|_{L^2(T)}^2 + \epsilon \| u-u_{hp} \|_{\widetilde{H}^{1/2}(\Gamma_\Sigma)}^2.
\end{align*}
The stabilization term is estimate as in \eqref{eq:apostEstStabTerm} but additionally with Youngs inequality and with the $H^{1/2}$-stability of the interpolation, i.e.
\begin{align*}
  {\rm III}  &  \leq \bar{\gamma}^{1/2} C (
 \| \gamma^{1/2} ( \lambda_{kq} + Wu_{hp} + (K + 1/2)^\prime \phi_{hp}) \|_{L^2(\Gamma_C)}    +   \| \gamma^{1/2}(\widetilde{W}-W)u_{hp}\|_{L^2(\Gamma_C)} + \| \gamma^{1/2}(\widetilde{K} -K)^\prime \phi_{hp} \|_{L^2(\Gamma_C)} ) \\
 & \quad \cdot \| u_{hp} - v_{hp} \|_{\widetilde{H}^{1/2}(\Gamma_\Sigma)}  \\
 & \leq \bar{\gamma} \frac{C}{\epsilon} (
 \| \gamma^{1/2} ( \lambda_{kq} + Wu_{hp} + (K + 1/2)^\prime \phi_{hp}) \|_{L^2(\Gamma_C)} +   \| \gamma^{1/2}(\widetilde{W}-W)u_{hp}\|_{L^2(\Gamma_C)} + \| \gamma^{1/2}(\widetilde{K} -K)^\prime \phi_{hp} \|_{L^2(\Gamma_C)} )^2\\
 & \quad + 2\epsilon \| u - u_{hp} \|_{\widetilde{H}^{1/2}(\Gamma_\Sigma)}^2 .
\end{align*}
From Lemma~\ref{lem:aposterioriContactPart}, Young's inequality and von Petersdorff inequality, c.f.~e.g.~\cite[Lem.~1]{Heuer2001Additive}, we obtain
\begin{align*}
 {\rm II} \leq  & 
\langle (\lambda^{kq}_n)^+,(g-u^{hp}_n)^+\rangle_{\Gamma_C}+ \frac{1}{4\epsilon} \| (g-u^{hp}_n)^-\|_{H^{1/2}(\Gamma_C)}^2 + \frac{1}{4\epsilon} \|(\lambda^{kq}_n)^-\|_{\widetilde{H}^{-1/2}(\Gamma_C)}^2  
+ \frac{1}{4\epsilon}\| (|\lambda^{kq}_t|-\mathcal{F} )^+ \|_{\widetilde{H}^{-1/2}(\Gamma_C)}^2  \\
& -\langle (|\lambda^{kq}_t|-\mathcal{F} )^-, |u^{hp}_t| \rangle_{\Gamma_C}  -\langle  \lambda^{kq}_t,u^{hp}_t \rangle_{\Gamma_C}
 +\langle |\lambda^{kq}_t| , |u^{hp}_t| \rangle_{\Gamma_C} 
  +\epsilon \|\lambda_{kq} -\lambda \|_{\widetilde{H}^{-1/2}(\Gamma_C)}^2 + 2\epsilon\|u_{hp}-u \|_{\widetilde{H}^{1/2}(\Gamma_\Sigma)}^2\\
\leq & \langle (\lambda^{kq}_n)^+,(g-u^{hp}_n)^+\rangle_{\Gamma_C}+ \frac{1}{4\epsilon} \| (g-u^{hp}_n)^-\|_{H^{1/2}(\Gamma_C)}^2 + \frac{C}{4\epsilon} \sum_{T\in \mathcal{T}_h \cap \Gamma_C} \left(\|(\lambda^{kq}_n)^-\|_{\widetilde{H}^{-1/2}(T)}^2  
+ \| (|\lambda^{kq}_t|-\mathcal{F} )^+ \|_{\widetilde{H}^{-1/2}(T)}^2 \right)  \\
& -\langle (|\lambda^{kq}_t|-\mathcal{F} )^-, |u^{hp}_t| \rangle_{\Gamma_C}  -\langle  \lambda^{kq}_t,u^{hp}_t \rangle_{\Gamma_C}
 +\langle |\lambda^{kq}_t| , |u^{hp}_t| \rangle_{\Gamma_C} 
  +\epsilon \|\lambda_{kq} -\lambda \|_{\widetilde{H}^{-1/2}(\Gamma_C)}^2 + 2\epsilon\|u_{hp}-u \|_{\widetilde{H}^{1/2}(\Gamma_\Sigma)}^2.
\end{align*}
$B(u_{hp},\phi_{hp};0,\phi - \phi_{hp})$ is estimated by Lemma~\ref{lem:localization_inverse_V_apost}.
Thus, with the coercivity of $W$ and $V$ and with Lemma~\ref{lem:aposteriori_lambda} to bound $\epsilon \|\lambda_{kq} -\lambda \|_{\widetilde{H}^{-1/2}(\Gamma_C)}^2$ we obtain
  \begin{align*}
  &  (\alpha_W - (C+5)\epsilon) \|u_{hp}-u \|_{\widetilde{H}^{1/2}(\Gamma_\Sigma)}^2 + (\alpha_V - (C+1) \epsilon) \| \phi - \phi_{hp}\|_{H^{-1/2}(\Gamma)}^2  \leq  \frac{C}{\epsilon}(1+\epsilon^2) \sum_{T \in \mathcal{T}_h \cap \Gamma_N} \eta_{N,T}^2 + \frac{C}{\epsilon}  \sum_{T\in \mathcal{T}_h} \eta_{V,T}^2 \\
  & \qquad
   + \frac{C}{\epsilon}(1+\epsilon^2) \sum_{T \in \mathcal{T}_h \cap \Gamma_C}   (1 + \bar{\gamma}^2) \eta_{C,T}^2 + \bar{\gamma}^2 \eta_{W,T}^2 + \bar{\gamma}^2 \eta_{K,T}^2 + 
   \langle (\lambda^{kq}_n)^+,(g-u^{hp}_n)^+\rangle_{\Gamma_C} + \frac{1}{4\epsilon} \| (g-u^{hp}_n)^-\|_{H^{1/2}(\Gamma_C)}^2 \\
   & \qquad + \frac{C}{4\epsilon} \sum_{T\in \mathcal{T}_h \cap \Gamma_C} \left(\|(\lambda^{kq}_n)^-\|_{\widetilde{H}^{-1/2}(T)}^2  
+ \| (|\lambda^{kq}_t|-\mathcal{F} )^+ \|_{\widetilde{H}^{-1/2}(T)}^2 \right)  -\langle (|\lambda^{kq}_t|-\mathcal{F} )^-, |u^{hp}_t| \rangle_{\Gamma_C}  -\langle  \lambda^{kq}_t,u^{hp}_t \rangle_{\Gamma_C}
 +\langle |\lambda^{kq}_t| , |u^{hp}_t| \rangle_{\Gamma_C}.
  \end{align*}
Choosing $\epsilon>0$ sufficiently small yields
  \begin{align*}
  & C \left( \|u_{hp}-u \|_{\widetilde{H}^{1/2}(\Gamma_\Sigma)}^2 + \| \phi - \phi_{hp}\|_{H^{-1/2}(\Gamma)}^2 \right) \leq   \sum_{T \in \mathcal{T}_h \cap \Gamma_N} \eta_{N,T}^2 + \sum_{T\in \mathcal{T}_h} \eta_{V,T}^2 
   +  \sum_{T \in \mathcal{T}_h \cap \Gamma_C}   (1 + \bar{\gamma}^2) \eta_{C,T}^2 + \bar{\gamma}^2 \eta_{W,T}^2 + \bar{\gamma}^2 \eta_{K,T}^2 \\
   & \qquad  +    \langle (\lambda^{kq}_n)^+,(g-u^{hp}_n)^+\rangle_{\Gamma_C} + \| (g-u^{hp}_n)^-\|_{H^{1/2}(\Gamma_C)}^2 + \sum_{T\in \mathcal{T}_h \cap \Gamma_C} \left(\|(\lambda^{kq}_n)^-\|_{\widetilde{H}^{-1/2}(T)}^2  
+ \| (|\lambda^{kq}_t|-\mathcal{F} )^+ \|_{\widetilde{H}^{-1/2}(T)}^2 \right)\\
   & \qquad   -\langle (|\lambda^{kq}_t|-\mathcal{F} )^-, |u^{hp}_t| \rangle_{\Gamma_C}  -\langle  \lambda^{kq}_t,u^{hp}_t \rangle_{\Gamma_C}
 +\langle |\lambda^{kq}_t| , |u^{hp}_t| \rangle_{\Gamma_C}.
  \end{align*}
Inserting this result into the estimate from Lemma~\ref{lem:aposteriori_lambda} yields the assertion.
 \end{proof}
\end{theorem}

 \section{Modifications for Coulomb friction}
 \label{sec:CoulombFriction}

 Tresca friction may yield unphysical behavior, namely non-zero tangential traction and stick-slip transition outside the actual contact zone. Therefore, in many applications the more realistic Coulomb friction is applied, in which the friction threshold $\mathcal{F}$ is replaced by $\mathcal{F}|\sigma_n(u)|$. In the discretization which we present here only the Lagrange multiplier set needs to be adapted, namely $M^+(\mathcal{F})$ and $M_{kq}^+(\mathcal{F})$ become
 \begin{align} 
 M^+(\mathcal{F}\lambda_n)&:=\left\lbrace \mu \in \widetilde{H}^{-1/2}(\Gamma_C): \langle \mu,v\rangle_{\Gamma_C} \leq \langle \mathcal{F}\lambda_n,|v_t| \rangle_{\Gamma_C} \forall v \in \widetilde{H}^{1/2}(\Gamma_\Sigma), v_n \leq 0 \right\rbrace\ , \label{eq:CoulombModLagrangeSet} \\
M_{kq}^+(\mathcal{F}\lambda^{kq}_n)&:= \left\{ \mu_{kq} \in L^2(\Gamma_C): \mu_{kq}|_E \circ \Psi_T \in \left[\mathbb{P}_{q_T}(\hat{T})\right]^d \; \forall T \in \mathcal{T}_k^C,\ \mu^{kq}_n(x) \geq 0, \ |\mu^{kq}_t(x)| \leq \mathcal{F}\lambda^{kq}_n(x)\; \text{for }x\in G_{kq} \right\}.
 \end{align}
 As in the a posteriori error estimate for the Tresca frictional case we can simply recall a partial result from \cite[Thm.~15]{Banz2014BEM}, namely:
 \begin{lemma}  \label{lem_coulomb_apost_part1}
  Let $(u,\phi,\lambda)$, $(u_{hp},\phi_{hp},\lambda_{kq})$ be a solution to \eqref{eq:compactWeakForm}, \eqref{eq:DiscreteProblem} with the Lagrange multiplier sets for the Coulomb frictional case, respectively, then there holds
 \begin{align*}
\langle \lambda-\lambda_{kq}, u_{hp}-u\rangle_{\Gamma_C} 
 \leq & 
 \langle (\lambda^{kq}_n)^+,(g-u^{hp}_n)^+\rangle_{\Gamma_C}+ \|\lambda_{kq} -\lambda \|_{\widetilde{H}^{-1/2}(\Gamma_C)} \| (g-u^{hp}_n)^-\|_{H^{1/2}(\Gamma_C)} + \|(\lambda^{kq}_n)^-\|_{\widetilde{H}^{-1/2}(\Gamma_C)} \|u_{hp}-u \|_{\widetilde{H}^{1/2}(\Gamma_\Sigma)} \\
 & +\| (|\lambda^{kq}_t|-\mathcal{F} \lambda^{kq}_n )^+ \|_{\widetilde{H}^{-1/2}(\Gamma_C)} \| u- u_{hp} \|_{\widetilde{H}^{1/2}(\Gamma_\Sigma)}  -\langle (|\lambda^{kq}_t|-\mathcal{F} \lambda^{kq}_n )^-, |u^{hp}_t| \rangle_{\Gamma_C} \\
& -\langle  \lambda^{kq}_t,u^{hp}_t \rangle_{\Gamma_C} +\langle |\lambda^{kq}_t| , |u^{hp}_t| \rangle_{\Gamma_C} 
+ \|\mathcal{F}\|_{L^\infty(\Gamma_C)}\|\lambda_{kq} - \lambda \|_{\widetilde{H}^{-1/2}(\Gamma_C)} \| u - u_{hp} \|_{\widetilde{H}^{1/2}(\Gamma_\Sigma)} 
\end{align*}
where $v^+=\max\{0,v\}$ and $v^-=\min\{0,v\}$, i.e.~$v=v^+ +v^-$.
 \end{lemma}
 \begin{flushleft}
Therewith the a posteriori error estimate becomes: \end{flushleft}
 \begin{theorem}
Let $(u,\phi,\lambda)$, $(u_{hp},\phi_{hp},\lambda_{kq})$ be a solution to \eqref{eq:compactWeakForm}, \eqref{eq:DiscreteProblem} with the Lagrange multiplier sets for the Coulomb frictional case, respectively. If $\| \mathcal{F}\|_{L^\infty(\Gamma_C)}$ is sufficiently small, then there holds
 \begin{align*}
  &C\left(\| u-u_{hp}\|_{\widetilde{H}^{1/2}(\Gamma_\Sigma)}^2 +  \| \phi - \phi_{hp}\|_{H^{-1/2}(\Gamma)}^2 + \| \lambda-\lambda_{kq} \|_{\widetilde{H}^{-1/2}(\Gamma_C)}^2  \right) \leq   \sum_{T \in \mathcal{T}_h \cap \Gamma_N} \eta_{N,T}^2    + \sum_{T \in \mathcal{T}_h \cap \Gamma_C}   (1 + \bar{\gamma}^2) \eta_{C,T}^2 + \bar{\gamma}^2 \eta_{W,T}^2 + \bar{\gamma}^2 \eta_{K,T}^2 \\ 
   & \qquad +  \sum_{T\in \mathcal{T}_h} \eta_{V,T}^2  + \langle (\lambda^{kq}_n)^+,(g-u^{hp}_n)^+\rangle_{\Gamma_C}  + \| (g-u^{hp}_n)^-\|_{H^{1/2}(\Gamma_C)}^2  +  \sum_{T\in \mathcal{T}_h \cap \Gamma_C} \left(\|(\lambda^{kq}_n)^-\|_{\widetilde{H}^{-1/2}(T)}^2  
+ \| (|\lambda^{kq}_t|-\mathcal{F}\lambda_n^{kq} )^+ \|_{\widetilde{H}^{-1/2}(T)}^2 \right) \\
   & \qquad -\langle (|\lambda^{kq}_t|-\mathcal{F}\lambda_n^{kq} )^-, |u^{hp}_t| \rangle_{\Gamma_C}  -\langle  \lambda^{kq}_t,u^{hp}_t \rangle_{\Gamma_C}
 +\langle |\lambda^{kq}_t| , |u^{hp}_t| \rangle_{\Gamma_C}
  \end{align*} 
with a constant $C>0$ independent of $h$, $p$, $k$, $q$ and $\bar{\gamma}$ and with the local error contributions defined in \eqref{eq:local_error_contributions}.
 
\begin{proof}
Since Lemma~\ref{lem:aposteriori_lambda} is independent of $\mathcal{F}$ it holds for the Coulomb frictional case as well. 
From Lemma~\ref{lem:aposteriori_lambda}, subadditivity of the square root $\sqrt{a^2+b^2}\leq a +b$ and Young's inequality we obtain  with 
\[
b^2:= \sum_{T \in \mathcal{T}_h \cap \Gamma_N}  \eta_{N,T}^2 + \sum_{T \in \mathcal{T}_h \cap \Gamma_C} \left ( (1+ \bar{\gamma}^2) \eta_{C,T}^2 + \bar{\gamma}^2 \eta_{W,T}^2 + \bar{\gamma}^2 \eta_{K,T}^2 \right ) 
\]
that
\begin{align*}
  & \|\mathcal{F}\|_{L^\infty(\Gamma_C)}\|\lambda_{kq} - \lambda \|_{\widetilde{H}^{-1/2}(\Gamma_C)} \| u - u_{hp} \|_{\widetilde{H}^{1/2}(\Gamma_\Sigma)}  
	\leq  C \|\mathcal{F}\|_{L^\infty(\Gamma_C)} \| u - u_{hp} \|_{\widetilde{H}^{1/2}(\Gamma_\Sigma)} \left [  b^2 + \|u-u_h\|^2_{\widetilde{H}^{1/2}(\Gamma_\Sigma)} + \|\phi-\phi_{hp}\|^2_{{H}^{-1/2}(\Gamma_\Sigma)}\right]^{1/2} \\
	& \qquad \quad \leq  C \|\mathcal{F}\|_{L^\infty(\Gamma_C)} \| u - u_{hp} \|_{\widetilde{H}^{1/2}(\Gamma_\Sigma)} \left[ \|u-u_h\|_{\widetilde{H}^{1/2}(\Gamma_\Sigma)} + \|\phi-\phi_{hp}\|_{{H}^{-1/2}(\Gamma_\Sigma)} +b \right]  \\
	& \qquad  \quad\leq  C \|\mathcal{F}\|_{L^\infty(\Gamma_C)} \| u - u_{hp} \|^2_{\widetilde{H}^{1/2}(\Gamma_\Sigma)} + \frac{C}{2} \|\mathcal{F}\|_{L^\infty(\Gamma_C)} \| u - u_{hp} \|^2_{\widetilde{H}^{1/2}(\Gamma_\Sigma)} + \frac{C}{2} \|\mathcal{F}\|_{L^\infty(\Gamma_C)} \| \phi - \phi_{hp} \|^2_{{H}^{-1/2}(\Gamma_\Sigma)} +\epsilon \| u - u_{hp} \|^2_{\widetilde{H}^{1/2}(\Gamma_\Sigma)} \\
	& \qquad  \qquad + \frac{C^2}{4\epsilon}\|\mathcal{F}\|^2_{L^\infty(\Gamma_C)} \left [\sum_{T \in \mathcal{T}_h \cap \Gamma_N}  \eta_{N,T}^2 + \sum_{T \in \mathcal{T}_h \cap \Gamma_C} \left ( (1+ \bar{\gamma}^2) \eta_{C,T}^2 + \bar{\gamma}^2 \eta_{W,T}^2 + \bar{\gamma}^2 \eta_{K,T}^2 \right ) \right]\\
	& \qquad  \quad = \left(\frac{3C}{2} \|\mathcal{F}\|_{L^\infty(\Gamma_C)} + \epsilon \right) \| u - u_{hp} \|_{\widetilde{H}^{1/2}(\Gamma_\Sigma)}^2+ \frac{C}{2}\|\mathcal{F}\|_{L^\infty(\Gamma_C)} \|\phi-\phi_{hp}\|_{{H}^{-1/2}(\Gamma_\Sigma)}^2\\
   & \qquad  \qquad+ \frac{C}{4\epsilon}\|\mathcal{F}\|_{L^\infty(\Gamma_C)}^2 \left(  \sum_{T \in \mathcal{T}_h \cap \Gamma_N}  \eta_{N,T}^2 + \sum_{T \in \mathcal{T}_h \cap \Gamma_C} (1+ \bar{\gamma}^2) \eta_{C,T}^2 + \bar{\gamma}^2 \eta_{W,T}^2 + \bar{\gamma}^2 \eta_{K,T}^2 \right).
\end{align*}
 Arguing exactly the same as in the proof of Theorem~\ref{thm:aposterioriErrorEst}, but using now Lemma~\ref{lem_coulomb_apost_part1} instead of Lemma \ref{lem:aposterioriContactPart} we obtain
\begin{align*}
  &  (\alpha_W - (C+5+1)\epsilon - \frac{3C}{2}\|\mathcal{F}\|_{L^\infty(\Gamma_C)}) \|u_{hp}-u \|_{\widetilde{H}^{1/2}(\Gamma_\Sigma)}^2 + (\alpha_V - (C+1) \epsilon- \frac{C}{2}\|\mathcal{F}\|_{L^\infty(\Gamma_C)}) \| \phi - \phi_{hp}\|_{H^{-1/2}(\Gamma)}^2  \\
	&\qquad \leq  \left[\frac{C}{\epsilon}(1+\epsilon^2) + \frac{C^2}{4\epsilon}\|\mathcal{F}\|^2_{L^\infty(\Gamma_C)} \right] \sum_{T \in \mathcal{T}_h \cap \Gamma_N} \eta_{N,T}^2 + \frac{C}{\epsilon}  \sum_{T\in \mathcal{T}_h} \eta_{V,T}^2 + \frac{1}{4\epsilon} \| (g-u^{hp}_n)^-\|_{H^{1/2}(\Gamma_C)}^2 \\
  & \qquad
   + \left[\frac{C}{\epsilon}(1+\epsilon^2) +\frac{C^2}{4\epsilon}\|\mathcal{F}\|^2_{L^\infty(\Gamma_C)} \right ] \sum_{T \in \mathcal{T}_h \cap \Gamma_C}  \left ( (1 + \bar{\gamma}^2) \eta_{C,T}^2 + \bar{\gamma}^2 \eta_{W,T}^2 + \bar{\gamma}^2 \eta_{K,T}^2 \right ) + 
   \langle (\lambda^{kq}_n)^+,(g-u^{hp}_n)^+\rangle_{\Gamma_C}  \\
   & \qquad \quad + \frac{C}{4\epsilon} \sum_{T\in \mathcal{T}_h \cap \Gamma_C} \left(\|(\lambda^{kq}_n)^-\|^2_{\widetilde{H}^{-1/2}(T)}  + \| (|\lambda^{kq}_t|-\mathcal{F}\lambda_n^{kq} )^+ \|^2_{\widetilde{H}^{-1/2}(T)} \right)  -\langle (|\lambda^{kq}_t|-\mathcal{F} \lambda_n^{kq})^-, |u^{hp}_t| \rangle_{\Gamma_C}  -\langle  \lambda^{kq}_t,u^{hp}_t \rangle_{\Gamma_C} \\
   & \qquad \quad  +\langle |\lambda^{kq}_t| , |u^{hp}_t| \rangle_{\Gamma_C}. 
  \end{align*} 
Hence, if $\epsilon$ and $\|\mathcal{F}\|_{L^\infty(\Gamma_C)}$ are sufficiently small, then
\begin{align*}
    & C \left( \|u_{hp}-u \|_{\widetilde{H}^{1/2}(\Gamma_\Sigma)}^2 +  \| \phi - \phi_{hp}\|_{H^{-1/2}(\Gamma)}^2  \right)
	 \leq \sum_{T \in \mathcal{T}_h \cap \Gamma_N} \eta_{N,T}^2 +  \sum_{T\in \mathcal{T}_h} \eta_{V,T}^2 + 
    \sum_{T \in \mathcal{T}_h \cap \Gamma_C}   (1 + \bar{\gamma}^2) \eta_{C,T}^2 + \bar{\gamma}^2 \eta_{W,T}^2 + \bar{\gamma}^2 \eta_{K,T}^2 \\
		&\qquad + \langle (\lambda^{kq}_n)^+,(g-u^{hp}_n)^+\rangle_{\Gamma_C} +  \| (g-u^{hp}_n)^-\|_{H^{1/2}(\Gamma_C)}^2 +  \sum_{T\in \mathcal{T}_h \cap \Gamma_C} \left(\|(\lambda^{kq}_n)^-\|^2_{\widetilde{H}^{-1/2}(T)}  + \| (|\lambda^{kq}_t|-\mathcal{F}\lambda_n^{kq} )^+ \|^2_{\widetilde{H}^{-1/2}(T)} \right)  \\
&\qquad -\langle (|\lambda^{kq}_t|-\mathcal{F} \lambda_n^{kq})^-, |u^{hp}_t| \rangle_{\Gamma_C}  
  +\langle |\lambda^{kq}_t| , |u^{hp}_t| \rangle_{\Gamma_C} -\langle  \lambda^{kq}_t,u^{hp}_t \rangle_{\Gamma_C}. 
  \end{align*}
Inserting this result into the estimate from Lemma~\ref{lem:aposteriori_lambda} yields the assertion.  
\end{proof}
\end{theorem}

\section{A possible approximation of the stabilization term}
\label{sec:approximation}
Choosing $\widetilde{W}=W$ and $(\widetilde{K}+1/2)^\prime=(K+1/2)^\prime$ works and has been done in \cite{Banz2015Stab}, but it requires to implement non-standard matrices and a lot of computational time. Computationalwise simple to handle is the composition of a projection operator and the boundary integral operator, e.g.~$\widetilde{W}v= \Pi_{HP} Wv$. Here $\Pi_{HP}$ is the $L^2(\Gamma_C)$-projection onto 
\begin{align} \label{eq:StabApproxSpace}
 Z_{HP}&:=\left\{\phi_{HP} \in  L^{2}(\Gamma_C) : \phi_{HP}|_T \circ \Psi_T \in \left[\mathbb{P}_{P_T}(\hat{T})\right]^d \; \forall T \in \mathcal{T}_{H} \right\},
\end{align}
i.e.
\begin{align}
  \langle  \Pi_{HP} v - v , \kappa_{HP} \rangle_{\Gamma_C} =0 \quad \forall\: \kappa_{HP} \in Z_{HP}.
\end{align}
Analogously, we set $\widetilde{K}^\prime v = \Pi_{HP} (K^\prime v)$. With these approximations, the algebraic representation of the core stabilization matrices become:
\begin{gather*}
 M_{ij}=\int_{\Gamma_C} \xi_i(x) \xi_j(x) \; ds_x, \qquad M^{(\gamma)}_{ij}=\int_{\Gamma_C} \gamma(x) \xi_i(x) \xi_j(x) \;
 ds_x, \qquad \widehat{W_{ij}} = \langle W \phi_j,\xi_i \rangle_{\Gamma_C}, \qquad \widehat{K_{ij}^\prime} = \langle K^\prime \psi_j,\xi_i \rangle_{\Gamma_C}, \\
 \langle \gamma  \widetilde{W}u_{hp},  \widetilde{W}v_{hp} \rangle_{\Gamma_C} = \left.\vec{v}\right.^\top \widehat{W}^\top M^{-1} M^{(\gamma)} M^{-1} \widehat{W} \vec{u}, \qquad 
 \langle \gamma \widetilde{W}u_{hp},  \widetilde{K}^\prime \psi_{hp} \rangle_{\Gamma_C} =\left.\vec{\psi}\right.^\top (\widehat{K^\prime})^\top M^{-1} M^{(\gamma)} M^{-1} \widehat{W} \vec{u}, \\ 
 \langle \gamma \widetilde{K}^\prime \phi_{hp}, \widetilde{K}^\prime \psi_{hp} \rangle_{\Gamma_C} = \left.\vec{\psi}\right.^\top (\widehat{K^\prime})^\top M^{-1} M^{(\gamma)} M^{-1} \widehat{K^\prime} \vec{\phi}
\end{gather*}
with $\xi_i$, $\phi_j$ and $\psi_j$ basis functions of $Z_{HP}$, $X_{hp}$ and $Y_{hp}$, respectively.
The well-posedness of \eqref{eq:DiscreteProblem} and a priori error convergence rate follow with the following lemma.

\begin{lemma} \label{lem:ApproxStabTerme}
Let $\mathcal{T}_H$ be obtained by mesh refinement of $\mathcal{T}_h|_{\Gamma_C}$ and $P_T \geq p_{T^\prime}$ for $T \in \mathcal{T}_H$ a child element of $T^\prime \in \mathcal{T}_h|_{\Gamma_C}$. For $\widetilde{W}= \Pi_{HP} W$ and $\widetilde{K}^\prime = \Pi_{HP} (K^\prime)$ the inverse inequalities and mapping properties \eqref{eq:theo_W_tilde} and \eqref{eq:theo_K_tilde} hold as well as the a priori error estimates
 \begin{align*}
 \| \gamma^{1/2} (\widetilde{W}v -Wv) \|_{L^2(\Gamma_C)} \leq C \bar{\gamma}^{1/2}\left(\frac{h}{p^2}\right)^{1/2} \left(\frac{H}{P}\right)^\beta \|Wv \|_{H^\beta(\Gamma)}  \text{ and } \| \gamma^{1/2} (\widetilde{K}^\prime\phi -K^\prime \phi) \|_{L^2(\Gamma_C)}   \leq C \bar{\gamma}^{1/2}\left(\frac{h}{p^2}\right)^{1/2} \left(\frac{H}{P}\right)^\beta \|K^\prime \phi \|_{H^\beta(\Gamma)}
\end{align*}
for all $v \in H^1(\Gamma)$, $\phi \in L^2(\Gamma)$, $\beta \geq 0$ and with a constant $C>0$ independent of $h$, $p$, $H$, $P$ and $\bar{\gamma}$.
 \begin{proof}
  By construction we have $\gamma\phi_{HP} \in Z_{HP}$ for all $\phi_{HP} \in Z_{HP}$ since $\gamma$ is piecewise constant w.r.t.~$\mathcal{T}_h|_{\Gamma_C}$ and is also piecewise constant w.r.t.~$\mathcal{T}_H$. Consequently, we have
  \begin{align}
   \| \gamma^{1/2} (\Pi_{HP}v -v) \|_{L^2(\Gamma_C)}^2 &= \langle \Pi_{HP}v -v, \gamma (\Pi_{HP}v -v) \rangle_{\Gamma_C} =\langle \Pi_{HP}v -v, \gamma (v_{HP} -v) \rangle_{\Gamma_C} \nonumber \\
   &\leq \| \gamma^{1/2} (\Pi_{HP}v -v) \|_{L^2(\Gamma_C)} \|\gamma^{1/2} (v_{HP} -v) \|_{L^2(\Gamma_C)} \label{eq:ProjectionStability}
  \end{align}
with $v_ {HP}\in Z_{HP}$ arbitrary and, hence,
  \begin{align} \label{eq:ProjectionErrorEst}
   \| \gamma^{1/2} (\Pi_{HP}v -v) \|_{L^2(\Gamma_C)} \leq C \bar{\gamma}^{1/2}\left(\frac{h}{p^2}\right)^{1/2} \left(\frac{H}{P}\right)^\beta \|v \|_{H^\beta(\Gamma_C)}, \quad \beta \geq 0
  \end{align}
  which gives the a priori error estimate assertion.
By the triangle inequality and \eqref{eq:ProjectionStability} with $v_{HP}=0$ we have
\[
\|\gamma^{1/2}\Pi_{HP}v\|_{L^2(\Gamma_C)}\leq \|\gamma^{1/2}(\Pi_{HP}v-v)\|_{L^2(\Gamma_C)}+ \|\gamma^{1/2} v\|_{L^2(\Gamma_C)} \leq 2  \|\gamma^{1/2} v\|_{L^2(\Gamma_C)}.
\]
Hence, analogously to \eqref{eq:ProjectionStability}, we have the $\gamma$-scaled $L^2$-stability and thus, we obtain the inverse inequalities \eqref{eq:theo_W_tilde} and \eqref{eq:theo_K_tilde}, i.e.
\begin{align*}
 \| \gamma^{1/2} \Pi_{HP} W v_{hp} \|_{L^2(\Gamma_C)} & \leq 2\| \gamma^{1/2} W v_{hp} \|_{L^2(\Gamma_C)} \leq \bar{\gamma}^{1/2} C \| v_{hp}\|_{\widetilde{H}^{1/2}(\Gamma_\Sigma)} \quad \forall v_{hp} \in X_{hp}\\
 \| \gamma^{1/2} \Pi_{HP} K^\prime \psi_{hp} \|_{L^2(\Gamma_C)} & \leq 2\| \gamma^{1/2}  K^\prime \psi_{hp} \|_{L^2(\Gamma_C)}  \leq \bar{\gamma}^{1/2} C \| \psi_{hp}\|_{H^{-1/2}(\Gamma)} \quad \forall v_{hp} \in Y_{hp},
\end{align*}
where we use for the final step the inverse inequalities in \cite[Thm.~5]{Banz2015Stab}. The mapping properties stated in \eqref{eq:theo_W_tilde} and \eqref{eq:theo_K_tilde} are satisfied trivially.
 \end{proof} 
\end{lemma}

For the proof of Lemma~\ref{lem:ApproxStabTerme} it is essential, that $\gamma\phi_{HP} \in Z_{HP}$ for all $\phi_{HP} \in Z_{HP}$. In particular this implies that $Z_{HP}$ consists of discontinuous piecewise polynomial functions, as suggested by the mapping properties of $W$ and $K^\prime$ themselves. Thus we cannot perform global integration by parts to compute the entires of the matrix $\widehat{W}$ but only elementwise integration by parts. In total the computational costs for the stabilization matrices are similar to the matrix representation of $B(\cdot,\cdot;\cdot,\cdot)$ and significantly better than in the case of no approximation in the stabilization term, i.e.~$\widetilde{W}=W$ and $(\widetilde{K}+1/2)^\prime=(K+1/2)^\prime$.\\
Coming back to the elementwise integration by parts, for the hypersingular integral operator in 3d it holds in a distributional sense that, c.f.~\cite[Eq.~2.12]{han1994boundary},
\begin{align} \label{eq:han_W_rep}
 W u(x) &= - \sum_{k=1}^3 \frac{\partial}{\partial S_k(x)} \int_\Gamma \frac{\mu}{4\pi} \frac{1}{|x-y|} \frac{\partial u(y)}{\partial S_k(y)}\, ds_y - D(\partial_x,n_x) \int_\Gamma \left(4 \mu^2 G(x,y) - \frac{\mu}{2\pi} I \frac{1}{|x-y|} \right) D(\partial_y,n_y) u(y) \,ds_y \nonumber \\
 & \quad +\frac{\mu}{4\pi} \int_\Gamma \left( D(\partial_x,n_x) \frac{1}{|x-y|} D(\partial_y,n_y) \right)^\top u(y) \, ds_y
\end{align}
with the differential operator
\begin{align*}
 D_{ij}(\partial_x,n_x)=n_j(x) \frac{\partial}{\partial x_i} - n_i(x) \frac{\partial}{\partial x_j}, \qquad 1 \leq i,j \leq 3.
\end{align*}
In particular,
\begin{gather*}
D_{11}(\partial_x,n_x)=D_{22}(\partial_x,n_x)=D_{33}(\partial_x,n_x)=0\\
D_{32}(\partial_x,n_x)=-D_{23}(\partial_x,n_x) = n_2(x) \frac{\partial}{\partial x_3} - n_3(x) \frac{\partial}{\partial x_2} =: \frac{\partial}{\partial S_1(x)}\\
D_{13}(\partial_x,n_x)=-D_{31}(\partial_x,n_x) = n_3(x) \frac{\partial}{\partial x_1} - n_1(x) \frac{\partial}{\partial x_3} =: \frac{\partial}{\partial S_2(x)}\\
D_{21}(\partial_x,n_x)=-D_{12}(\partial_x,n_x) = n_1(x) \frac{\partial}{\partial x_2} - n_2(x) \frac{\partial}{\partial x_1}=: \frac{\partial}{\partial S_3(x)}.
\end{gather*}

\begin{lemma} \label{lem:part_int_s_k_term}
 Let $T$ be a smooth, open surface with piecewise smooth boundary $\partial T$. Let $\vec{t}\,(x)$ be the unit tangential vector along the boundary of $T$ with mathematically positive orientation and let $u$, $v$ be sufficiently smooth scalar functions. Then there holds
 \begin{align}
  \int_T \frac{\partial u(x)}{\partial S_k(x)} v(x) \, ds_x = \int_{\partial T} u(x) v(x) \vec{t}_k(x) \, d\sigma_x - \int_T u(x) \frac{\partial v(x)}{\partial S_k(x)} \, ds_x.
  \end{align}
 \begin{proof}
It holds
  \begin{align*}
   \int_T \frac{\partial u(x)}{\partial S_k(x)} v(x) \, ds_x &= \int_T \left( \vec{n}_x \times \nabla u(x) \right)_k v(x) \,ds_x =  \int_T \left( \vec{n}_x \times \nabla u(x) \right) \left(\vec{e}_k v(x) \right) \,ds_x \\
   &= \int_{\partial T} u(x) v(x) \vec{e}_k \vec{t}(x) \, d\sigma_x - \int_T u(x) \vec{n}_x \left( \nabla \times \vec{e}_k v(x) \right) \,ds_x \\
   & = \int_{\partial T} u(x) v(x) \vec{t}_k(x) \, d\sigma_x - \int_T u(x) \frac{\partial v(x)}{\partial S_k(x)} \, ds_x
  \end{align*}
where the integration by parts follows from \cite[Proof of Lem.~6.16]{steinbach2008numerical}.
 \end{proof}
\end{lemma}

\begin{theorem} \label{thm:elementPartIntW}
 There holds the elementwise integrated by parts representation of the hypersingular boundary integral operator in 3d:
 \begin{align}
  \langle Wu,v \rangle_{\Gamma_C}  = \sum_{T \in \mathcal{T}_H} & \frac{\mu}{4 \pi} \int_\Gamma \int_T  \frac{1}{|x-y|} \sum_{k=1}^3 \frac{\partial \vec{u}(y)}{\partial S_k(y)} \cdot \frac{\partial \vec{v}(x)}{\partial S_k(x)}\,ds_x\, ds_y -  \frac{\mu}{4 \pi} \int_\Gamma \int_{\partial T} \frac{1}{|x-y|} \sum_{k=1}^3 \vec{v}(x) \frac{\partial \vec{u}(y)}{\partial S_k(y)} \vec{t}_k(x) \, d\sigma_x \, ds_y \nonumber \\
 & - \int_\Gamma \int_T \left( D(\partial_x,n_x)\vec{v}(x) \right)^\top \left(4 \mu^2 G(x,y) - \frac{\mu}{2\pi} I \frac{1}{|x-y|} \right) D(\partial_y,n_y)\vec{u}(y) \, ds_x\, ds_y\nonumber \\
 &  - \int_\Gamma \int_{\partial T} \vec{v}^\top(x) T(x)
 \left(4 \mu^2 G(x,y) - \frac{\mu}{2\pi} I \frac{1}{|x-y|} \right) D(\partial_y,n_y)\vec{u}(y) \, d\sigma_x\,ds_y \nonumber \\
 &+ \frac{\mu}{4 \pi} \int_\Gamma \int_T \sum_{k=1}^3 \sum_{l=1}^3 \sum_{j=1}^3 \left(D_{jk}(\partial_y,n_y) \vec{u}_l(y) \right)  \frac{1}{|x-y|} \left(D_{jl}(\partial_x,n_x) \vec{v}_k(x) \right) \,ds_x\, ds_y \nonumber \\
&  + \frac{\mu}{4 \pi}  \int_\Gamma\int_{\partial T} \sum_{k=1}^3 \sum_{l=1}^3 \sum_{j=1}^3 \frac{1}{|x-y|}  \vec{v}_k(x) t_{lj}(x) \left(D_{jk}(\partial_y,n_y) \vec{u}_l(y) \right) \,d\sigma_x \, ds_y
 \end{align}
for all $u \in H^1(\Gamma) $, $v \in \prod_{T \in \mathcal{T}_H} H^1(T)$ and with tangential matrix
\begin{align} \label{eq:tangential_Matrix}
 T_{kj}(x)= t_{kj}(x)=\begin{pmatrix}0 & -\vec{t}_3(x) & \vec{t}_2(x)\\ \vec{t}_3(x) & 0 & -\vec{t}_1(x)\\ -\vec{t}_2(x) & \vec{t}_1(x) & 0  \end{pmatrix}_{kj} .
\end{align}

\begin{proof}
We perform elementwise integration by parts for the three terms of $\langle Wu,v \rangle$, c.f.~\eqref{eq:han_W_rep}, individually. For the first summand we obtain with Lemma~\ref{lem:part_int_s_k_term} that
 \begin{align}
  & \quad -\frac{\mu}{4 \pi} \int_T \int_\Gamma \vec{v}(x) \sum_{k=1}^3 \frac{\partial}{\partial S_k(x)} \frac{1}{|x-y|} \frac{\partial \vec{u}(y)}{\partial S_k(y)} \, ds_y \, ds_x 
  = -\frac{\mu}{4 \pi} \int_\Gamma \int_T \sum_{k=1}^3 \sum_{l=1}^3 \vec{v}_l(x) \frac{\partial \vec{u}_l(y)}{\partial S_k(y)}  \frac{\partial}{\partial S_k(x)} \frac{1}{|x-y|}  \, ds_x \, ds_y \nonumber \\
  &= \frac{\mu}{4 \pi} \int_\Gamma \int_T \sum_{k=1}^3 \frac{1}{|x-y|}  \frac{\partial}{\partial S_k(x)} \left( \sum_{l=1}^3 \vec{v}_l(x) \frac{\partial \vec{u}_l(y)}{\partial S_k(y)} \right)   \, ds_x \, ds_y  - \frac{\mu}{4 \pi} \int_\Gamma \int_{\partial T} \sum_{k=1}^3 \frac{1}{|x-y|} \sum_{l=1}^3 \vec{v}_l(x) \frac{\partial u_l(y)}{\partial S_k(y)} \vec{t}_k(x) \, d\sigma_x \, ds_y \nonumber \\
  &= \frac{\mu}{4 \pi} \int_\Gamma \int_T  \frac{1}{|x-y|} \sum_{k=1}^3 \frac{\partial \vec{u}(y)}{\partial S_k(y)} \cdot \frac{\partial \vec{v}(x)}{\partial S_k(x)} \,ds_x\, ds_y -  \frac{\mu}{4 \pi} \int_\Gamma \int_{\partial T} \frac{1}{|x-y|} \sum_{k=1}^3 \vec{v}(x) \frac{\partial \vec{u}(y)}{\partial S_k(y)} \vec{t}_k(x) \, d\sigma_x \, ds_y
 \end{align}
since the bilinearity of the cross-product yields
\begin{align*}
 \frac{\partial}{\partial S_k(x)}\left( \sum_{l=1}^3  \vec{v}_l(x) \frac{\partial \vec{u}_l(y)}{\partial S_k(y)} \right)   &=\left[ \vec{n}(x) \times \nabla \left( \sum_{l=1}^3  \vec{v}_l(x) \frac{\partial \vec{u}_l(y)}{\partial S_k(y)} \right)  \right]_k  =  \sum_{l=1}^3 \frac{\partial \vec{u}_l(y)}{\partial S_k(y)} \left[ \vec{n}(x) \times \nabla v_l(x) \right]_k = \sum_{l=1}^3 \frac{\partial \vec{u}_l(y)}{\partial S_k(y)} \frac{\partial \vec{v}_l(x)}{\partial S_k(x)} \\
 &= \frac{\partial \vec{u}(y)}{\partial S_k(y)} \cdot \frac{\partial \vec{v}(x)}{\partial S_k(x)}.
\end{align*}
The second summand has the form (coming from the vector-matrix-matrix-vector multiplication)
\begin{align*}
 -\int_\Gamma \int_T \sum_{k=1}^3 \sum_{l=1}^3 \sum_{j=1}^3 \vec{v}_k(x) D_{kj}(\partial_x,n_x) B_{jl}(x,y) \left( D(\partial_y,n_y) \vec{u}(y) \right)_l \, ds_x\, ds_y
\end{align*}
with $B(x,y)=4 \mu^2 G(x,y) - \frac{\mu}{2\pi} I \frac{1}{|x-y|}$.
Assuming $D_{kj}(\partial_x,n_x) = \frac{\partial}{ \partial S_z(x)}$ it holds with Lemma~\ref{lem:part_int_s_k_term} that
\begin{align*}
 & \quad -\int_\Gamma \int_T \vec{v}_k(x) \frac{\partial}{ \partial S_z(x)} B_{jl}(x,y) \left( D(\partial_y,n_y) \vec{u}(y) \right)_l \, ds_x\, ds_y\\
 &= \int_\Gamma \int_T B_{jl}(x,y) \frac{\partial}{ \partial S_z(x)} \left[ \vec{v}_k(x) \left( D(\partial_y,n_y) \vec{u}(y) \right)_l \right] \, ds_x\, ds_y - \int_\Gamma \int_{\partial T} B_{jl}(x,y)  \left[ \vec{v}_k(x) \left( D(\partial_y,n_y) \vec{u}(y) \right)_l \right] \vec{t}_z(x) \, d\sigma_x\,ds_y \\
 &=-\int_\Gamma \int_T D_{jk}(\partial_x,n_x)\vec{v}_k(x) B_{jl}(x,y) \left( D(\partial_y,n_y) \vec{u}(y) \right)_l   \, ds_x\, ds_y - \int_\Gamma \int_{\partial T} \vec{v}_k(x) \vec{t}_z(x) B_{jl}(x,y)  \left( D(\partial_y,n_y) \vec{u}(y) \right)_l  \, d\sigma_x\,ds_y
\end{align*}
since
\begin{align*}
 \frac{\partial}{ \partial S_z(x)} & \left[ \vec{v}_k(x) \left( D(\partial_y,n_y) \vec{u}(y) \right)_l \right] = \left[n_x \times \nabla \vec{v}_k(x) \left( D(\partial_y,n_y) \vec{u}(y) \right)_l \right]_z = \left( D(\partial_y,n_y) \vec{u}(y)) \right)_l \left(n_x \times \nabla \vec{v}_k(x) \right)_z\\
 &=\left( D(\partial_y,n_y) \vec{u}(y) \right)_l \frac{\partial \vec{v}_k(x)}{\partial S_z(x)}
 =\left( D(\partial_y,n_y) \vec{u}(y) \right)_l D_{kj}(\partial_x,n_x) \vec{v}_k(x)
 =-\left( D(\partial_y,n_y) \vec{u}(y) \right)_l D_{jk}(\partial_x,n_x) \vec{v}_k(x).
\end{align*}
Hence,
\begin{align}
-\int_T &\int_\Gamma  v(x)D(\partial_x,n_x)B(x,y)D(\partial_y,n_y)u(y)\, ds_y \, ds_x \nonumber \\
 &= -\int_\Gamma \int_T \sum_{k=1}^3 \sum_{l=1}^3 \sum_{j=1}^3 \vec{v}_k(x) D_{kj}(\partial_x,n_x) B_{jl}(x,y) \left( D(\partial_y,n_y) \vec{u}(y) \right)_l \, ds_x\, ds_y \nonumber \\
 &= -\int_\Gamma \int_T \sum_{k=1}^3 \sum_{l=1}^3 \sum_{j=1}^3 D_{jk}(\partial_x,n_x)\vec{v}_k(x) B_{jl}(x,y) \left( D(\partial_y,n_y) \vec{u}(y) \right)_l   \, ds_x\, ds_y \nonumber \\
 & \quad - \int_\Gamma \int_{\partial T} \sum_{k=1}^3 \sum_{l=1}^3 \sum_{j=1}^3 \vec{v}_k(x) T_{kj}(x) B_{jl}(x,y)  \left( D(\partial_y,n_y) \vec{u}(y) \right)_l  \, d\sigma_x\,ds_y \nonumber \\
 &=- \int_\Gamma \int_T \left( D(\partial_x,n_x)\vec{v}(x) \right)^\top B(x,y) D(\partial_y,n_y)\vec{u}(y) \, ds_x\, ds_y - \int_\Gamma \int_{\partial T} \vec{v}^\top(x) T(x)
 B(x,y) D(\partial_y,n_y)\vec{u}(y) \, d\sigma_x\,ds_y
\end{align}
with $T_{kj}(x)$ defined in \eqref{eq:tangential_Matrix}.
For the third summand we get by arguing as above
\begin{align*}
 & \quad \int_T\vec{v}_k(x) D_{lj}(\partial_x,n_x) \frac{1}{|x-y|} D_{jk}(\partial_y,n_y) \vec{u}_l(y) \, ds_x \\
 & = -\int_T \left(D_{lj}(\partial_x,n_x) \vec{v}_k(x) \right) \frac{1}{|x-y|} \left(D_{jk}(\partial_y,n_y) \vec{u}_l(y) \right) \, ds_x + \int_{\partial T} \frac{1}{|x-y|}  \vec{v}_k(x) t_{lj}(x) \left(D_{jk}(\partial_y,n_y) \vec{u}_l(y) \right) \,d\sigma_x.
\end{align*}
Thus, by the anti-symmetry of $D$ we obtain
\begin{align}
 \frac{\mu}{4 \pi} \int_\Gamma \int_T & \vec{v}(x) \left( D(\partial_x,n_x) \frac{1}{|x-y|} D(\partial_y,n_y) \right)^\top u(y) \,ds_x \, ds_y \nonumber \\
& = \frac{\mu}{4 \pi} \int_\Gamma \int_T \sum_{k=1}^3 \sum_{l=1}^3 \sum_{j=1}^3 \left(D_{jk}(\partial_y,n_y) \vec{u}_l(y) \right)  \frac{1}{|x-y|} \left(D_{jl}(\partial_x,n_x) \vec{v}_k(x) \right) \,ds_x\, ds_y \nonumber \\
& \quad + \frac{\mu}{4 \pi}  \int_\Gamma\int_{\partial T} \sum_{k=1}^3 \sum_{l=1}^3 \sum_{j=1}^3 \frac{1}{|x-y|}  \vec{v}_k(x) t_{lj}(x) \left(D_{jk}(\partial_y,n_y) \vec{u}_l(y) \right) \,d\sigma_x \, ds_y.
\end{align}
Combing the rearrangement of the three summands and summing over all elements $T \in \mathcal{T}_{H}$ yields the assertion.
\end{proof}
\end{theorem}

\begin{remark}
 If the test function $v(x)$ is continuous and the mesh $\mathcal{T}_{H}$ describes a closed surface $\Gamma_C$ or $v|_{\partial \Gamma_C}=0$, the elementwise integration by parts of $W$ in Theorem~\ref{thm:elementPartIntW} coincides with global integration by parts of $W$ in \cite[Eq.~3.12]{han1994boundary}.
\end{remark}

For the computations it is convenient to represent the double layer potential $K$ as a sum of the harmonic double layer potential, the harmonic single layer potential and the single layer potential, c.f.~e.g.~\cite[Lem.~2.2]{han1994boundary}. As that representation is only valid on a closed boundary with continuous ansatz functions, we redo the arguments of Theorem~\ref{thm:elementPartIntW} for the double layer potential to be able to implement the dual of $\widehat{K^\prime}$.

\begin{lemma}
 For $u \in \prod_{T \in \mathcal{T}_H} H^1(T)$ there holds
 \begin{align}
  Ku(x)&=\sum_{T \in \mathcal{T}_H} \int_T u(y) \frac{\partial}{\partial n_y} \frac{1}{4\pi |x-y|}\, ds_y 
  -\int_T \frac{1}{4\pi |x-y|} D(\partial_y,n_y)u(y) \, ds_y + \int_{\partial T} \frac{1}{4\pi |x-y|} T(y)  \vec{u}(y) \, d\sigma_y \nonumber \\
  & \quad +2\mu \int_T G(x,y) D(\partial_y,n_y)u(y) \, ds_y -2\mu \int_{\partial T}  G(x-y) T(y)  \vec{u}(y) \, d\sigma_y.
 \end{align}
\end{lemma}

\section{Numerical Results}
\label{sec:Numerics}
Common to all numerical results are: We choose $h=k$ and $p=q$, i.e.~a pairing which does not satisfy the discrete inf-sup condition, and we choose $H=h$ and $P=p$ in \eqref{eq:StabApproxSpace} for the approximation of the stabilization term. For the adaptive scheme we use the Algorithm~25 in \cite{Banz2015Stab}. Namely, we use D\"orfler marking with the parameter $\theta\in (0,1)$ to decide which elements are to be refined, i.e.~those $T \in \mathcal{M}_h \subset \mathcal{T}_h$ such that $ \sum_{T \in \mathcal{M}_h} \eta_T^2 \geq \theta \sum_{T \in \mathcal{T}_h} \eta_T^2$ and $\mathcal{M}_h$ has minimal cardinality. Based on the decay rate of the Legendre expansion coefficient \cite{houston2005note} with the parameter $\delta \in (0,1)$, we decide between $h$- or $p$-refinement of the marked element.
We solve the discrete problem \eqref{eq:DiscreteProblem} with a semi-smooth Newton method, for which we write \eqref{eq:DiscreteContactConstraints} as the projection problem \eqref{eq:projection_problem}. We stop the Newton scheme once the square root of the merit function is less than $10^{-12}$. The Coulomb frictional problem is solved analogously. 

\subsection{Two dimensional Tresca frictional contact problem}
Within this section we present numerical results on the reduction of the approximation error for a uniform $h$, $h$-adaptive and $hp$-adaptive scheme, as well as on the error induced by approximating the stabilization term and on the sensitivity of the error estimate on the scaling constant $\bar{\gamma}$. The considered problem is the same as in \cite{Banz2015Stab}, that is: The domain is $\Omega=[-\frac{1}{2},\frac{1}{2}]^2$ with $\Gamma_C=[-\frac{1}{2},\frac{1}{2}]\times \left\{-\frac{1}{2}\right\}$, $\Gamma_D=[\frac{1}{4},\frac{1}{2}]\times \left\{\frac{1}{2}\right\}$ and $\Gamma_N=\partial \Omega \setminus \left(\Gamma_C \cup \Gamma_D\right)$. The material parameters are $E=500$ and $\nu=0.33$, the gap function is $g=1-\sqrt{1-0.01x_1^2}$ and the Tresca friction function is $\mathcal{F}=0.211+0.412x_1$. The Neumann force $f$ is zero except for
\begin{align*}
 f=  \left(
 \begin{array}{c}
 	 -(\frac{1}{2}-x_2)(-\frac{1}{2}-x_2)\\
 	0
 \end{array}\right) \text{ on } \left\{-\frac{1}{2}\right\}\times \left[-\frac{1}{2},\frac{1}{2}\right],
 \qquad
 f=  \left(
 \begin{array}{c}
 	0\\
 	20(-\frac{1}{2}-x_1)(-\frac{1}{4}-x_1)
 \end{array}\right) \text{ on } \left[-\frac{1}{2},-\frac{1}{4}\right]\times \left\{\frac{1}{2}\right\}
\end{align*}
and the stabilization scaling constant is $\bar{\gamma}=10^{-3}$ unless stated otherwise.
Similar examples can be found in \cite{schroder2012posteriori} for FEM and in \cite{Banz2014BEM} for BEM with biorthogonal basis functions. The solution is characterized by two singular points at the interface from Neumann to Dirichlet boundary condition. These singularities are more severe than the loss of regularity from the contact conditions. At the contact boundary the solution has a long interval in which it is sliding, i.e.~where $|\sigma_t|=\mathcal{F}$ and $u_t=-\alpha \sigma_t$ for some $\alpha \geq 0$, and in which the absolute value of the tangential Lagrange multiplier increases linearly like $\mathcal{F}$. The actual contact set, i.e.~where $u_n=g$, is centered slightly to the right of the midpoint of $\Gamma_C$.\\
 
\begin{figure}[tbp]
   \centering \mbox{
   \subfigure[a posteriori error estimation]{
   \includegraphics[trim = 15mm 2mm 19mm 8mm, clip,width=60.0mm, keepaspectratio]{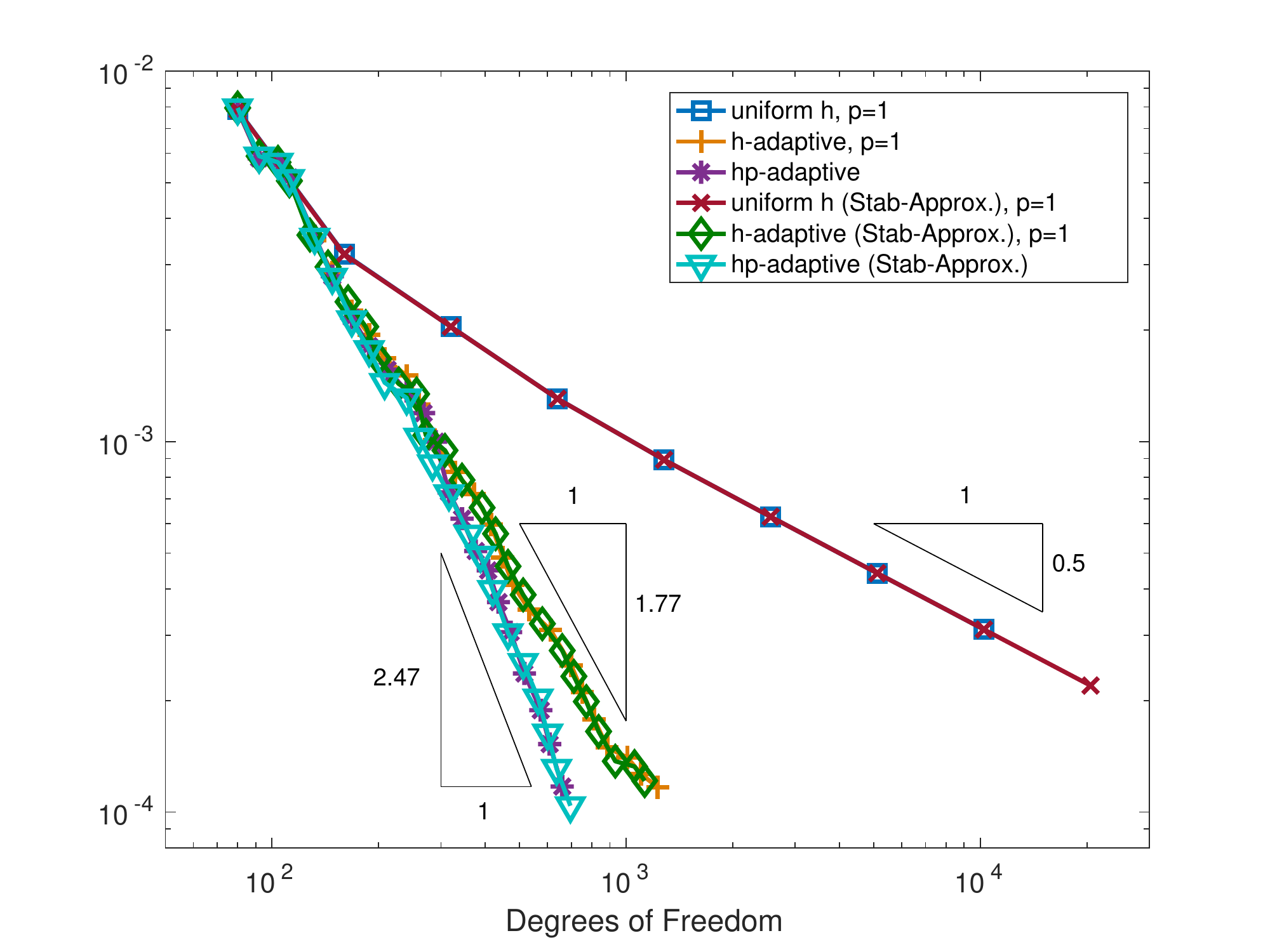}}
   \quad
      \subfigure[approximate error in energy norm]{
   \includegraphics[trim = 15mm 2mm 19mm 8mm, clip,width=60.0mm, keepaspectratio]{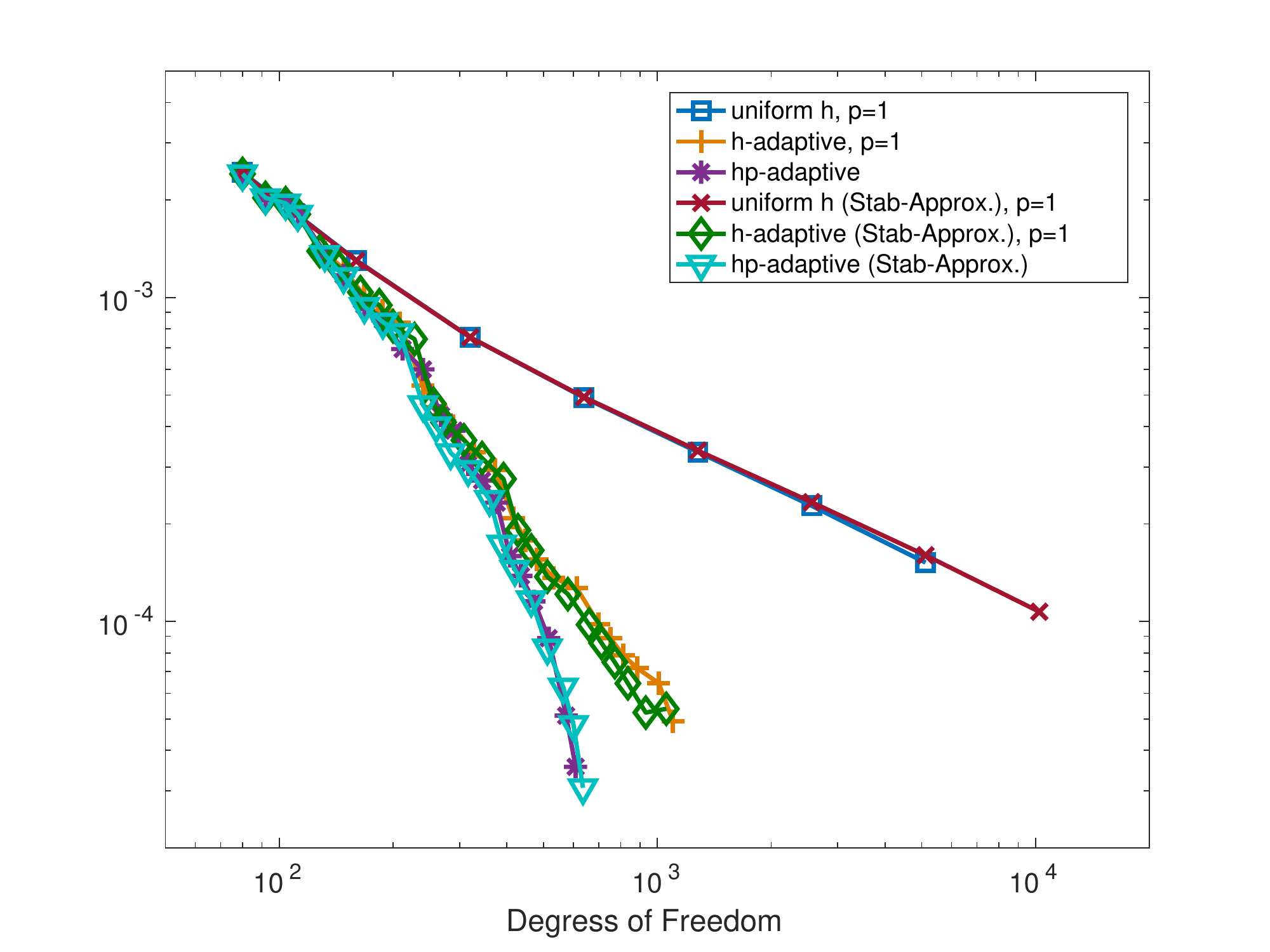}}
   }
   \caption{Error estimates for different families of discrete solutions (Tresca-friction case).}
   \label{fig:Tresca:errorDiri}
 \end{figure}

\begin{figure}[tbp]
   \centering \mbox{
   \subfigure[$h$-adaptive]{   
  \begin{overpic}[trim = 280mm 74mm 250mm 51mm,clip,width=50.0mm, keepaspectratio]{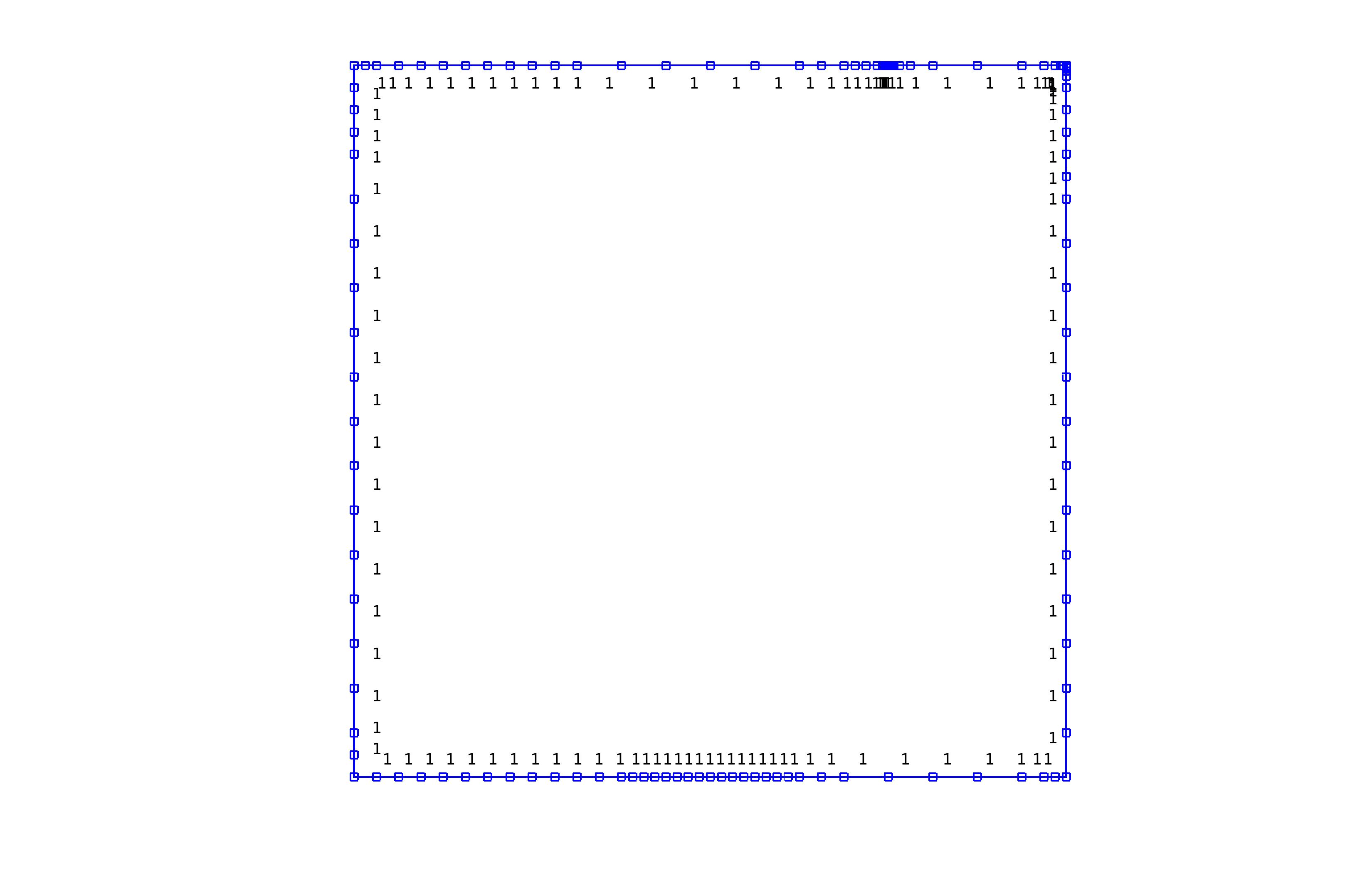}
     \put(8.5,10.5){ \includegraphics[trim = 280mm 74mm 250mm 51mm, clip, width=40.0mm, keepaspectratio]{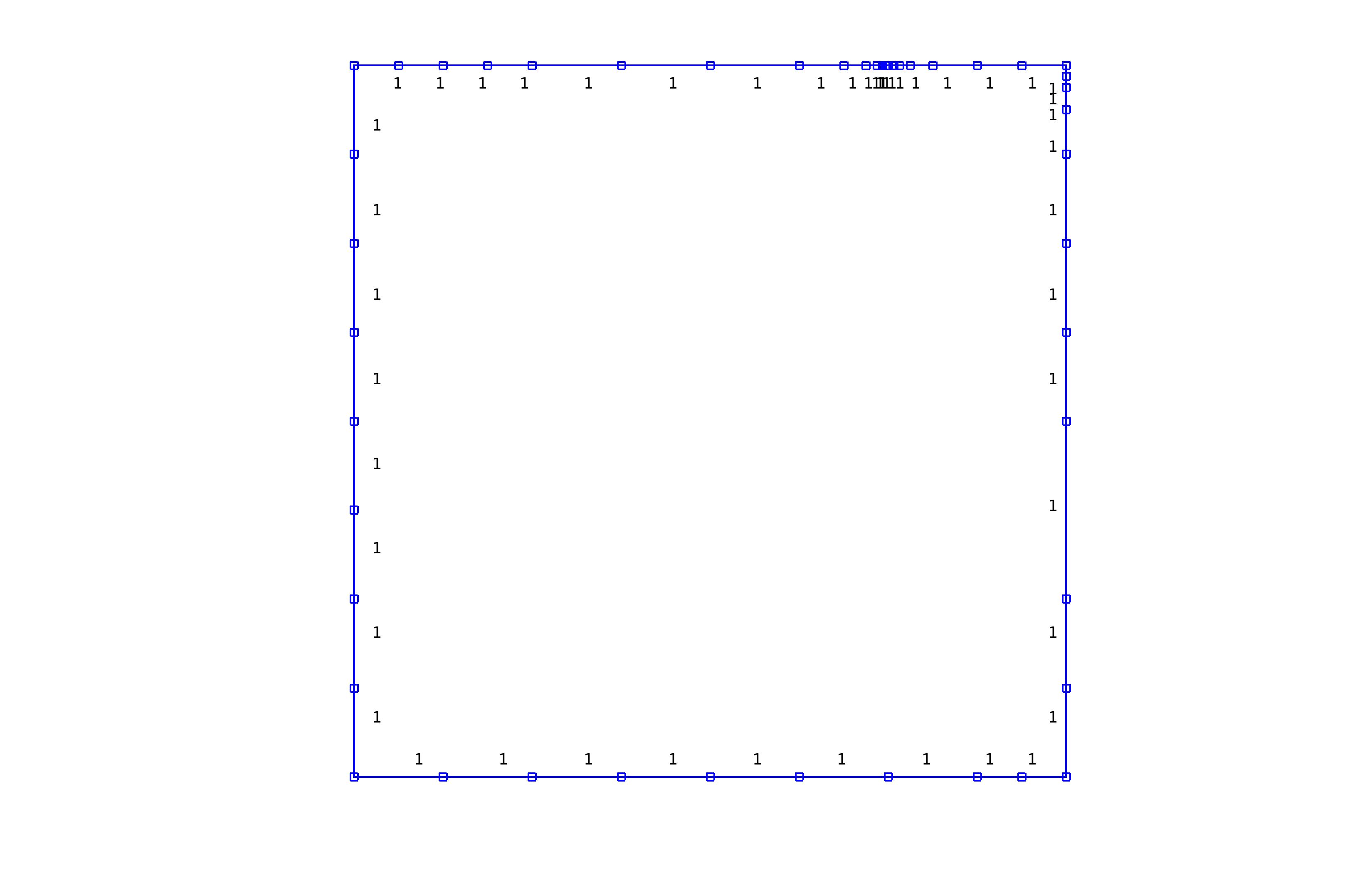}  }
   \end{overpic}   
 	} \qquad 
  \subfigure[$hp$-adaptive]{   
  \begin{overpic}[trim = 280mm 74mm 250mm 51mm,clip,width=50.0mm, keepaspectratio]{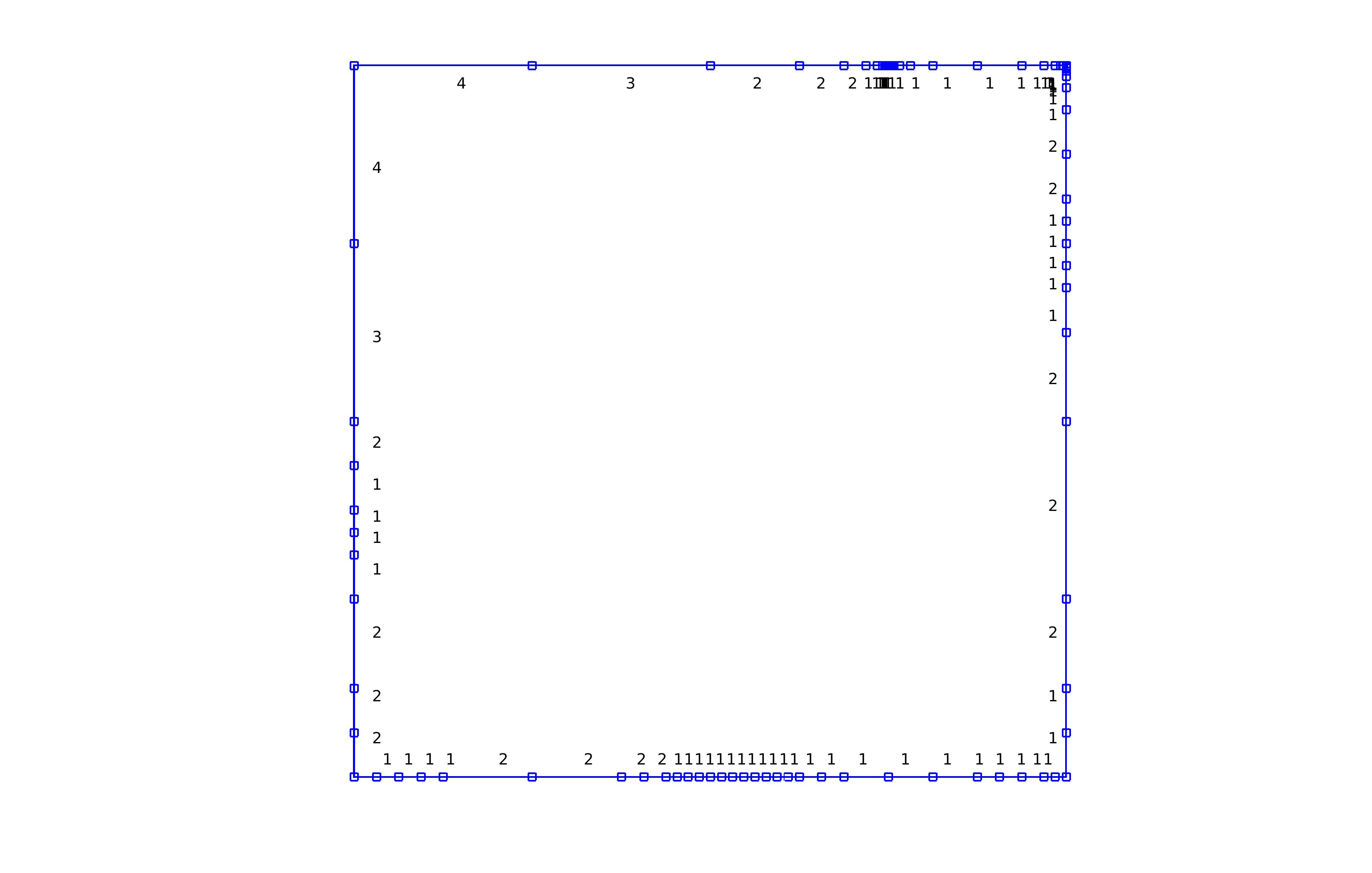}
     \put(8.5,10.5){ \includegraphics[trim = 280mm 74mm 250mm 51mm, clip, width=40.0mm, keepaspectratio]{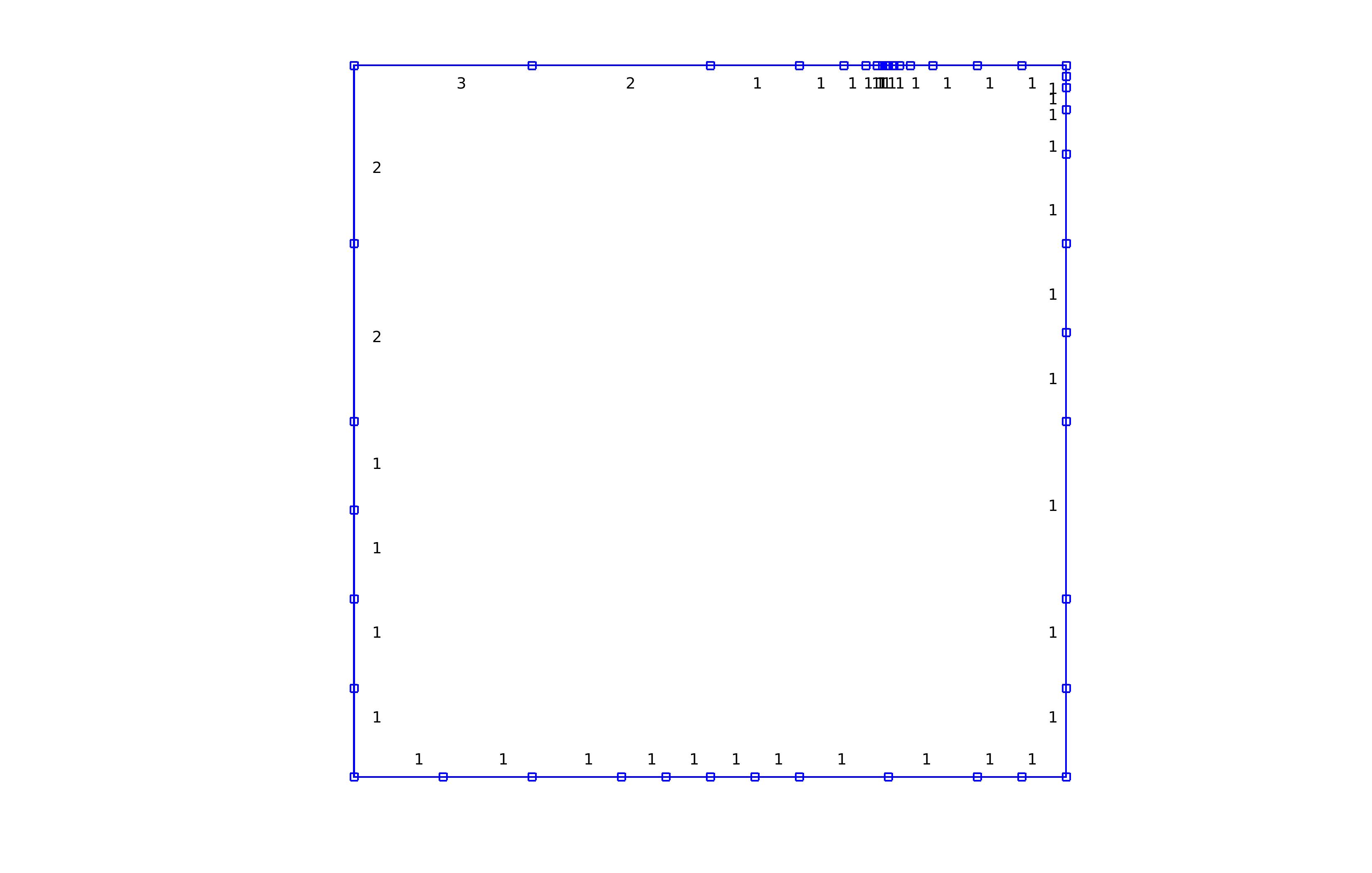}  }
   \end{overpic}   }
   }
   \caption{Adaptively generated meshes (Tresca-friction case), mesh nr. 10 (inner), nr. 20 (outer).}
   \label{fig:Tresca:AdaptiveMeshes}
 \end{figure}

\begin{figure}[tbp]
   \centering \mbox{\!
   \subfigure[uniform $h$-version with $p=1$]{
 	\includegraphics[trim = 13mm 2mm 15mm 8mm, clip,width=53.0mm, keepaspectratio]{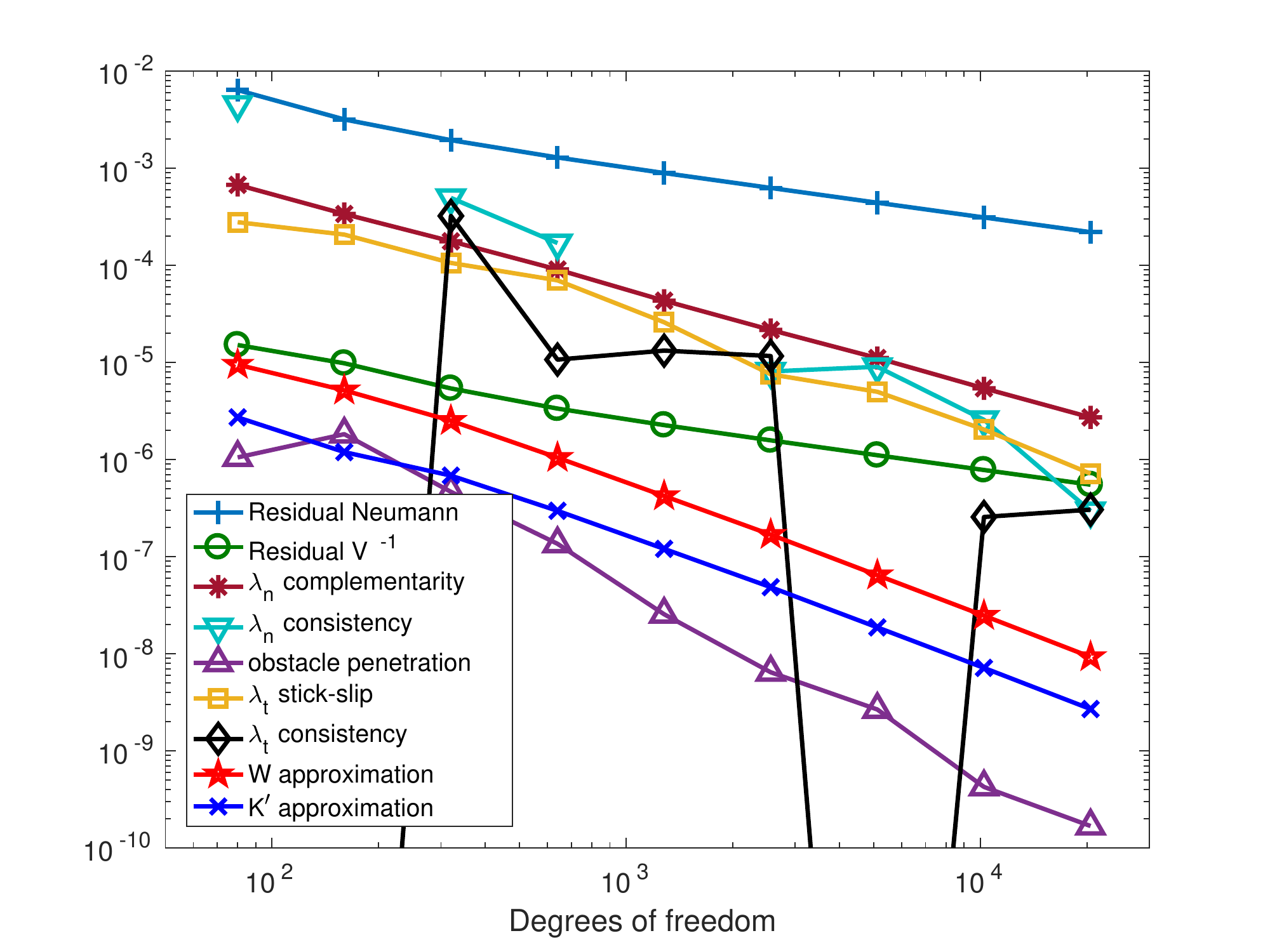}} 
   \subfigure[$h$-adaptive with $p=1$]{
 	\includegraphics[trim = 13mm 2mm 15mm 8mm, clip,width=53.0mm, keepaspectratio]{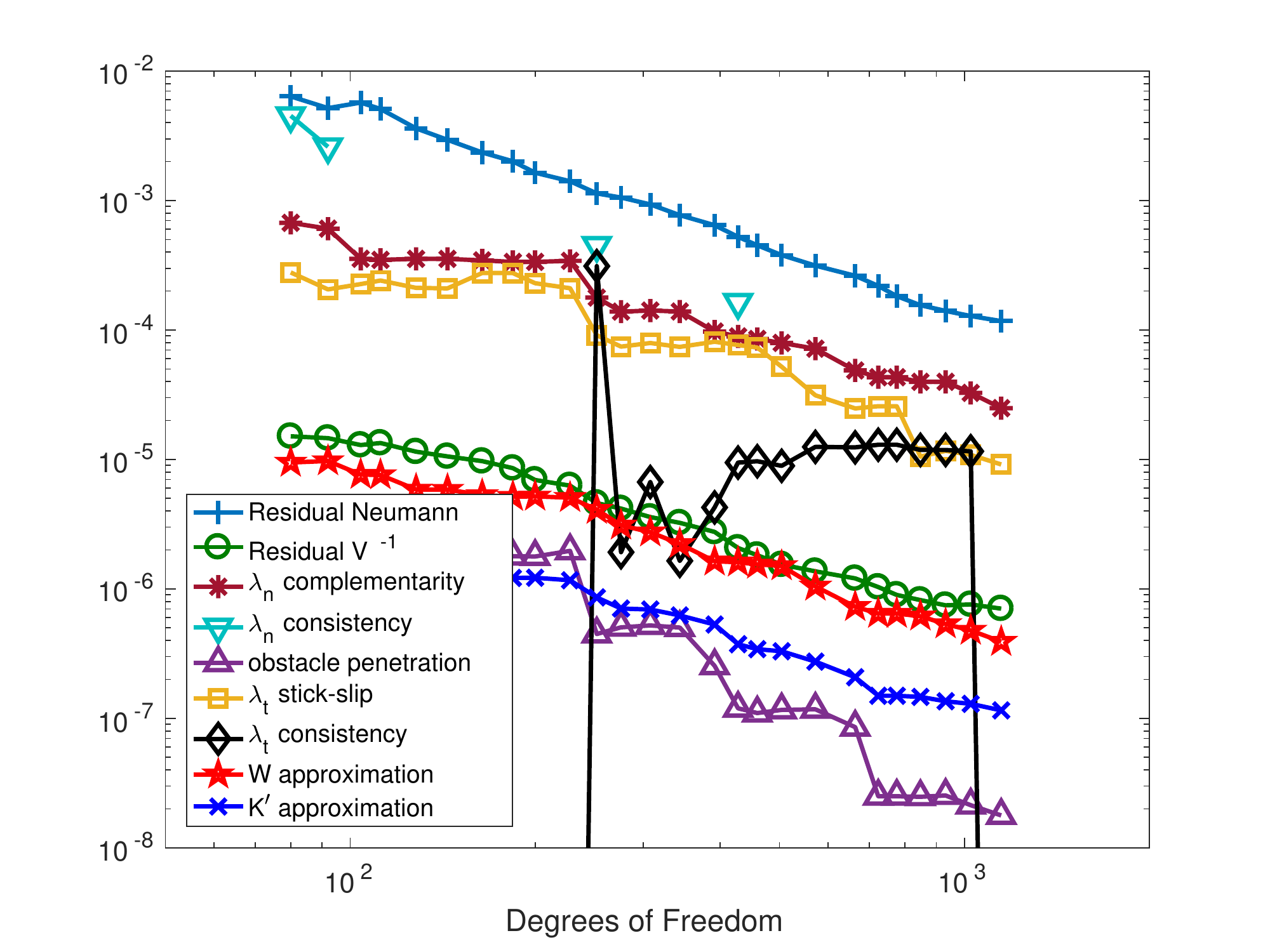}}  
	\subfigure[$hp$-adaptive]{
 	\includegraphics[trim = 13mm 2mm 15mm 8mm, clip,width=53.0mm, keepaspectratio]{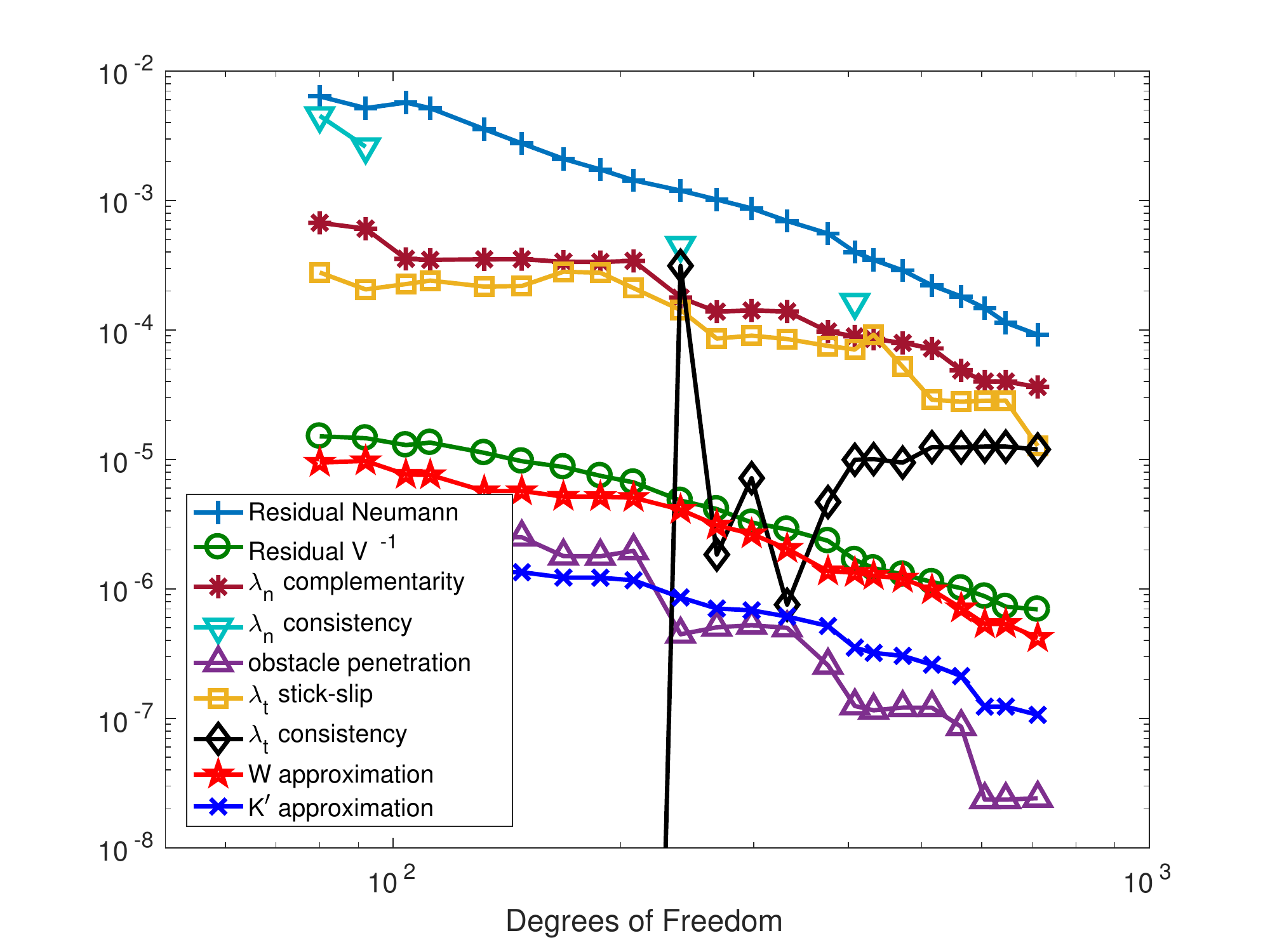}}
	}	
   \caption{Error contributions of the residual based a posteriori error estimate (Tresca-friction case).}
   \label{fig:Tresca:errorContributions}
 \end{figure}

\begin{figure}[tbp]
    \centering
    \includegraphics[trim = 15mm 2mm 19mm 8mm, clip,width=60.0mm, keepaspectratio]{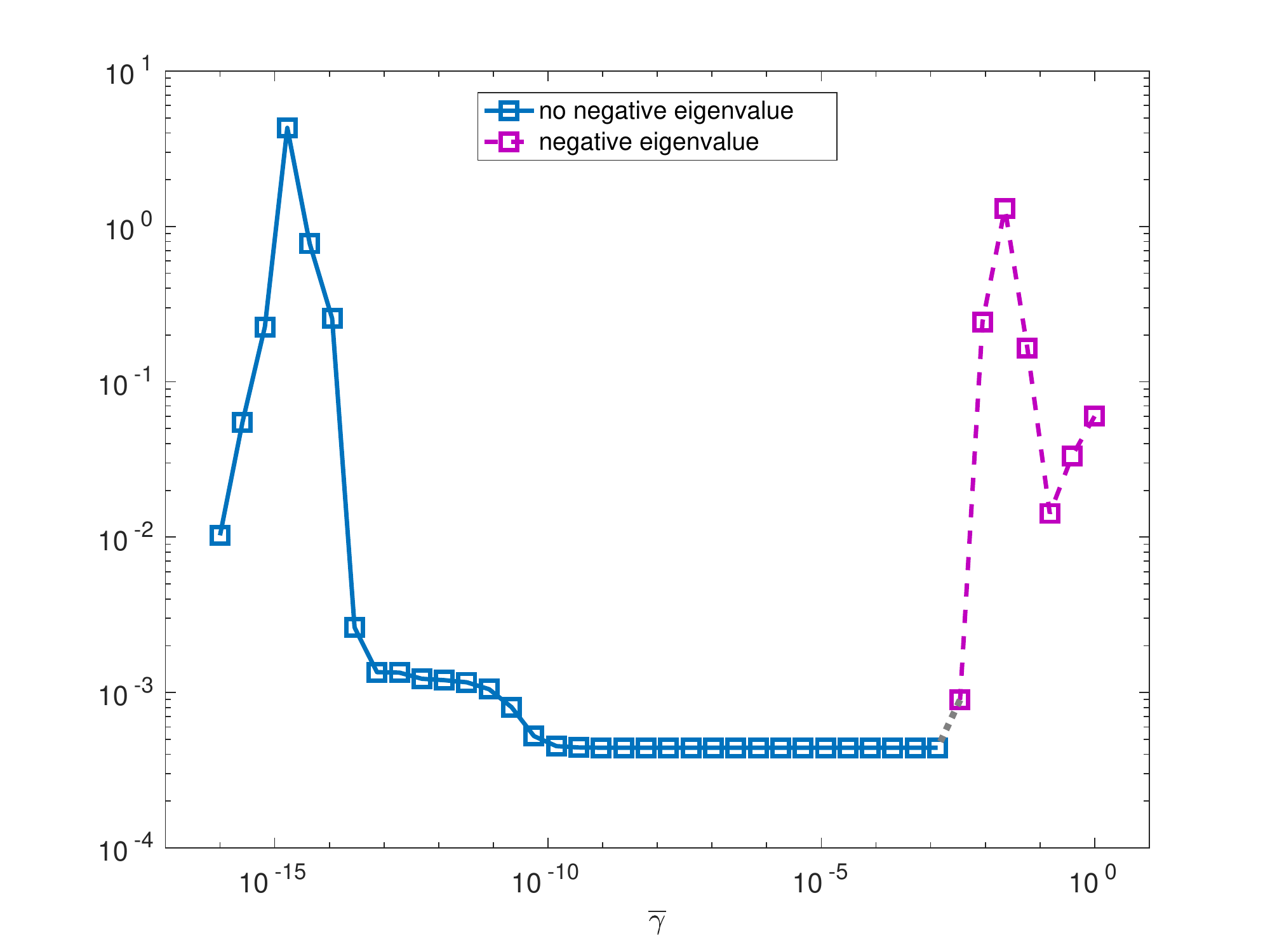}
    \caption{Dependency of the a posteriori error estimate on $\bar{\gamma}$ for a uniform mesh with 256 elements and $p=1$ (Tresca-friction case).}
    \label{fig:gamma_tresca}
  \end{figure}

 Figure~\ref{fig:Tresca:errorDiri} (a) displays the reduction of the a posteriori error estimate with an increase of the degrees of freedom for the uniform $h$-version with $p=1$, $h$-adaptive scheme with $p=1$ and $\theta=0.3$, as well as for the $hp$-adaptive scheme with $\theta=0.33$ and $\delta = 0.5$. Figure~\ref{fig:Tresca:errorDiri} (b) displays the same information but for the approximate error. For that we replace the error with the approximation
 \begin{align*}
  \left( \| u - u_N\|_W^2 + \| \phi - \phi_N \|_V^2 + \| hp^{-2} (\lambda - \lambda_N)\|_{L^2(\Gamma_C)}^2\right)^{1/2} \approx  \left( \| u_{\text{fine}} - u_N\|_W^2 + \| \phi_{\text{fine}} - \phi_N \|_V^2 + \| hp^{-2} (\lambda_{\text{fine}} - \lambda_N)\|_{L^2(\Gamma_C)}^2\right)^{1/2}
 \end{align*}
 in which $(u_{\text{fine}},\phi_{\text{fine}},\lambda_{\text{fine}})$ is the finest, last solution of the respective sequence of discrete solutions.
 Figure~\ref{fig:Tresca:errorDiri} indicates, that the error induced by approximating the stabilization term is of several magnitudes smaller than the approximation error, which is consistent with the finding in \cite{Banz2015Stab}. In particular, even the error curves for the $hp$-adaptive schemes are almost the same, but the sequence of iterates itself is not. The experimentally determined efficiency indices (error estimate divided by error) is between two and four (we have divided $\eta_{N,T}^2$ and $\eta_{V,T}^2$ by 100 to have roughly the same absolute value as if we were to apply a bubble error indicator). The computed convergence rates are $0.5$ for the uniform $h$-version with $p=1$, $1.77$ for the $h$-adaptive scheme which is expected to drop to 1.5 asymptotically, and $2.47$ for the $hp$-adaptive scheme.\\
 In both adaptive cases the transition from Neumann to Dirichlet boundary condition as well as the contact zone and the stick-slip transition are resolved by mesh refinements. Contrary to $h$-adaptivity, the $hp$-adaptive scheme increases the polynomial degree moving away from those singularities, c.f.~Figure~\ref{fig:Tresca:AdaptiveMeshes}.\\
 The for the adaptivity responsible a posteriori error estimate is decomposed into nine error contributions, as depicted in Figure~\ref{fig:Tresca:errorContributions}. In all cases the classical residual error contribution is the dominate source of error, followed equally by the violation of the complementarity condition and stick-slip condition. The error contributions associated with the approximation of the stabilization term are of several orders of magnitudes smaller the three largest error contributions. The consistency contribution of $\lambda$ are highly fluctuate, more on that in the next section.\\
 The drawback of the stabilization method is that the user has to choose the stabilization scaling parameter $\bar{\gamma}$ which must be sufficiently small compared to the ratio of ellipticity constant and polynomial inverse estimate constant. Our computations in Figure~\ref{fig:gamma_tresca} indicate that the ``best`` choice range for $\bar{\gamma}$ is with $3.7 \cdot 10^{-10} \leq \bar{\gamma} \leq 1.3 \cdot 10^{-3}$ very large. If it is chosen too large, negative eigenvalues occur, and if it is chosen too small, the stabilization is effectively switched off by the inaccurate computations of computers. In both cases the a posteriori error estimate increases dramatically.

\subsection{Two dimensional Coulomb frictional contact problem}
For the following numerical experiments taken from \cite{Banz2015Stab}, see also \cite{hueber2012equilibration} for FEM results of a very similar problem, the domain is $\Omega=[-\frac{1}{2},\frac{1}{2}]^2$ with $\Gamma_C=[-\frac{1}{2},\frac{1}{2}]\times \left\{-\frac{1}{2}\right\}$ and $\Gamma_N=\partial \Omega \setminus \Gamma_C$. Since no Dirichlet boundary is prescribed, the kernel of the Steklov operator consists of the three rigid body motions $\operatorname{ker}(S)=\operatorname{span}\left\{ (x_1,0)^\top, (0,x_2)^\top, (x_2,-x_1)^\top \right\}$. To obtain a unique solution nevertheless, we fix the point $(0,-1/2)$ as the problem is symmetric. Now, the contact conditions prevent rotations of the body $\Omega$. The material parameters are $E=5$ and $\nu=0.45$, and the Coulomb friction coefficient is $0.3$. The Neumann force is
\begin{alignat*}{2}
 f&=  \left(
 \begin{array}{c}
 	-10\operatorname{sign}(x_1)(\frac{1}{2}+x_2)(\frac{1}{2}-x_2)\exp(-10(x_2+\frac{4}{10})^2)\\
 	\frac{7}{8}(\frac{1}{2}+x_2)(\frac{1}{2}-x_2)
 \end{array}\right) &\quad \text{ on }& \left\{\pm\frac{1}{2}\right\}\times \left[-\frac{1}{2},\frac{1}{2}\right],
 \\
 f&=  \left(
 \begin{array}{c}
 	0\\
 	-\frac{25}{2}(\frac{1}{2}-x_1)^2(\frac{1}{2}+x_1)^2
 \end{array}\right) &\quad \text{ on }& \left[-\frac{1}{2},\frac{1}{2}\right]\times \left\{\frac{1}{2}\right\},
\end{alignat*}
the gap to the obstacle is zero and $\bar{\gamma}=10^{-3}$.\\
  
\begin{figure}[tbp]
   \centering \mbox{
   \subfigure[a posteriori error estimation]{
 	\includegraphics[trim = 15mm 2mm 15mm 8mm, clip,width=60.0mm, keepaspectratio]{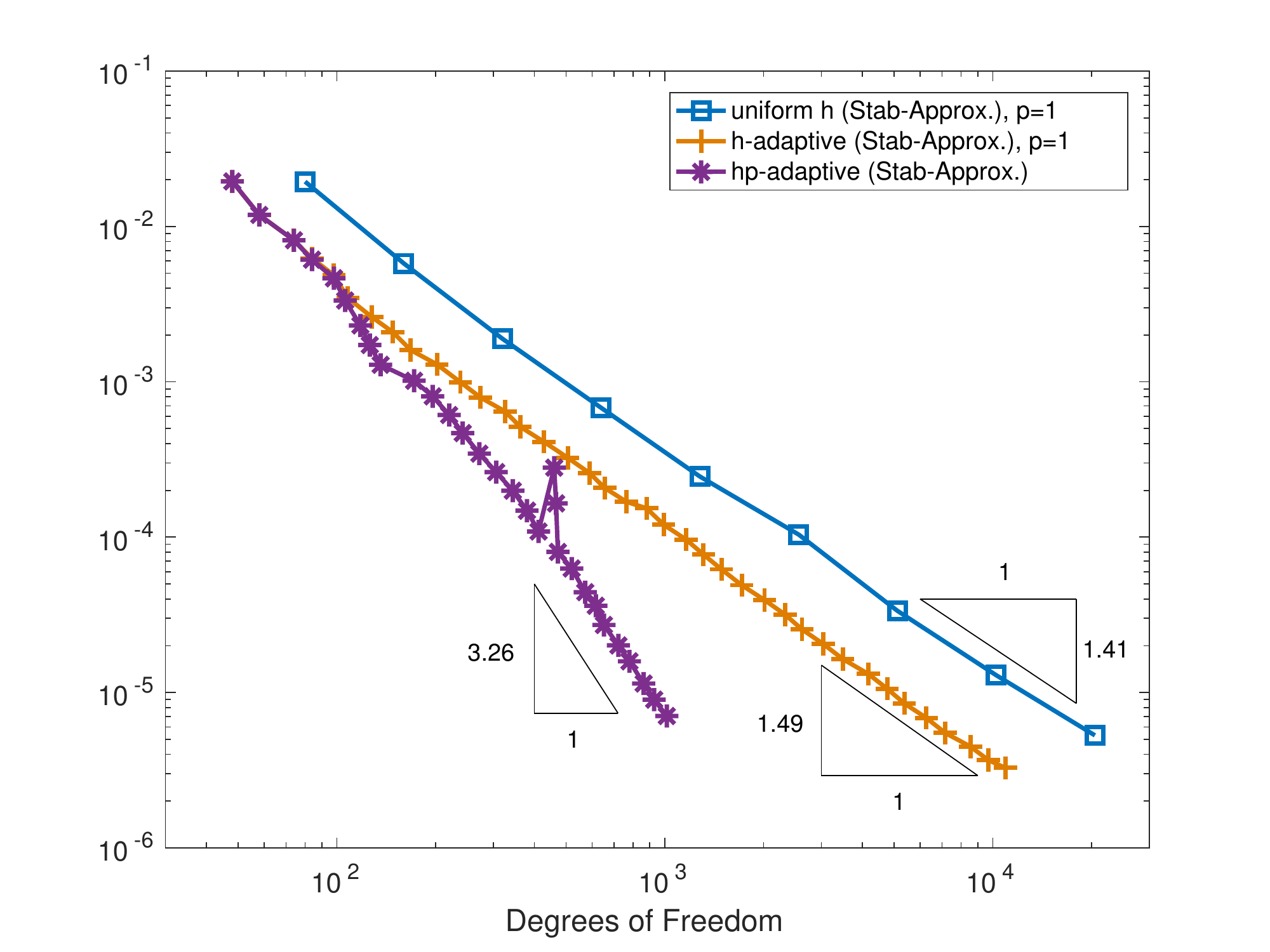}}  \quad 
   \subfigure[approximate error in energy norm]{
 	\includegraphics[trim = 15mm 2mm 15mm 8mm, clip,width=60.0mm, keepaspectratio]{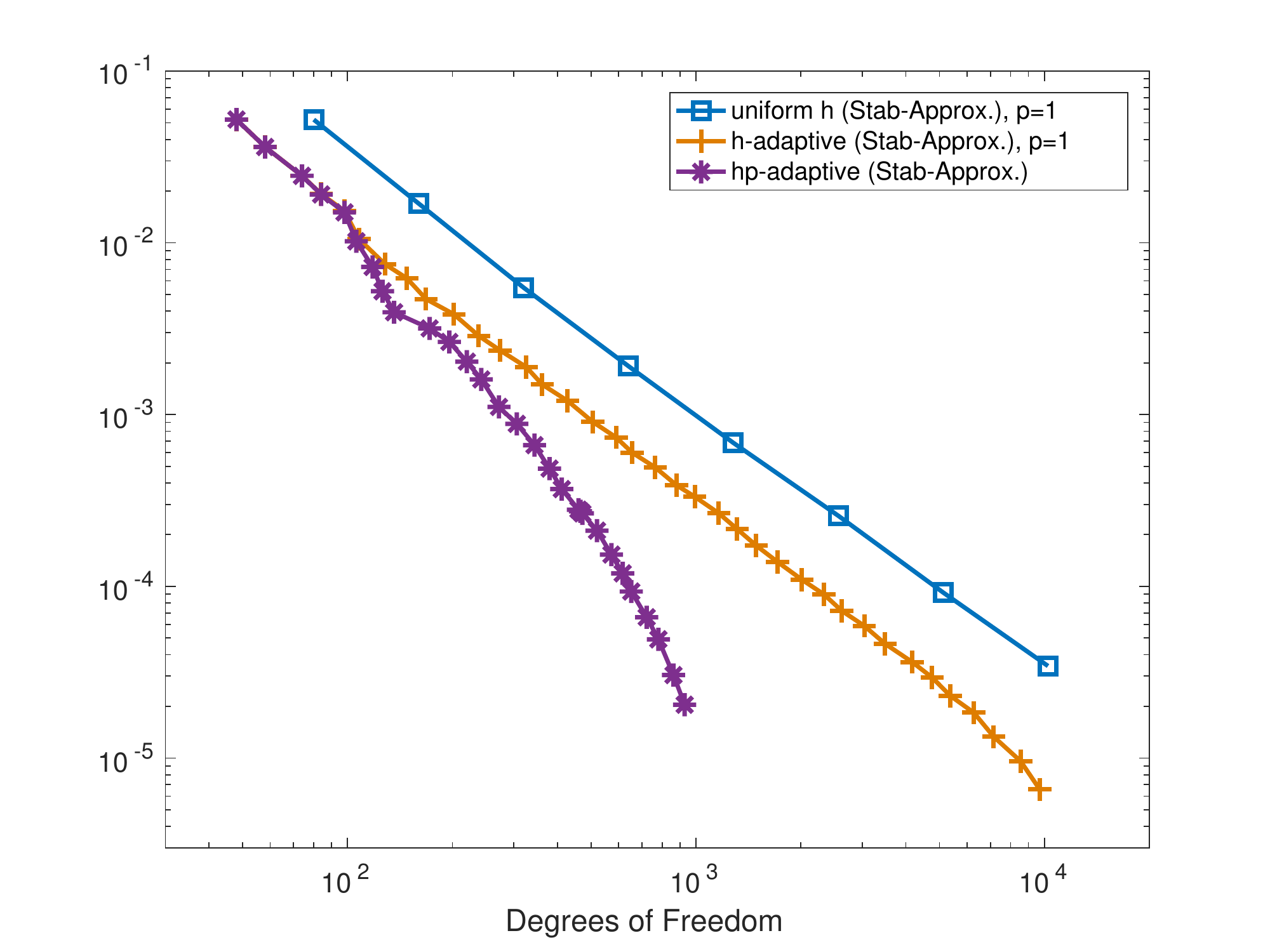}}  }
   \caption{Error estimates for different families of discrete solutions (Coulomb-friction case).}
   \label{fig:Coulomb:errorEstimates}
 \end{figure}

\begin{figure}[tbp]
   \centering \mbox{\!
   \subfigure[uniform $h$-version with $p=1$]{
 	\includegraphics[trim = 13mm 2mm 15mm 8mm, clip,width=53.0mm, keepaspectratio]{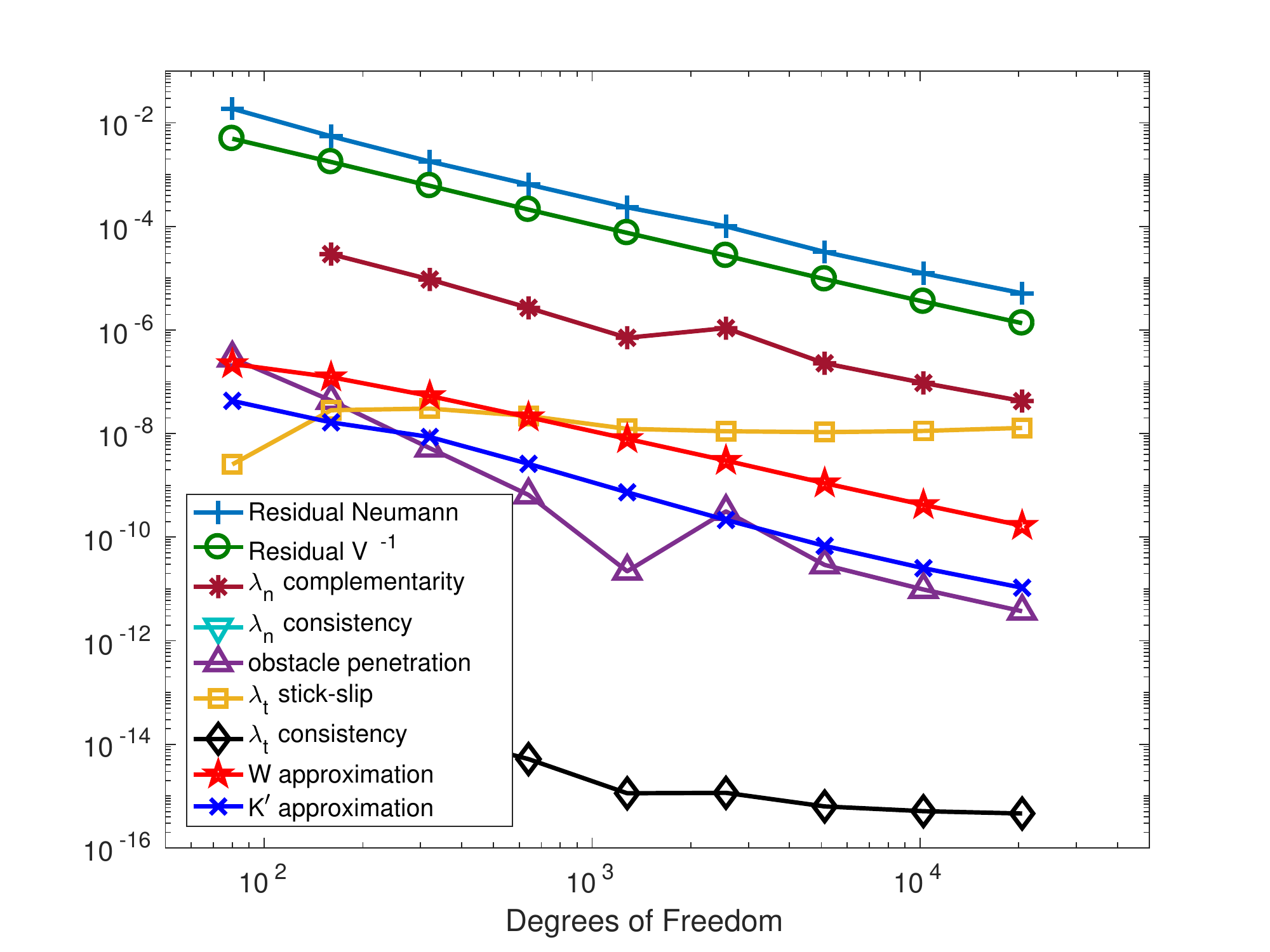}}
   \subfigure[$h$-adaptive with $p=1$]{
 	\includegraphics[trim = 13mm 2mm 15mm 8mm, clip,width=53.0mm, keepaspectratio]{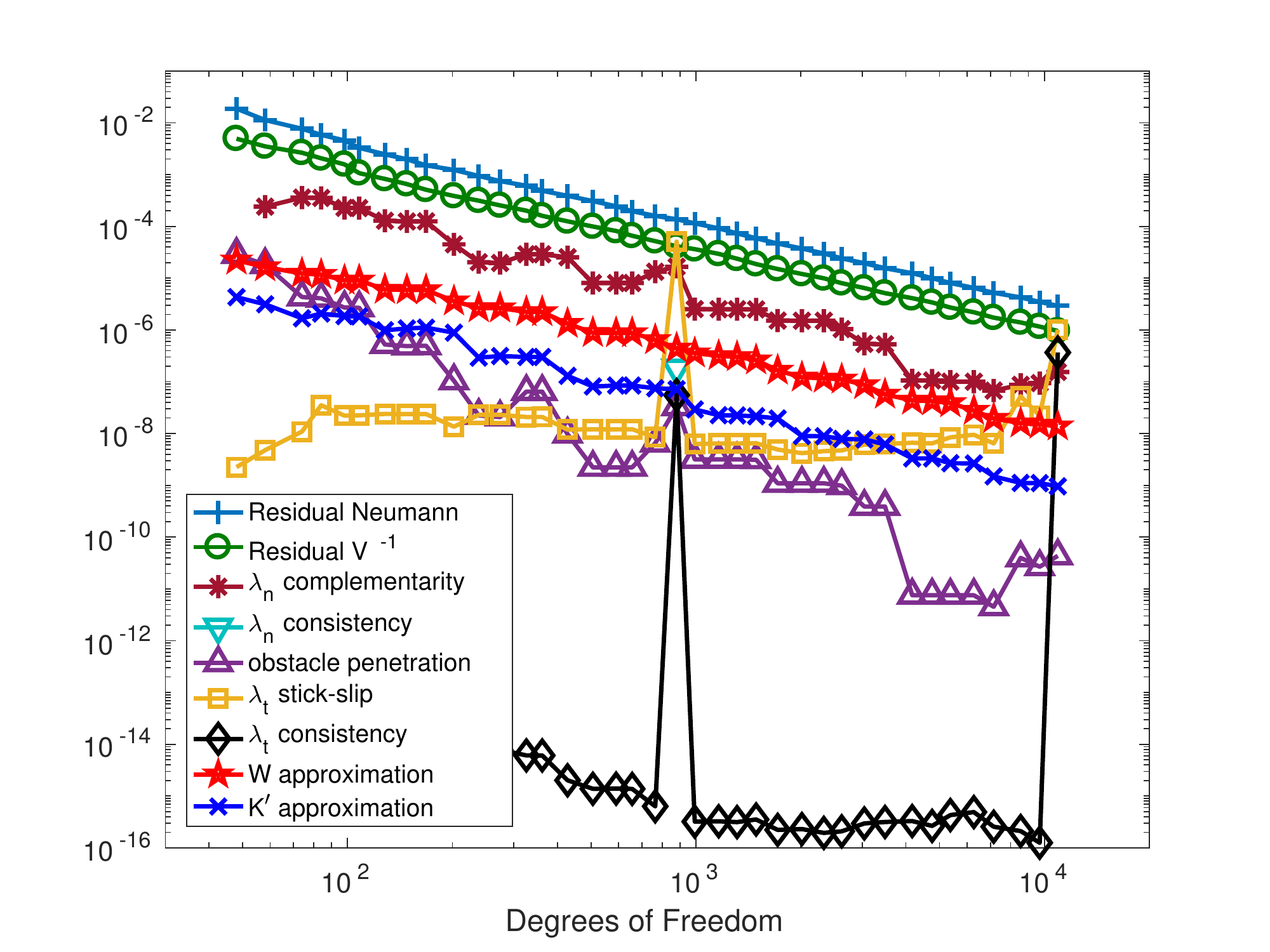}}  	
  \subfigure[$hp$-adaptive]{
 	\includegraphics[trim = 13mm 2mm 15mm 8mm, clip,width=53.0mm, keepaspectratio]{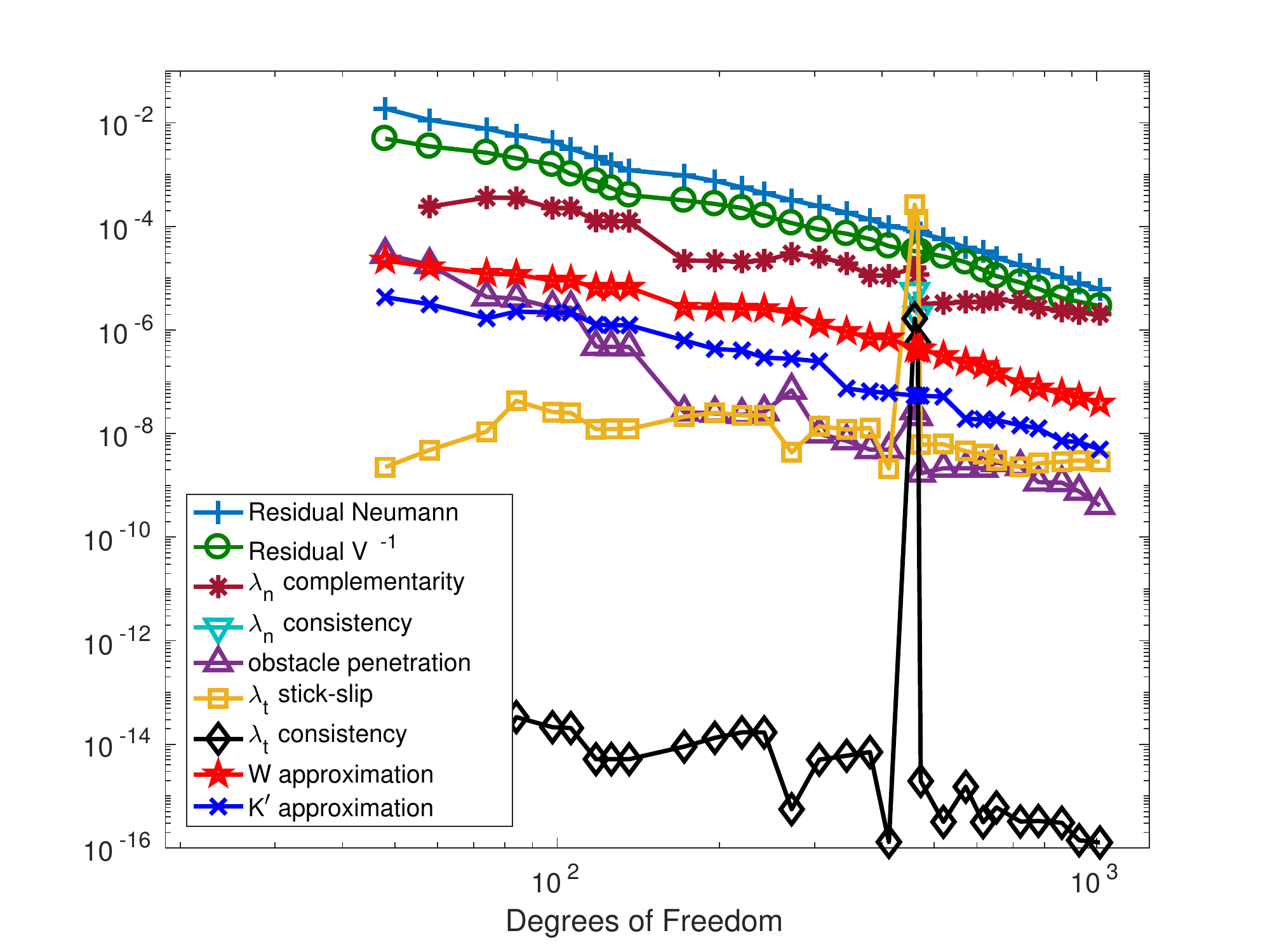}} } 
		
   \caption{Error contributions of the residual based a posteriori error estimate (Coulomb-friction case).}
   \label{fig:Coulomb:errorContributions}
 \end{figure}

Figure~\ref{fig:Coulomb:errorEstimates} shows the reduction of the error estimate and of the approximate error as in the previous section. The efficiency indices range between 2.6 and 3. The convergence rates are $1.41$ for the uniform $h$-version with $p=1$, (optimal) $1.49$ for the $h$-adaptive scheme with $p=1$ and $\theta=0.4$, and $3.26$ for the $hp$-adaptive scheme with $\theta=0.4$ and $\delta=0.5$. The curve of the a posteriori error estimate for the $hp$-adaptive scheme shows an upwards jump for the duration of two refinements, which is not there in the approximate error curve. That jump is a result of the stick-slip and consistency contribution of $\lambda_t$ in the a posteriori error estimate, visuable in Figure~\ref{fig:Coulomb:errorContributions}. That extreme behavior does not seem to be a result of an inaccurate solution of the discrete problem (decreasing the tolerance and changing the iterative solver itself did not affect the upward jumps, and also the residual Neumann contribution of the a posteriori error estimate is not effected), but more a non-conformity problem of the discretization method itself. Enforcing the sign-, box-constraints of $\lambda_n$, $\lambda_t$, respectively, in the Gauss-Legendre quadrature nodes is very handy for proving convergence rates, but leads to a non-conforming discretization even for $p=1$. If such a constraint violation occurs for the discrete solution, it occurs on a subinterval with positive measure, and is thus picked up strongly by the a posteriori error. These peaks might be avoidable if Bernstein polynomials are used for $\lambda_{hp}$ as suggested in \cite{Banz2015Stab}.\\
The adaptively generated meshes show typical refinement pattern and are omitted for brevity. 

\subsection{Three dimensional Coulomb frictional contact problem}

\begin{figure}[tbp]
   \centering \mbox{
   \subfigure[deformed body]{
   \includegraphics[trim = 30mm 5mm 40mm 13mm, clip,width=50.0mm, keepaspectratio]{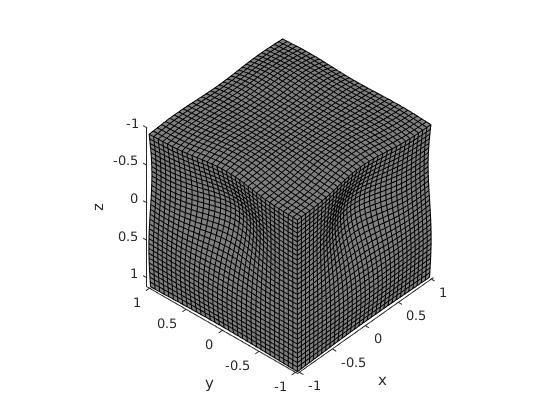}}
   \qquad
      \subfigure[negative tangential contact forces]{
   \includegraphics[trim = 20mm 4mm 30mm 5mm, clip,width=50.0mm, keepaspectratio]{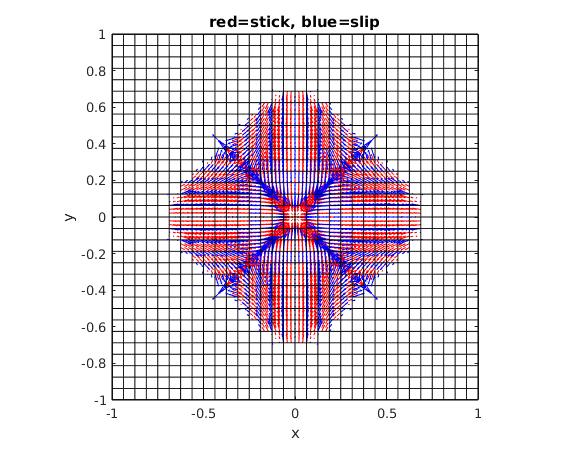}}
   }
   \caption{Discrete solution on uniform mesh with 6144 elements and $p=1$ (3D Coulomb-friction case).}
   \label{fig:3d:solution}
 \end{figure}

Let $\Omega=[-1,1]^3$, $\Gamma_C=[-1,1]^2 \times \{-1\}$, $\Gamma_N = \partial \Omega \setminus \Gamma_C$. The remaining data are $\bar{\gamma}=10^{-3}$, $E=5$, $\nu=0.45$,  $g=0$, $\mathcal{F}=0.3$ and 
\begin{alignat*}{2}
 f&=  \left(
 \begin{array}{c}
 -c_3\operatorname{sign}(x)(1-z^2)(1-y^2)\exp( -c_4( (z+c_5)^2+y^2 ) )\\
 0 \\
 	c_2(1-z^2)(1-y^2)
 \end{array}\right) &\quad \text{ on }& \left\{\pm 1\right\}\times \left[-1,1\right]\times \left[-1,1\right],
 \\
 f&=  \left(
 \begin{array}{c}
 	0\\ 
 	-c_3\operatorname{sign}(y)(1-z^2)(1-x^2)\exp( -c_4( (z+c_5)^2+x^2 ) ) \\
 	c_2(1-z^2)(1-x^2)
 \end{array}\right) &\quad \text{ on }& \left[-1,1\right]\times \left\{\pm 1 \right\} \times \left[-1,1\right].
 \\
 f&=  \left(
 \begin{array}{c}
 	0\\ 0 \\
 	-c_1(1-x^2)(1-y^2)\exp( -1.5c_4( x^2+y^2 ) )
 \end{array}\right) &\quad \text{ on }& \left[-1,1\right]\times\left[-1,1\right]\times \left\{1\right\}
\end{alignat*}
with $c_1=2$, $c_2=0.11$, $c_3=1.5$, $c_4=4$ and $c_5=0.5$.
The six rigid body motions are removed by limiting the movements of the midpoint of the six faces of the cube according to the symmetry of the domain and of the applied forces. The deformed body with its diamond shaped actual contact set and the negative tangential contact forces are shown in Figure~\ref{fig:3d:solution}.

\begin{figure}[tbp]
   \centering \mbox{
   \subfigure[a posteriori error estimation]{
 	\includegraphics[trim = 10mm 0mm 17mm 6mm, clip,width=60.0mm, keepaspectratio]{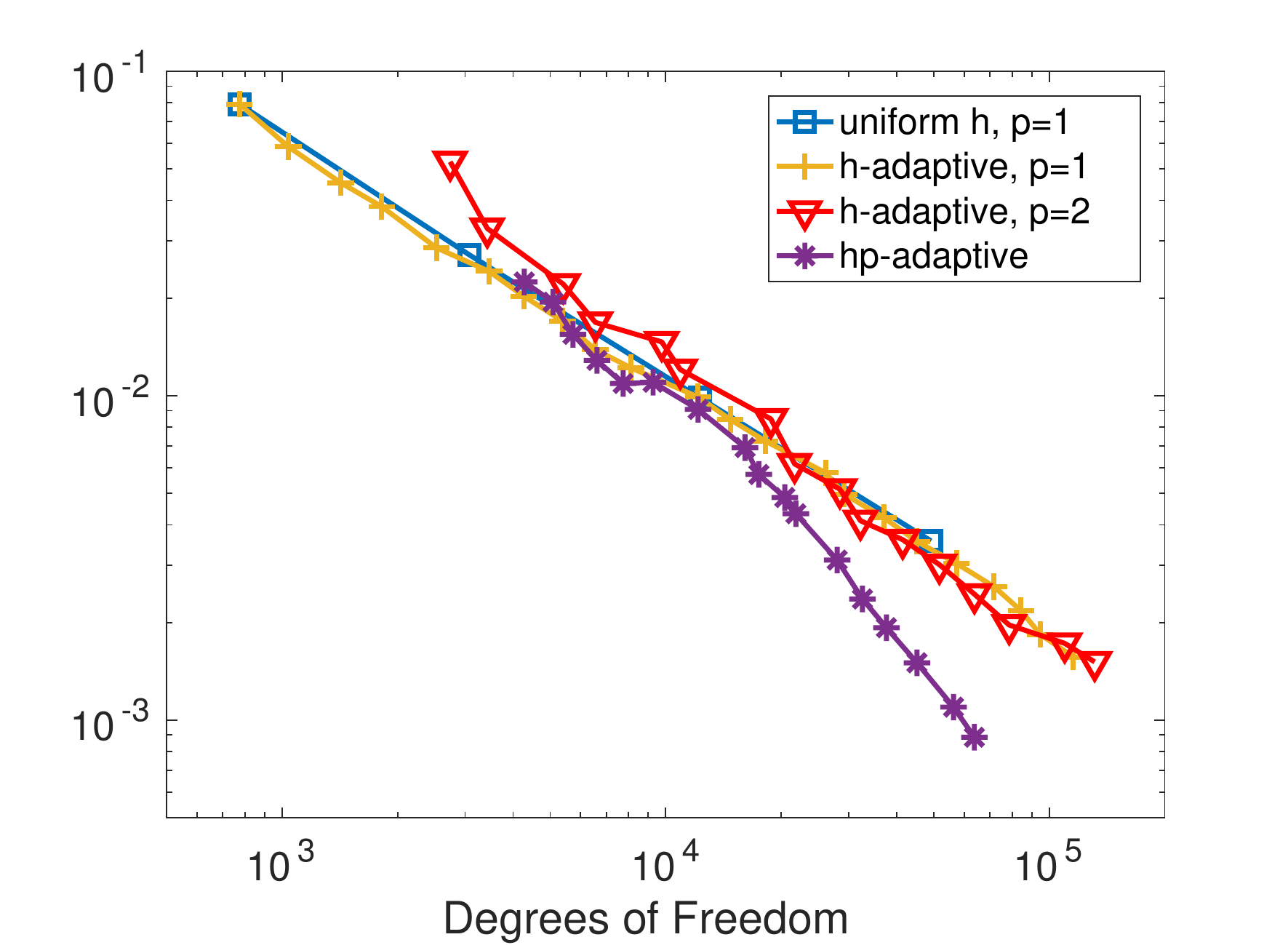}}  \qquad 
   \subfigure[approximate error in energy norm]{
 	\includegraphics[trim = 10mm 0mm 17mm 6mm, clip,width=60.0mm, keepaspectratio]{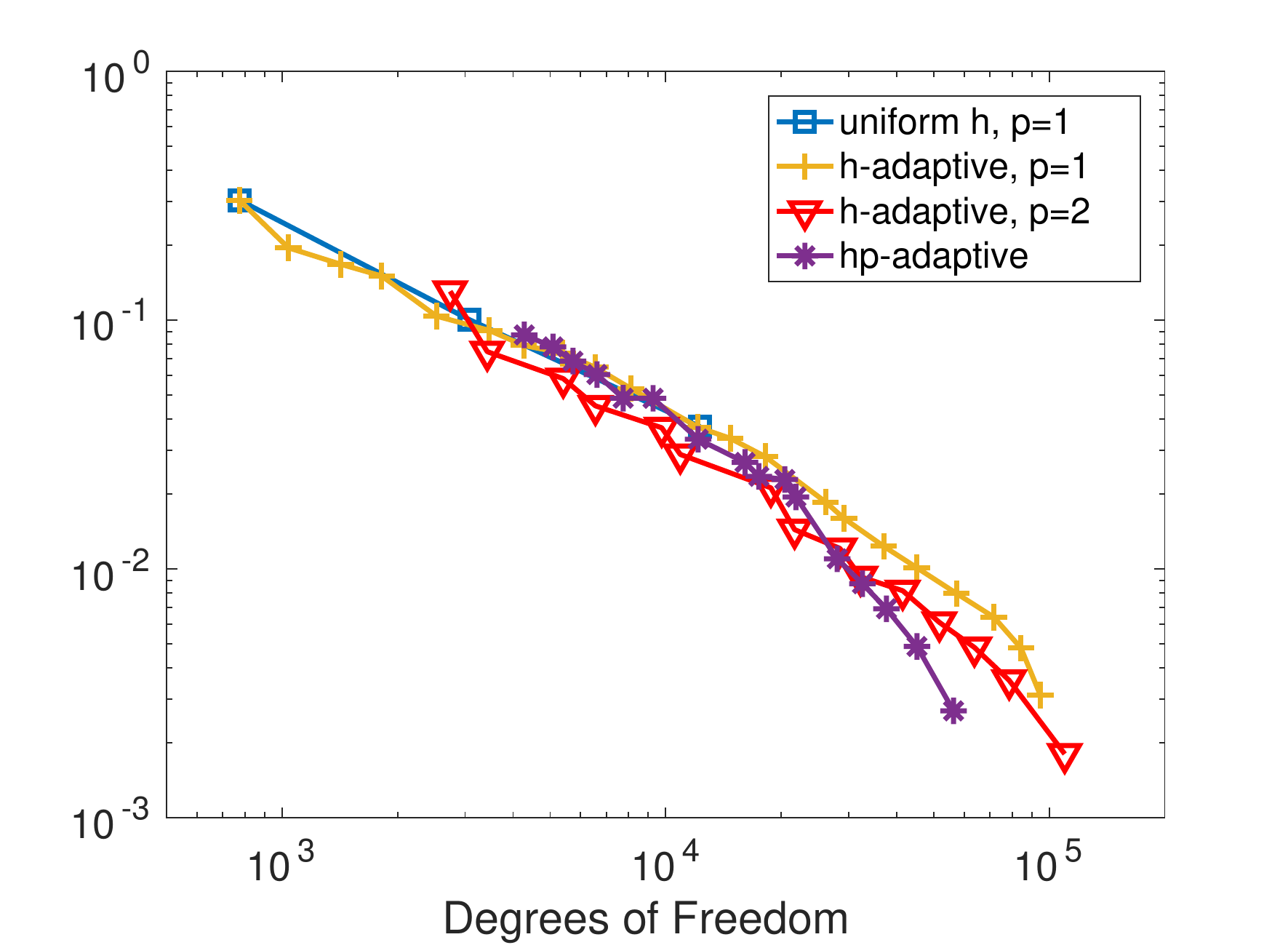}}  }
   \caption{Error estimates for different families of discrete solutions (3D Coulomb-friction case).}
   \label{fig:3d:errorEstimates}
 \end{figure}

The decay of the a posteriori error estimate and of the approximate error is displayed in Figure~\ref{fig:3d:errorEstimates}. The convergence rate of the error and its estimate for the lowest order uniform $h$-version is with $0.75$ w.r.t~the degrees of freedom optimal and thus cannot be improved by the corresponding $h$-adaptive method. The $h$-adaptive method with $p=2$ achieves an experimental convergence rate of $0.9$ due to the inability of resolving the free boundary which goes diagonally though the boundary elements by isotropic refinements. The best convergence rate of $1.46$ is obtained with the $hp$-adaptive scheme. Again the isotropic refinement strategy limits the order of convergence, see e.g.~\cite{banz2015biorthogonal,banz2013posteriori} for the same observation in a 2d obstacle problem with FEM.

\begin{figure}[tbp]
   \centering \mbox{\!
   \subfigure[$h$-adaptive with $p=1$]{
 	\includegraphics[trim = 20mm 2mm 19mm 10, clip,width=53.0mm, keepaspectratio]{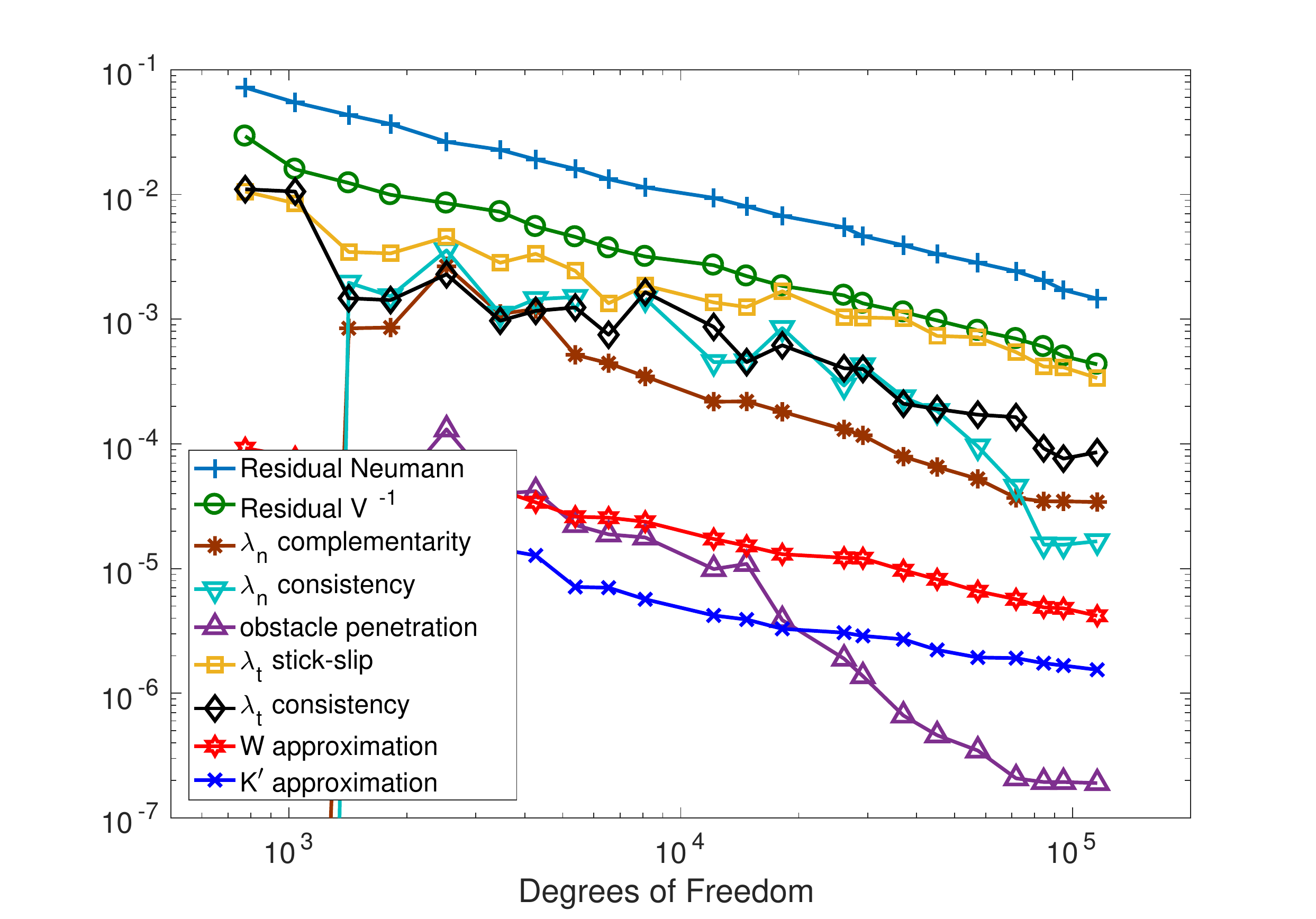}}
   \subfigure[$h$-adaptive with $p=2$]{
 	\includegraphics[trim = 20mm 2mm 19mm 10, clip,width=53.0mm, keepaspectratio]{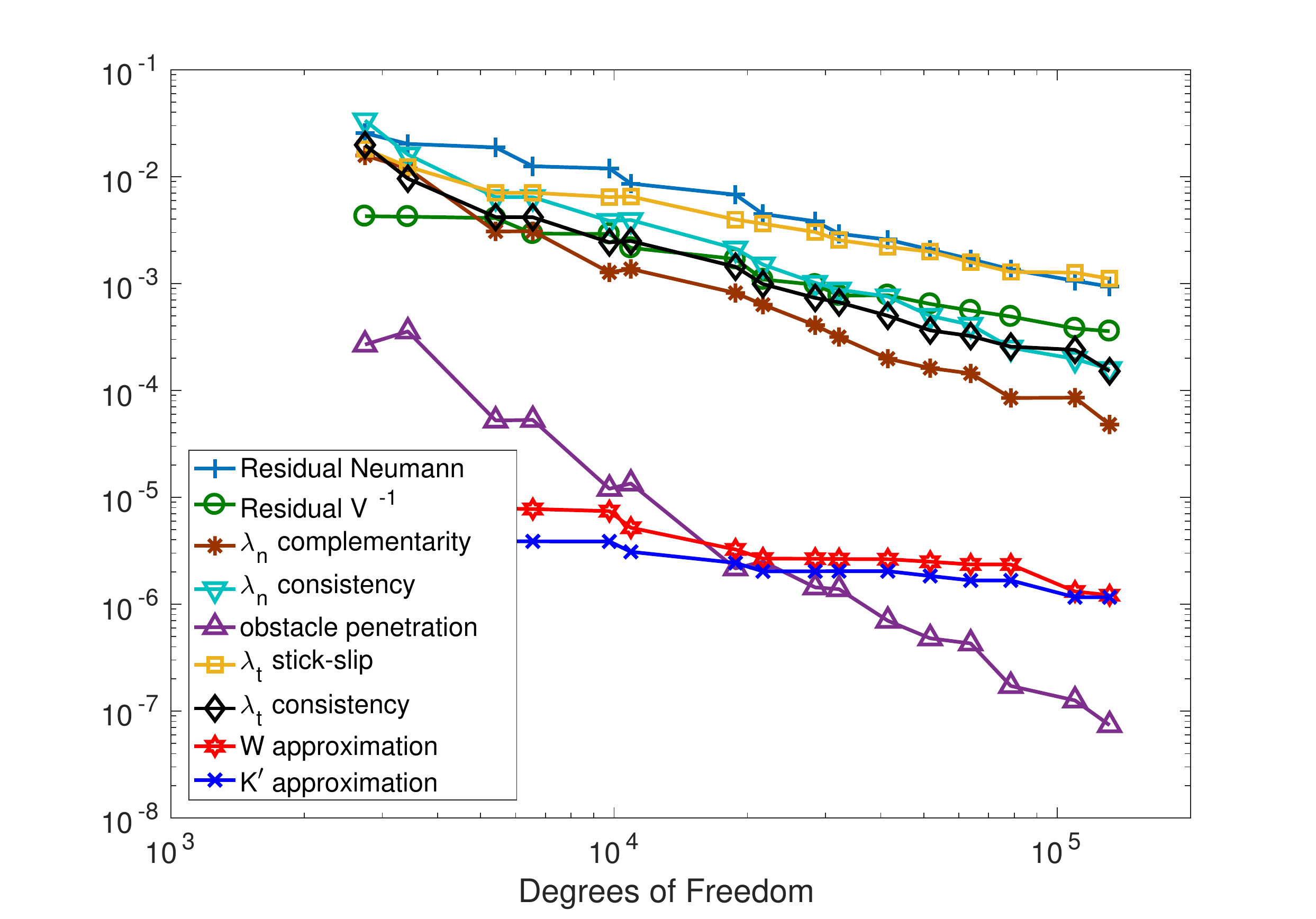}}  	
  \subfigure[$hp$-adaptive]{
 	\includegraphics[trim = 20mm 2mm 19mm 10mm, clip,width=53.0mm, keepaspectratio]{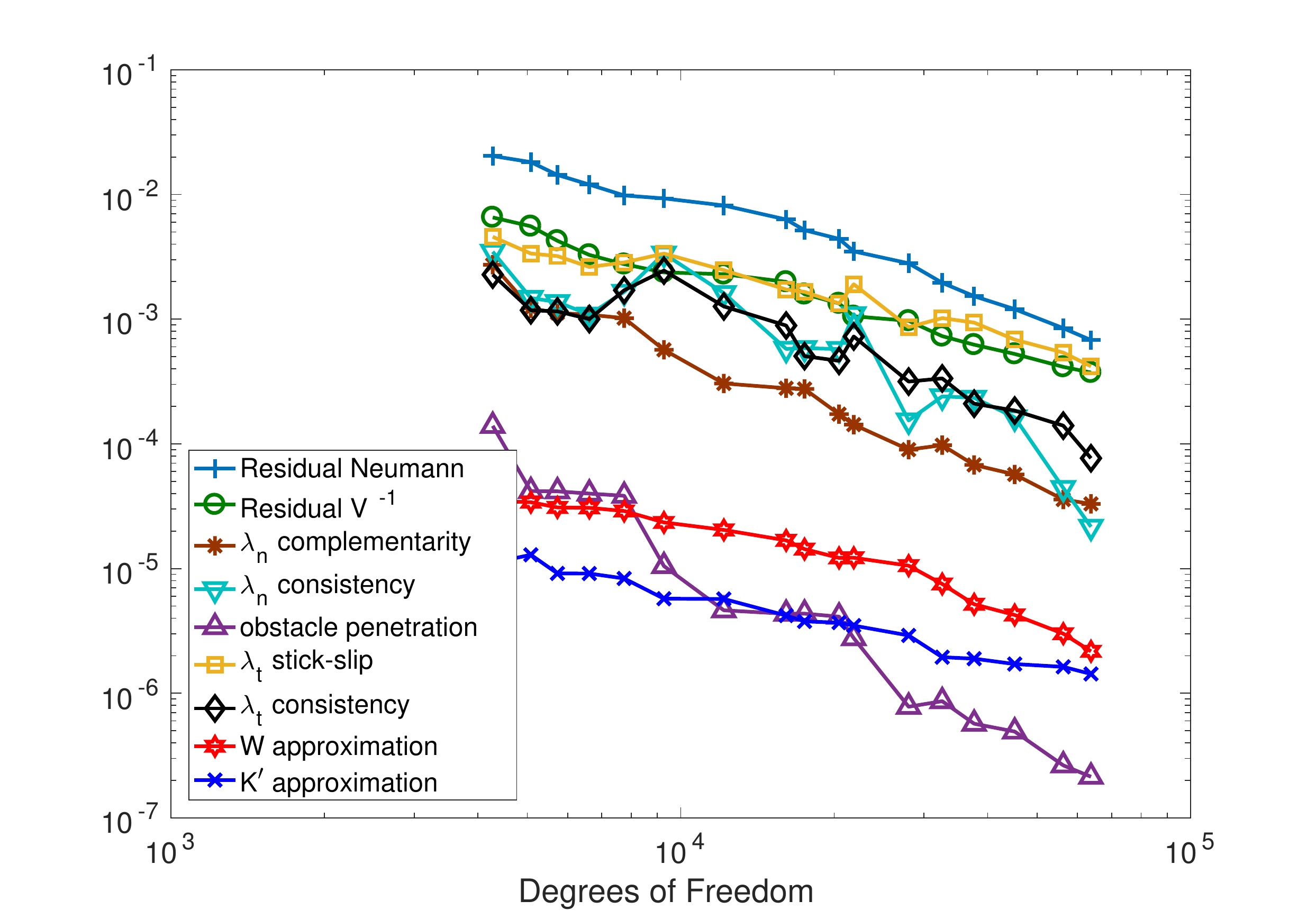}} } 
		
   \caption{Error contributions of the residual based a posteriori error estimate (3D Coulomb-friction case).}
   \label{fig:3d:errorContributions}
 \end{figure}

The individual contributions of the a posteriori error estimate are depicted in Figure~\ref{fig:3d:errorContributions}. Once again the dominant source of error is the classical residual error contribution. Equally important is the violation of the stick-slip condition for the $h$-adaptive scheme with $p=2$ and the $hp$-adaptive scheme.

\begin{figure}[tbp]
   \centering \mbox{\!
   \subfigure[$h$-adaptive with $p=1$]{
 	\includegraphics[trim = 20mm 10mm 21mm 11mm, clip,width=53.0mm, keepaspectratio]{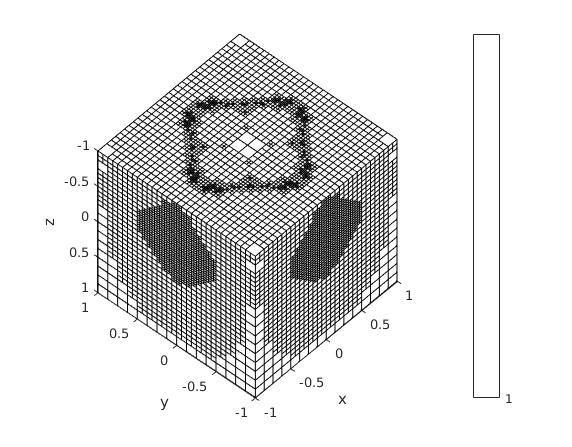}}
   \subfigure[$h$-adaptive with $p=2$]{
 	\includegraphics[trim = 20mm 10mm 21mm 11mm, clip,width=53.0mm, keepaspectratio]{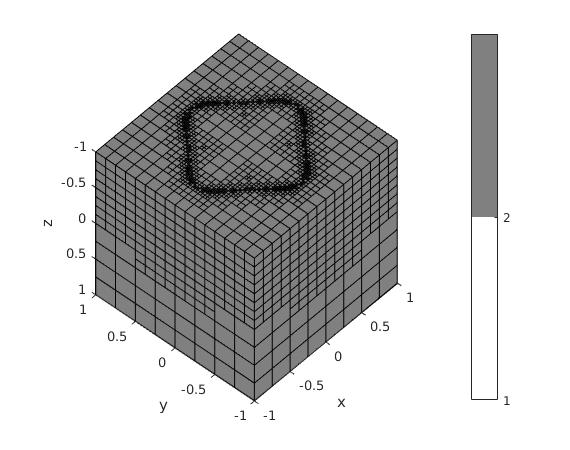}}  	
  \subfigure[$hp$-adaptive]{
 	\includegraphics[trim = 20mm 10mm 21mm 11mm, clip,width=53.0mm, keepaspectratio]{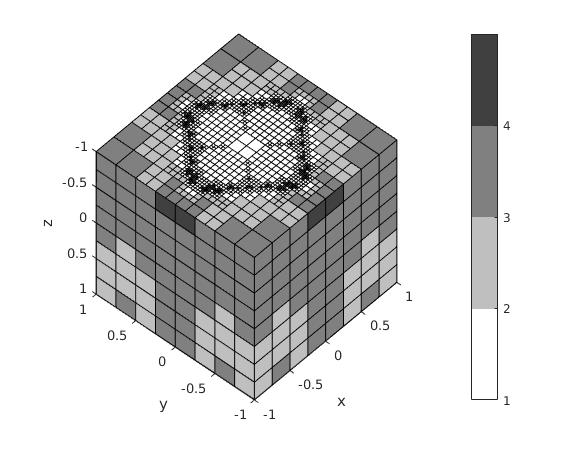}} } 
		
   \caption{Adaptively generated meshes (3D Coulomb-friction case).}
   \label{fig:3d:meshes}
 \end{figure}

All three types of adaptive scheme identify the free boundary and perform mesh refinements towards it, c.f.~Figure~\ref{fig:3d:meshes}. Additionally, the $hp$-meshes show an increase in the polynomial degree away from the free boundary.

%
%
%
%
%
  \bibliographystyle{model1b-num-names}
  \bibliography{lit2}
%
%
%
%
%
%
%

\end{document}